\newcommand{\xBc}{\langle}
\newcommand{\xBe}{\rangle}

\newcommand{\xbF}{\Phi}
\newcommand{\xbG}{\Gamma}

\newcommand{\xbP}{\Pi}

\newcommand{\xba}{\alpha}
\newcommand{\xbb}{\beta}

\newcommand{\xbe}{\in}
\newcommand{\xbf}{\phi}
\newcommand{\xbg}{\gamma}

\newcommand{\xbo}{\omega}

\newcommand{\xbr}{\rho}
\newcommand{\xbs}{\sigma}
\newcommand{\xbt}{\tau}

\newcommand{\xCK}{\times}
\newcommand{\xCL}{\pm}

\newcommand{\xCN}{\neg}
\newcommand{\xCO}{ }
\newcommand{\xCQ}{\emptyset}

\newcommand{\xCe}{>}

\newcommand{\xcA}{\forall}

\newcommand{\xcE}{\exists}

\newcommand{\xcH}{\not\Rightarrow}

\newcommand{\xcO}{\bigvee}
\newcommand{\xcP}{\not\rightarrow}

\newcommand{\xcT}{\bot}
\newcommand{\xcU}{\bigwedge}

\newcommand{\xcg}{\geq}
\newcommand{\xch}{\Rightarrow}

\newcommand{\xck}{\leq}

\newcommand{\xco}{\vee}
\newcommand{\xcp}{\rightarrow}

\newcommand{\xcr}{\leftrightarrow}
\newcommand{\xcs}{\cap}
\newcommand{\xct}{\top}
\newcommand{\xcu}{\wedge}
\newcommand{\xcv}{\cup}

\newcommand{\xcz}{\Box}

\newcommand{\xDH}{\item }

\newcommand{\xDL}{\oplus}

\newcommand{\xDN}{\ominus}
\newcommand{\xDO}{\circ}

\newcommand{\xDc}{\ll}

\newcommand{\xEI}{\begin{itemize}}
\newcommand{\xEJ}{\end{itemize}}

\newcommand{\xEd}{\neq}

\newcommand{\xEh}{\begin{enumerate}}

\newcommand{\xEj}{\end{enumerate}}

\newcommand{\xEn}{\begin{description}}
\newcommand{\xEp}{\end{description}}

\newcommand{\xeA}{\nabla}

\newcommand{\xex}{\upharpoonright}

\newcommand{\Xl}{\ldots}

\newcommand{\bl}{\begin{lemma} \rm}
\newcommand{\el}{\end{lemma}}
\newcommand{\br}{\begin{remark} \rm}
\newcommand{\er}{\end{remark}}
\newcommand{\be}{\begin{example} \rm}
\newcommand{\ee}{\end{example}}
\newcommand{\bco}{\begin{corollary} \rm}
\newcommand{\eco}{\end{corollary}}
\newcommand{\bc}{\begin{claim} \rm}
\newcommand{\ec}{\end{claim}}
\newcommand{\bfa}{\begin{fact} \rm}
\newcommand{\efa}{\end{fact}}
\newcommand{\bp}{\begin{proposition} \rm}
\newcommand{\ep}{\end{proposition}}
\newcommand{\bd}{\begin{definition} \rm}
\newcommand{\ed}{\end{definition}}
\newcommand{\bcs}{\begin{construction} \rm}
\newcommand{\ecs}{\end{construction}}
\newcommand{\bcd}{\begin{condition} \rm}
\newcommand{\ecd}{\end{condition}}
\newcommand{\bt}{\begin{theorem} \rm}
\newcommand{\et}{\end{theorem}}
\newcommand{\bn}{\begin{notation} \rm}
\newcommand{\en}{\end{notation}}
\newcommand{\bfi}{\begin{bild} \rm}
\newcommand{\efi}{\end{bild}}
\newcommand{\bsta}{\begin{statement} \rm}
\newcommand{\esta}{\end{statement}}
\newcommand{\bcom}{\begin{comment} \rm}
\newcommand{\ecom}{\end{comment}}
\newcommand{\bdia}{\begin{diagram} \rm}
\newcommand{\edia}{\end{diagram}}

\newcommand{\bfc}{\begin{figure}[htb] \begin{center}}
\newcommand{\efc}{\end{center} \end{figure}}

\documentclass{article}

\usepackage{amssymb,latexsym,epic,eepic,rotating}

\title{Comments on Yablo's Construction
\thanks{File: XAB, [Sch24]
}
}

\author{Karl Schlechta
\thanks{
schcsg@gmail.com - https://sites.google.com/site/schlechtakarl/ -
Koppeweg 24, D-97833 Frammersbach, Germany}
\thanks{
Retired, formerly: Aix-Marseille Universit\'{e}, CNRS, LIF UMR 7279, F-13000
Marseille, France
}
}

\begin{document}

\newtheorem{lemma}{Lemma}[section]
\newtheorem{theorem}[lemma]{Theorem}
\newtheorem{proposition}[lemma]{Proposition}
\newtheorem{corollary}[lemma]{Corollary}
\newtheorem{claim}[lemma]{Claim}
\newtheorem{fact}[lemma]{Fact}
\newtheorem{remark}[lemma]{Remark}
\newtheorem{definition}{Definition}[section]
\newtheorem{construction}{Construction}[section]
\newtheorem{condition}{Condition}[section]
\newtheorem{example}{Example}[section]
\newtheorem{notation}{Notation}[section]
\newtheorem{bild}{Figure}[section]
\newtheorem{comment}{Comment}[section]
\newtheorem{statement}{Statement}[section]
\newtheorem{diagram}{Diagram}[section]

\renewcommand{\labelenumi}
  {(\arabic{enumi})}
\renewcommand{\labelenumii}
  {(\arabic{enumi}.\arabic{enumii})}
\renewcommand{\labelenumiii}
  {(\arabic{enumi}.\arabic{enumii}.\arabic{enumiii})}
\renewcommand{\labelenumiv}
  {(\arabic{enumi}.\arabic{enumii}.\arabic{enumiii}.\arabic{enumiv})}

\maketitle

\setcounter{secnumdepth}{3}
\setcounter{tocdepth}{3}

\begin{abstract}

We analyse Yablo's coding of the liar paradox by infinite acyclic graphs,
\cite{Yab82}, and show that his construction is in a certain way minimal.

This leads to a very limited representation result.

We then argue that the "right" level of description of more complicated
constructions is probably the level of paths, and not of single arrows.

We also look at more complicated cells of basic contradiction ("diamonds"), and
show that, at least under certain conditions, they fail to achieve
Yablo's result.

\end{abstract}

\tableofcontents
\clearpage

%
%
%

$ \xCO $
\clearpage
\section{
Introduction
}

\label{Section IntroductionGlobal}
\subsection{
Overview
}

Yablo's idea (upon which all of the present material is based),
appeared in  \cite{Yab82}.

The author's interest was kindled by  \cite{RRM13},
commented in  \cite{Sch18c}, and extended in  \cite{Sch22}.
(It was shwon there that conjecture 15 in  \cite{RRM13} is wrong.)

We liberally use material from  \cite{Sch22}, in particular in
Section 
\ref{Section IntroductionGlobal} (page 
\pageref{Section IntroductionGlobal}),
Section \ref{Section Or-And} (page \pageref{Section Or-And}),
and Section 
\ref{Section Saw-Blades-All} (page 
\pageref{Section Saw-Blades-All}).

Our approach is ``bottom-up''. We analyse Yablo's construction, an
intricate composition of contradictions in ``Yablo Cells'', of the
graph form $x \xcP y,$ $y \xcP z,$ $x \xcP z,$ meaning $x= \xCN y \xcu
\xCN z,$ $y= \xCN z.$ We distinguish
$y$ and $z,$ calling $y$ the knee of the cell, $z$ the foot (and $x$ the
head).
In Yablo's construction, $y$ and $z$ will be the head of new cells (among
other
arrows and constructions), etc.

Thus, $x$ is contradictory if true (as a formula), and if $x$ is false,
all
arrows starting at $x$ will lead to new contradictions (as then $y$ and
$z$ will
be true).
We take this as construction principle, and show that the property is
recursively true, i.e., will hold for $y$ and $z,$ etc., too.
We can base Yablo's and similar constructions on this principle.
In particular, Yablo's construction is minimal to satisfy these
principles.

It is natural to look now at alternative contradictory cells as building
blocks. We discuss several candidates, the only interesting one is the
``diamond'' of the form $x \xcP y,$ $y \xcP z,$ $x \xcP y',$ $y \xcp z,$
with the meaning
$x= \xCN y \xcu \xCN y',$ $y= \xCN z,$ $y' =z.$
But the recursive construction, appending
new diamonds at $y$ and $y' $ fails, at least under the restriction of
using
only conjunctions for the logical form in the positive case (x true).

We also
 \xEI
 \xDH argue that the right level of abstraction is to consider paths, not
single arrows,
 \xDH give some examples of Yablo-like graphs with more general formulas,
i.e. of disjunctive normal forms, and different basic orders (i.e. not
necessarily that of the natural numbers),
 \xDH illustrate postponing contradictions.
 \xEJ
\subsubsection{
Related Work
}

The ideas in
 \cite{BS17} and  \cite{Wal23}
are more ``top down'', whereas our ideas are more ``naive'', starting
from an analysis of Yablo's construction and the basic contradictions.
Thus, we needed an argument to consider paths, and not combinations
of arrows as the ``right'' level of abstraction, whereas paths are the
basic level of investigation in above articles.

The closest common points are perhaps in Section 4.5 of
 \cite{BS17}. In addition, our
Section \ref{Section Procrastination} (page \pageref{Section Procrastination})
can be seen as an allusion to a game theoretic interpretation.
\subsection{
More Detailed Overview
}

Some aspects of Yablo's construction
are hidden behind its elegance. Our perhaps somewhat pedantic analysis
leads to the concepts of ``head'', ``knee'', and ``foot'', and then to ``saw
blades''
in Section 
\ref{Section Saw-Blades-All} (page 
\pageref{Section Saw-Blades-All}).

Section \ref{Section Or-And} (page \pageref{Section Or-And})
generalizes Yablo's construction to
arbitrary formulas of the type $ \xcO \xcU $ (disjunctive normal
forms), and offers a number of easy variations of such
structures by modifying the order of the graph.

Section 
\ref{Section Saw-Blades-All} (page 
\pageref{Section Saw-Blades-All})  uses our idea of finer analysis
of contradictory cells to build somewhat different
structures - though the distinction is blurred by the
necessarily recursive and infinite construction of contradictions.
This is used in our limited (semi) representation results in
Section 
\ref{Section Saw-Blades-Discussion} (page 
\pageref{Section Saw-Blades-Discussion})  and
Section \ref{Section Paths-Arrows} (page \pageref{Section Paths-Arrows}).

In Section 
\ref{Section AbstractionLevel} (page 
\pageref{Section AbstractionLevel}),
we show that we can make contradictions arbitrarily complicated when
seen on the level of arrows. This leads us to think that the level of
arrows
is too detailed, it is not a useful level of abstraction.
We therefore work with paths, trivialising some intermediate nodes.

In Section 
\ref{Section ContradictionCells} (page 
\pageref{Section ContradictionCells})
we look at various candidates for basic contradiction cells,
one of them (Diamond) seems much
more complicated, and embedding it into a recursive construction as in
Yablo's work leads to not so obvious problems, unsolvable at least in the
present framework.

In Section 
\ref{Section Procrastination} (page 
\pageref{Section Procrastination}),
we analyse the finite procrastination of contradictions using Yablo's form
of
contradictions.

The author of the present text did not study the literature
systematically. So,
if some examples are already discussed elsewhere, the author would ask to
be
excused for not quoting previous work.
\subsubsection{
General Strategy
}

Our general strategy (as in Yablo's original paper) will be as follows:

 \xEh
 \xDH
We want to show that $x=TRUE$ is impossible, and $x=FALSE$ is impossible,
too.

 \xDH
Say $x$ is the root of a graph, so for every possibility in this graph,
$x=TRUE$ is impossible, as every set of paths from $x$ leads to a
contradiction in
case $x=TRUE.$

 \xDH
Conversely, if $x=FALSE,$ then every path from $x$ leads leads to some $x'
$
such that $x' =TRUE$ and this is impossible for the same reasons as above
for $x=TRUE.$

(Thus, the construction is essentially recursive,
see Section 
\ref{Section Saw-Blades-Discussion} (page 
\pageref{Section Saw-Blades-Discussion})).

 \xDH
Moreoever, we will not postpone contradictions, but e.g. in above case
for $x=TRUE$ we will require the graph to branch at $x$ (and not to be
some kind of preamble, where we first have to go to some $x'' $ where the
branching happens). - This is a minor point.

See Section 
\ref{Section Procrastination} (page 
\pageref{Section Procrastination})  for some exceptions.

 \xDH
Yablo uses conjunctions of negated propositional variables,
$x= \xcU \xCN x_{i},$ this takes
care of above ideas, conjunctions provide the bases for contradictions,
and the negated $x$ makes at least one $x_{i}$ positive, this takes care
of
making $x=FALSE$ impossible.

 \xEj

(Graphically, Yablo has negative arrows from $x$ to the $x_{i},$ we
sometimes
consider negative paths instead.
In Section \ref{Section Or-And} (page \pageref{Section Or-And})
we also consider mixed Or/And formulas, and in
Section 
\ref{Section ContradictionCells} (page 
\pageref{Section ContradictionCells})
we consider more complicated contradictions.)
\subsubsection{
Main Ideas of the Paper
}

 \xEh
 \xDH
Analysis of Yablo's construction in
Section 
\ref{Section Saw-Blades-Discussion} (page 
\pageref{Section Saw-Blades-Discussion}), and the
resulting partial representation result in
Condition \ref{Condition Yablo} (page \pageref{Condition Yablo}).

 \xDH
Discussion of the suitable level of abstraction for the description
of more complicated constructions, see
Section \ref{Section NegationCells} (page \pageref{Section NegationCells})  and
Section \ref{Section Paths-Arrows} (page \pageref{Section Paths-Arrows}).

 \xDH
Discussion of more complicated basic contradiction units
(``diamonds'') and the result that they are unsuitable
(at least under the present, limited, conditions), see
Section 
\ref{Section Diamonds} (page 
\pageref{Section Diamonds}), and in particular
Section \ref{Section Nested-Diamonds} (page \pageref{Section Nested-Diamonds}).

 \xEj
\subsection{
Preliminaries
}

\label{Section Preliminaries-XBA}

We consider the problem from the graph side and the logics side.
\subsubsection{
The Logical Side
}

On the logics side,
we work with propositional formulas, which may, however, be infinite.

We will see that we need infinite formulas (of infinite depth) to
have nodes to which we cannot attribute truth values.
See Fact \ref{Fact Infinite} (page \pageref{Fact Infinite}).

We will work here with disjunctive normal forms, i.e. with formulas of the
type
$x:= \xcO \{ \xcU x_{i}:i \xbe I\},$ where $x_{i}:=\{x_{i,j}:j \xbe
J_{i}\},$ and the $x_{i,j}$ are propositional
variables or negations thereof - most of the time pure conjunctions of
negations of propositional variables.

\bfa

$\hspace{0.01em}$


\label{Fact DNF}

Let $x:= \xcO \{ \xcU \{x_{i,j}:j \xbe J_{i}\}:i \xbe I\},$ where the
$x_{i,j}$ are propositional
variables or negations thereof.

 \xEh

 \xDH

Let $F:= \xbP \{x_{i}:i \xbe I\}$ - the set of choice functions in the
sets $x_{i},$ where
$x_{i}:=\{x_{i,j}:j \xbe J_{i}\}.$

Then $ \xCN x$ $=$ $ \xcO \{ \xcU \{ \xCN x_{i,j}:$ $x_{i,j} \xbe
ran(f)\}:f \xbe F\}.$

(Argue semantically with the sets of models and general distributivity
and complements of sets.)

 \xDH

Contradictions will be between two formulas only, one a propositional
variable, the other the negation of the former.

 \xEj

$ \xcz $
\\[3ex]
\subsubsection{
The Graph Side
}

\efa

We work with directed, acyclic graphs. They will usually have one root,
often called $x_{0}.$ In diagrams, the graphs may grow upward, downward,
or
sideways, we will say what is meant.

\bd

$\hspace{0.01em}$


\label{Definition Graph}

Nodes stand for propositional variables.

If a node $x$ is not terminal, it has also a propositional formula $
\xbf_{x}$
attached to it, sometimes written $d(x)= \xbf_{x},$
with the meaning $x \xcr \xbf_{x},$ abbreviated $x= \xbf_{x}.$
The successors of $x$ are the propositional variables occuring in $
\xbf_{x}.$
Thus, if $x \xcp x' $ and $x \xcp x'' $ are the only successors of $x$ in
$ \xbg,$ $ \xbf_{x}$ may be
$x' \xco x'',$ $x' \xcu \xCN x'',$ but not $x' \xcu y.$

Usually, the $ \xbf_{x}$ are (possibly infinite) conjunctions of
propositional
variables or (in most cases) their negations, which we write $ \xcU \xCL
x_{i}$ etc.
We often indicate the negated variables in the graph with negated arrows,
like $x \xcP y,$ etc. Thus, $x \xcP x',$ $x \xcp x'' $ usually stands for
$ \xbf_{x}= \xCN x' \xcu x''.$

\ed

\be

$\hspace{0.01em}$


\label{Example Basic}

Consider the basic construction of a contradiction (used by Yablo and
here,
defined ``officially'' in
Definition \ref{Definition YC} (page \pageref{Definition YC})).

$ \xbG:=\{x \xcP y \xcP z,$ $x \xcP z\}.$ $ \xbG $ stands for $x= \xCN y
\xcu \xCN z,$ $y= \xCN z,$ so $x=z \xcu \xCN z,$ which
is impossible.

If we negate $x,$ then $ \xCN x=y \xco z= \xCN z \xco z,$ so $ \xCN x$ is
possible.

From the graph perspective, we have two paths in $ \xbG $ from $x$ to $z,$
$ \xbs:=x \xcP y \xcP z,$
and $ \xbs':=x \xcP z.$

We add now $y \xcP y' $ to $ \xbG,$ so $ \xbG':=\{x \xcP y \xcP z,$ $x
\xcP z,$ $y \xcP y' \},$ thus
$x= \xCN y \xcu \xCN z,$ $y= \xCN z \xcu \xCN y',$ so $ \xCN y=z \xco y'
,$ and $x=(z \xco y') \xcu \xCN z=(z \xcu \xCN z) \xco (y' \xcu \xCN z),$
and $x$ is not contradictory any more.

\ee

\bd

$\hspace{0.01em}$


\label{Definition Path-Value}

We can attribute a value to a path $ \xbs,$ $val(\xbs),$ expressing
whether it
changes a truth value from the beginning to the end. $ \xbs:=x \xcP y
\xcP z$ does not change
the value of $z$ compared to that of $x,$ $ \xbs':=x \xcP z$ does. We
say $val(\xbs)=+,$
$val(\xbs')=-,$ or positive (negative) path.

Formally:

Let $ \xbs,$ $ \xbs' $ be paths as usual.
 \xEh
 \xDH
If $ \xbs:=a \xcp b,$ then $val(\xbs)=+,$ if $ \xbs:=a \xcP b,$ then
$val(\xbs)=-.$
 \xDH
Let $ \xbs \xDO \xbs' $ be the concatenation of $ \xbs $ and $ \xbs'.$
Then $val(\xbs \xDO \xbs')=+$ iff $val(\xbs)=val(\xbs'),$ and -
otherwise.
 \xEj

\ed

If all arrows are negative, then $val(\xbs)=+$ iff the length of $ \xbs
$ is even.

\bd

$\hspace{0.01em}$


\label{Definition Contradictory-Paths}

We call two paths $ \xbs,$ $ \xbs' $ with common start and end
contradictory, and the pair a contradictory cell iff $val(\xbs) \xEd
val(\xbs').$ The
structures considered here
will be built with contradictory cells.

\ed

\br

$\hspace{0.01em}$


\label{Remark Contradictory-Paths}

 \xEh

 \xDH
Note that the fact that $ \xbs $ and $ \xbs' $ are contradictory or not
is
independent of how we start, whether $x=TRUE$ or $x=FALSE.$

 \xDH
We saw already in Example 
\ref{Example Basic} (page 
\pageref{Example Basic})  that it is not sufficient
for a ``real'' contradiction to have two contradictory paths.

We need
 \xEh
 \xDH (somewhere) an ``AND'', so both have to be valid together, an ``OR''
is not sufficient,
 \xDH we must not have a branching with an ``OR'' on the way as in $ \xbG
',$ an
``escape'' path, unless this leads again to a contradiction.
 \xEj
See also Definition 
\ref{Definition Terminology} (page 
\pageref{Definition Terminology}).

 \xEj

\er

\bn

$\hspace{0.01em}$


\label{Notation Plus-And}

 \xEh
 \xDH
When we give nodes a truth value, we will use $x+$ (and $x \xcU,$ $x+
\xcU,$ etc.
if $ \xbf_{x}$ has the form $ \xcU \xCL x_{i})$
to denote the case $x=TRUE,$ and $x-,$ $x \xcO,$ $x- \xcO,$ etc. for the
case $x=FALSE.$

 \xDH
We sometimes use the notation $x$ $ \xDL $ for shorthand that $x+$ is
contradictory,
and $x$ $ \xDN $ for shorthand that every path (mostly: every arrow) from
$x$ will lead to some $y$ such that $y$ $ \xDL,$ and $y+$ if $x-.$
 \xEj

\en

\bfa

$\hspace{0.01em}$


\label{Fact Value2}

(Simplified).

Consider three paths, $ \xbr,$ $ \xbs,$ $ \xbt,$ for simplicity with
same origin, i.e.
$ \xbr (0)$ $=$ $ \xbs (0)$ $=$ $ \xbt (0).$

 \xEh

 \xDH No contradiction loops of length 3.

 \xEh

 \xDH
Suppose they meet at a common point, i.e. $ \xbr (m_{ \xbr })$ $=$ $ \xbs
(m_{ \xbs })$ $=$ $ \xbt (m_{ \xbt }).$
Then it is impossible that
$ \xbr $ contradicts $ \xbs $ contradicts $ \xbt $ contradicts $ \xbr $
(at $m_{ \xbr }).$
(``$ \xba $ contradicts $ \xbb $'' means here that for some $i,j$ $ \xba
(i)= \xbb (j),$ but
one has value $+,$ the other value -.)

(Trivial, we have only 2 truth values).
 \xDH
Suppose, first $ \xbr $ and $ \xbs $ meet, then $ \xbr $ (or $ \xbs)$ and
$ \xbt $ meet, but
once they meet, they will continue the same way (e.g., if $ \xbr (i)= \xbs
(j),$
then for all $k>0$ $ \xbr (i+k)= \xbs (j+k)$).
Then it is again impossible that
$ \xbr $ contradicts $ \xbs $ contradicts $ \xbt $ contradicts $ \xbr.$
$(\xbr $ and $ \xbs $ continue to be the same but with different truth
values
until they meet $ \xbt,$ so it is the same situation as above.)
 \xEj

 \xDH
Above properties generalize to any loops of odd length (with more paths).
 \xEj

See  \cite{Sch22} for more details.

\efa

This does not hold when the paths may branch again after meeting, as
the next Example shows.

\be

$\hspace{0.01em}$


\label{Example Tower2}

(Example 7.4.2 in  \cite{Sch22}.)

Let $ \xbs_{0}:x_{0} \xcP x_{1} \xcp x_{2} \xcP x_{3} \xcp x_{4},$
$ \xbs_{1}:x_{0} \xcP x_{1} \xcp x_{2} \xcp x_{4},$
$ \xbs_{2}:x_{0} \xcp x_{2} \xcP x_{3} \xcp x_{4},$
$ \xbs_{3}:x_{0} \xcp x_{2} \xcp x_{4},$

then $ \xbs_{0}$ contradicts $ \xbs_{1}$ in the lower part, $ \xbs_{2}$
and $ \xbs_{3}$ in the upper part,
$ \xbs_{1}$ contradicts $ \xbs_{2}$ and $ \xbs_{3}$ in the upper part,
$ \xbs_{2}$ contradicts $ \xbs_{3}$ in the lower part.

Obviously, this may be generalized to $2^{ \xbo }$ paths.

\ee

Consider Yablo's original construction:

\be

$\hspace{0.01em}$


\label{Example Yablo}

Let the nodes be $\{x_{i}:i< \xbo \},$ and the arrows for $x_{i}$ $\{x_{i}
\xcP x_{j}:i<j\},$ expressed
as a relation by $\{x_{i}<x_{j}:i<j\},$ and as a logical formula by
$x_{i}= \xcU \{ \xCN x_{j}:i<j\}.$

Thus $ \xCN x_{i}= \xcO \{x_{j}:i<j\}.$
For any $x_{i},$ we have a contradiction by $x_{i}= \xcU \{ \xCN
x_{j}:i<j\}$ and $ \xCN x_{i+1}= \xcO \{x_{j}:i+1<j\}$
for any $x_{i}+,$ and for any $x_{k}-$ for a suitable $k' >k$ $x_{k' }+.$

It is important that, although we needed to show the property $ \xDN $ for
$x_{0}$
only, it holds for all $x_{i},$ thus it is a recursive construction.
\subsubsection{
Interplay of the Graph and the Logical Side
}

\ee

We can either think on the logical level with formulas, or graphically
with conflicting paths, as in the following Fact.

We need infinite depth and width in our constructions:

\bfa

$\hspace{0.01em}$


\label{Fact Infinite}

 \xEh
 \xDH The construction needs infinite depth,

 \xDH the logic as used in Yablo's construction is not compact,

 \xDH the construction needs infinite width, i.e. the logic cannot be
classical.
 \xEj

\efa

\subparagraph{
Proof
}

$\hspace{0.01em}$


 \xEh
 \xDH
Let $x_{i}$ be a minimal element, then we can chose an arbitrary truth
value
$x_{i},$ and propagate this value upwards. If there are no infinite
descending
chains, we can do this for the whole construction.

 \xDH
The logic as used in Yablo's construction is not compact:

Trivial.
(Take $\{ \xcO \{ \xbf_{i}:i \xbe \xbo \}\} \xcv \{ \xCN \xbf_{i}:i \xbe
\xbo \}.$ This is obviously inconsistent,
but no finite subset is.).

 \xDH
It is impossible to construct a Yablo-like structure with classical logic:

Take an acyclic graph, and interpret it as in Yablo's construction.
Wlog., we may assume the graph is connected. Suppose it shows that
$x_{0}$ cannot be given a truth value. Then the set of formulas showing
this does not have a model, so it is inconsistent. If the formulas were
classical, it would have a finite, inconsistent subset, $ \xbF.$ Define
the depth
of a formula as the shortest path from $x_{0}$ to this formula.
There is a (finite) $n$ such that all formulas in $ \xbF $ have depth $
\xck n.$
Give all formulas of depth $n$ (arbitrary) truth values, and work
upwards using truth functions. As the graph is acyclic, this is possible.
Finally, $x_{0}$ has a truth value.

Thus, we need the infinite $ \xcU / \xcO.$
 \xEj

$ \xcz $
\\[3ex]
\subsubsection{
Further Remarks on Contradictory Structures
}

\label{Section Remarks}

See Section 7.3 in  \cite{Sch22} for more material.

\be

$\hspace{0.01em}$


\label{Example Simple-Cells}

We consider now some simple, contradictory cells. They should not only
be contradictory for the case $x+,$ but also be a potential start for
the case $x-,$ without using more complex cells. (Otherwise, we
postpone the solution, and may forget the overly simple start.)

We will consider in
Section 
\ref{Section ContradictionCells} (page 
\pageref{Section ContradictionCells})
a more complicated contradiction cell (``Diamond''), and discuss here
only simple cells.

For this, we order the complexity of the cases by $(1)<(2.1)<(2.3)$ below,
(2.2) is not contradictory, so it is excluded.

See Diagram \ref{Diagram YC} (page \pageref{Diagram YC}), center part.

 \xEh
 \xDH The cell with 2 arrows.

It corresponds to the formula
$x=y \xcu \xCN y,$ graphically, it has a positive and a negative arrow
from
$x$ to $y,$ so exactly one of $ \xba $ and $ \xbb $ is negative.

If $x$ is positive, we have a contradiction.

If $x$ is negative, however, we have a problem. Then, we have
$x=y \xco \xCN y.$ Let $ \xba $ be the originally positive path, $ \xbb $
the
originally negative path. Note that $ \xba $ is now negative, and $ \xbb $
is positive.
The $ \xba $ presents no problem, as $y$ is positive,
and we can append the same construction to $y,$ and have a contradiction.
However, $ \xbb $ has to lead to a contradiction, too, and, as we will not
use more complicated cells, we face the same problem again, $y$ is
negative. So we have an ``escape path'', assigning $ \xcT $ to all elements
in one
branch.

(Consider
$x_{0} \xch_{ \xCN }x_{1} \xch_{ \xCN }x_{2} \xch_{ \xCN }x_{3}$  \Xl,
setting all $x_{i}:=-$ is a consistent valuation.
So combining this cell with itself does not result in a contradictory
structure.)

Of course, appending at $y$ a Yablo Cell (see below, case (2.3))
may be the beginning of a contradictory structure, but this is
``cheating'', we use a more complex cell.

 \xDH Cells with 3 arrows.

Note that the following examples are not distinguished in the graph!

Again, we want a contradiction for $x$ positive, so we need an $ \xcu $ at
$x,$
and have the possibilities
(up to equivalence) $x=x' \xcu y,$ $x=x' \xcu \xCN y,$ $x= \xCN x' \xcu
y,$ $x= \xCN x' \xcu \xCN y.$

Again, if $x$ is negative, all paths $ \xba:x \xcp y,$ $ \xbb:x \xcp x'
$ (or $ \xbb \xbg:x \xcp x' \xcp y)$
have to lead to a contradiction.

 \xEh
 \xDH one negative arrow, with $x= \xCL x' \xcu \xCL y$
 \xEh
 \xDH
$x \xcp x' \xcP y,$ $x \xcp y,$ corresponding to $x=x' \xcu y,$ so if $x$
is negative, this
is $ \xCN x' \xco \xCN y.$

But, at $y,$ this possible paths ends, and we have the same situation
again,
with a negative start, as in Case (1).

 \xDH
$x \xcp x' \xcp y,$ $x \xcP y$

Here, we have again a positive path to $y,$ through $x',$ so for $x$
negative
both $x' $ and $y$
will be negative, and neither gives a start for a new contradiction.
 \xDH
$x \xcP x' \xcp y,$ $x \xcp y$

This case is analogous to case (2.1.1).
 \xEj

(Similar arguments apply to more complicated cells with an even number
of negative arrows until the first branching point - see
Remark \ref{Remark YC} (page \pageref{Remark YC})  below.

 \xDH 2 negative arrows: not contradictory.

 \xDH The original type of contradiction in Yablo's construction

$x \xcP x' \xcP y,$ $x \xcP y.$

This will be discussed in detail in the rest of the paper.
But we see already that both paths, $x \xcP x' $ and $x \xcP y$ change
sign,
so $x' $ and $y$ will be positive, appending the same type of cell
at $x' $ and $y$ solves the problem (locally), and offers no escape.

 \xEj
 \xEj

\ee

\bd

$\hspace{0.01em}$


\label{Definition YC}

 \xEh
 \xDH
We call the contradiction of the type (2.3) of
Example \ref{Example Simple-Cells} (page \pageref{Example Simple-Cells}), i.e.
$x \xcP x' \xcP y,$ $x \xcP y,$
a Yablo Cell, YC,
and sometimes $x$ its head, $y$ its foot, and $x' $ its knee.

 \xDH
We sometimes abbreviate a Yablo Cells simply by $ \xeA,$ without going
into any further details.

 \xEj

If we combine Yablo Cells, the knee for one cell may become the head
for another Cell, etc.

See Diagram \ref{Diagram YC} (page \pageref{Diagram YC}), upper part.

\clearpage

\begin{diagram}

Simple Contradictions

\label{Diagram YC}
\index{Diagram YC}

\unitlength0.6mm
\begin{picture}(230,310)(0,0)

\put(0,0){\line(1,0){230}}
\put(0,0){\line(0,1){310}}
\put(230,310){\line(-1,0){230}}
\put(230,310){\line(0,-1){310}}

\put(5,300){Horizontal lines indicate arrows pointing to the right,}
\put(5,292){the other lines downward pointing arrows}

\put(5,280){{\rm\bf Yablo Cells }}

\put(15,265){Head}
\put(39,265){Knee}
\put(15,235){Foot}


\put(20,260){$x$}
\put(20,240){$y$}
\put(26,261){\line(1,0){12}}
\put(32,260){\line(0,1){2}}
\put(22,243){\line(0,1){14}}
\put(21,250){\line(1,0){2}}
\put(24,243){\line(1,1){16}}
\put(31,251.5){\line(1,-1){1.5}}


\put(40,260){$x'$}


\put(75,265){Head}
\put(93,271){Knee=}
\put(93,265){Head'}
\put(75,235){Foot}
\put(120,265){Knee'}
\put(95,235){Foot'}


\put(80,260){$x$}
\put(80,240){$y$}
\put(86,261){\line(1,0){12}}
\put(92,260){\line(0,1){2}}
\put(82,243){\line(0,1){14}}
\put(81,250){\line(1,0){2}}
\put(84,243){\line(1,1){16}}
\put(91,251.5){\line(1,-1){1.5}}


\put(100,260){$x'$}
\put(100,240){$y'$}
\put(106,261){\line(1,0){12}}
\put(112,260){\line(0,1){2}}
\put(102,243){\line(0,1){14}}
\put(101,250){\line(1,0){2}}
\put(104,243){\line(1,1){16}}
\put(111,251.5){\line(1,-1){1.5}}


\put(120,260){$x''$}


\put(5,180){{\rm\bf Contradictory Cells }}


\put(20,160){$x$}
\put(20,140){$y$}
\put(19,143){\line(0,1){14}}
\put(22,143){\line(0,1){14}}

\put(15,150){$\xba$}
\put(24,150){$\xbb$}

\put(50,160){$x$}
\put(50,140){$y$}
\put(56,161){\line(1,0){12}}
\put(51,143){\line(0,1){14}}
\put(54,143){\line(1,1){16}}
\put(70,160){$x'$}

\put(47,150){$\xba$}
\put(62,164){$\xbb$}
\put(64,150){$\xbg$}


\put(5,80){{\rm\bf See Remark $\ref{Remark YC}$ }}

\put(5,60){Head}
\put(42,54){Knee}
\put(15,35){Foot}


\put(20,60){$x$}
\put(20,40){$y$}
\put(26,61){\line(1,0){12}}
\put(32,60){\line(0,1){2}}
\put(22,43){\line(0,1){14}}
\put(21,50){\line(1,0){2}}
\put(24,43){\line(1,1){16}}
\put(31,51.5){\line(1,-1){1.5}}


\put(40,60){$x'$}
\put(60,60){$z$}

\put(46,61){\line(1,0){12}}
\put(52,60){\line(0,1){2}}

\put(22,63){\line(0,1){6}}
\put(62,63){\line(0,1){6}}
\put(22,69){\line(1,0){40}}
\put(42,68){\line(0,1){2}}

\end{picture}

\end{diagram}

\clearpage

\ed

\br

$\hspace{0.01em}$


\label{Remark YC}

 \xEh
 \xDH
The distinction between $x' $ and $y,$ i.e. between knee and foot, is
very important. In the case $x+,$ we have at $y$ a complete contradiction,
at $x',$ we have not yet constructed a contradiction. Thus, if we have
at $x' $ again an $ \xcu $ (as at $x),$ this becomes an $ \xco,$ and we
have to
construct a contradiction for all $x' \xcp z$ (or $x' \xcP z),$ not only
for
$x' \xcP y.$ Otherwise, we have an escape possibility for $x+.$ Obviously,
the contradiction need not be immediate at $z,$ the important property
is that ALL paths through all $z$ lead to a contradiction, and the
simplest
way is to have the contradiction immediately at $z$ - as in Yablo's
construction, and our saw blades.

See Diagram \ref{Diagram YC} (page \pageref{Diagram YC}), lower part.

See also the discussion in
Section 
\ref{Section Saw-Blades-Discussion} (page 
\pageref{Section Saw-Blades-Discussion})  and
the construction of saw blades in
Section \ref{Section Saw-Blades} (page \pageref{Section Saw-Blades}).

 \xDH
As we work for a contradiction in the case x-, too, the simplest
way to achieve this is to have a negative arrow $x \xcP x',$ and at $x' $
again
an $ \xcu.$ This gives a chance to construct a contradiction at $x'.$
Of course, we have to construct a contradiction at $y,$ too, as in
the case $x-,$ we have an $ \xco $ at $x.$

Of course, we may have branches originating at $x',$ which all lead to
contradictions in the case $x-,$ so $x \xcp x' $ (resp. $ \xco $ at $x')$
is
possible, too.

But, for the simple construction, we need $ \xcP $ between $x$ and $x',$
and
$ \xcu $ at $x'.$ And this leads to the construction of contradictions
for all $x' \xcp z$ (or $x' \xcP z)$ as just mentioned above.

 \xEj

\er

\be

$\hspace{0.01em}$


\label{Example Procrastination}

This Example shows that
infinitely many finitely branching points cannot always
replace infinite branching - there is an infinite
``procrastination branch'' or ``escape branch''.
This modification of the Yablo structure has one acceptable
valuation for $Y_{1}:$

Let $Y_{i},$ $i< \xbo $ as usual, and introduce new $X_{i},$ $3 \xck i<
\xbo.$

Let $Y_{i} \xcP Y_{i+1},$ $Y_{i} \xcp X_{i+2},$ $X_{i} \xcP Y_{i},$ $X_{i}
\xcp X_{i+1},$ with

$Y_{i}:=$ $ \xCN Y_{i+1} \xcu X_{i+2},$ $X_{i}:= \xCN Y_{i} \xcu X_{i+1}.$

If $Y_{1}= \xct,$ then $ \xCN Y_{2} \xcu X_{3},$ by $X_{3},$ $ \xCN Y_{3}
\xcu X_{4},$ so, generally,

if $Y_{i}= \xct,$ then $\{ \xCN Y_{j}:$ $i<j\}$ and $\{X_{j}:$ $i+1<j\}.$

If $ \xCN Y_{1},$ then $Y_{2} \xco \xCN X_{3},$ so if $ \xCN X_{3},$
$Y_{3} \xco \xCN X_{4},$ etc., so, generally,

if $ \xCN Y_{i},$ then $ \xcE j(i<j,$ $Y_{j})$ or $ \xcA j\{ \xCN X_{j}:$
$i+1<j\}.$

Suppose now $Y_{1}= \xct,$ then $X_{j}$ for all $2<j,$ and $ \xCN Y_{j}$
for all $1<j.$
By $ \xCN Y_{2}$ there is $j,$ $2<j,$ and $Y_{j},$ a contradiction, or $
\xCN X_{j}$ for all $3<j,$
again a contradiction.

But $ \xCN Y_{1}$ is possible, by setting $ \xCN Y_{i}$ and $ \xCN X_{i}$
for all $i.$

Thus, replacing infinite branching by an infinite number of finite
branching does not work for the Yablo construction, as we can always chose
the
``procrastinating'' branch.

See Diagram 
\ref{Diagram Procrastination} (page 
\pageref{Diagram Procrastination}).

\clearpage

\begin{diagram}

Escape Path

\label{Diagram Procrastination}
\index{Diagram Procrastination}

\unitlength1.0mm
\begin{picture}(140,200)(0,0)

\put(0,0){\line(1,0){140}}
\put(0,0){\line(0,1){200}}
\put(140,200){\line(-1,0){140}}
\put(140,200){\line(0,-1){200}}

\put(30,175){{\rm\bf Diagram for Example \ref{Example Procrastination}}}

\put(30,160){$Y_1$}
\put(32,158){\line(0,-1){32}}
\put(33,158){\line(3,-4){23}}
\put(30,132){\line(1,0){4}}

\put(50,160){Lines represent downward (negated) arrows}

\put(30,120){$Y_2$}
\put(32,118){\line(0,-1){32}}
\put(55,120){$X_3$}
\put(57,118){\line(0,-1){32}}
\put(33,118){\line(3,-4){23}}
\put(56,118){\line(-3,-4){23}}
\put(30,92){\line(1,0){4}}
\put(36,94){\line(4,-3){3}}

\put(30,80){$Y_3$}
\put(32,78){\line(0,-1){32}}
\put(55,80){$X_4$}
\put(57,78){\line(0,-1){32}}
\put(33,78){\line(3,-4){23}}
\put(56,78){\line(-3,-4){23}}
\put(30,52){\line(1,0){4}}
\put(36,54){\line(4,-3){3}}

\put(30,40){$Y_4$}
\put(32,38){\line(0,-1){32}}
\put(55,40){$X_5$}
\put(57,38){\line(0,-1){32}}
\put(33,38){\line(3,-4){23}}
\put(56,38){\line(-3,-4){23}}
\put(30,12){\line(1,0){4}}
\put(36,14){\line(4,-3){3}}

\put(30,2){$Y_5$}
\put(55,2){$X_6$}

\end{picture}

\end{diagram}

\vspace{4mm}

\clearpage



\clearpage
\section{
A Generalization of Yablo's Construction to Formulas in
Disjunctive Normal Form - Examples
}

\label{Section Or-And}
\subsection{
Introduction
}

\ee

We discuss here a general strategy and a number of examples
(in Example \ref{Example Or-And} (page \pageref{Example Or-And})).

They have in common that the formulas are of the type $ \xcO \xcU,$ i.e.
in
disjunctive normal form. The limiting cases are pure conjunctions
(as in Yablo's original approach) or pure disjunctions.

The examples are straightforward generalizations of Yablo's construction,
as we have here several columns, in Yablo's construction just one
column, and our choice functions $g$ (in all columns) correspond to the
choice
of one element in Yablo's construction.

Consider the first example below.

More precisely, as Yablo works with $ \xcU x_{i},$ there is one uniform
set of $x_{i}.$
We work with $ \xcO \xcU x_{i,j},$ so we have to distinguish the elements
in the $ \xcU $
from the sets in $ \xcO.$ We define for this purpose columns, whose
elements
are the elements in the $ \xcU,$ and the set of columns are the sets in $
\xcO.$
Negation is now slightly more complicated, not just OR of negated
elements,
but
OR of choice functions of negated elements in the columns.
We also need some enumeration of $ \xbo \xCK \xbo,$ or of a suitable set.

It is perhaps easiest to see the following example geometrically.
We have a vertical column $C_{i_{0}}$ where a certain property holds,
and a horizontal line (the choice function $g),$ where the opposite
holds, and column and horizontal line meet.

More precisely, the property will not necessarily hold in all of
$C_{i_{0}},$
but only from a certain height onward. As the choice functions $g$ will
chose even higher in the columns, the clash is assured.

For more details of the general strategy, see
Definition \ref{Definition Or-And} (page \pageref{Definition Or-And}),
Case \ref{Case DOA} (page \pageref{Case DOA}).
and
Example \ref{Example Or-And} (page \pageref{Example Or-And}),
in particular
Case \ref{Case OA6} (page \pageref{Case OA6}).
\subsection{
The Examples
}

We first introduce some notation and definitions.

\bd

$\hspace{0.01em}$


\label{Definition Or-And}

 \xEh
 \xDH
The examples will differ in the size of the $ \xcU $ and $ \xcO $ -
where the case of both finite
is trivial, no contradictions possible - and the relation in the
graph, i.e. which formulas are ``visible'' from a given formula.
(In Yablo's construction, the $x_{j},$ $j>i,$ are visible from $x_{i}.)$
We will write the relation by $ \xcP,$ $x \xcP y,$ but for simplicity
sometimes $x<y,$ too,
reminding us that variables inside the
formulas will be negated.
We consider linear and ranked orders, leaving general partial orders
aside.
All relations will be transitive and acyclic.

For the intuition (and beyond) we order the formulas in sets of
columns. Inside a column, the formulas are connected by $ \xcU,$ the
columns themselves are connected by $ \xcO.$ This is the logical
ordering, it is different from above order relation in the graph.

 \xDH The basic structure.

 \xEh
 \xDH
So, we have columns $C_{i},$ and inside the columns variables $x_{i,j}.$
As we might have finitely or countably infinitely many columns, of finite
or countably infinite size, we
write the set of columns $C:=\{C_{i}:i< \xba \},$ where $ \xba < \xbo +1,$
and
$C_{i}:=\{x_{i,j}:j< \xbb_{i}\},$ where $ \xbb_{i}< \xbo +1$ again. By
abuse of language, $C$ will
also denote the whole construction, $C:=\{x_{i,j}:i,j< \xbo \}.$
 \xDH
Given $C$ and $x_{i,j},$ $C \xex x_{i,j}:=\{x_{i',j' }:x_{i,j} \xcP x_{i'
,j' }\},$ the part of $C$ visible from
$x_{i,j}.$
 \xDH
Likewise, $C_{i' } \xex x_{i,j}$ $:=$ $(C \xex x_{i,j}) \xcs C_{i' }.$
 \xEj

 \xDH

Let $I(x_{i,j})$ $:=$ $\{i' < \xba:$ $C_{i' } \xex x_{i,j} \xEd \xCQ \}.$
(In some $i',$ there might be no $x_{i',j' }$ s.t. $x_{i,j} \xcP x_{i'
,j' }.)$

Given $i' \xbe I(x_{i,j}),$ let
$J(i',x_{i,j})$ $:=$ $\{j' < \xbb_{i' }:x_{i,j} \xcP x_{i',j' }\}.$
By $i' \xbe I(x_{i,j}),$ $J(i',x_{i,j}) \xEd \xCQ.$

 \xDH
Back to logic. Let $x_{i,j}:= \xcO \{ \xcU \{ \xCN x_{i',j' }:j' \xbe
J(i',x_{i,j})\}:i' \xbe I(x_{i,j})\},$
$(\xcu $ inside columns, $ \xco $ between columns).
We will sometimes abbreviate $d(x_{i,j})$ by $x_{i,j}.$
(Note that the $x_{i',j' }$ are exactly the elements visible from
$x_{i,j}.)$

 \xEh
 \xDH
If $x_{i,j}$ is true (written $x_{i,j}+),$ then all elements in one of the
$C_{i' },$
$i' \xbe I(x_{i,j}),$
and visible from $x_{i,j},$ must all be false - but we do not know in
which $C_{i' }.$
We denote this $C_{i' }$ by $C_{i(x_{i,j})},$ and define
$C[x_{i,j}]:=C_{i(x_{i,j})} \xex x_{i,j}.$
 \xDH
Conversely, suppose $x_{i,j}$ is false, $x_{i,j}-.$
Again, we consider only elements in $C \xex x_{i,j},$ i.e. visible from
$x_{i,j}.$
By distributivity, $x_{i,j}-$ $=$ $ \xcO \{ \xcU ran(g):g \xbe \xbP
\{C_{i' } \xex x_{i,j}:i' \xbe I(x_{i,j})\}\}.$

(g is a choice function chosing in all columns $C_{i' }$ the ``sufficiently
big''
elements $x_{k,m},$ i.e. above $x_{i,j},$
$ran(g)$ its range or image.)
See Fact \ref{Fact DNF} (page \pageref{Fact DNF}).

Note that the elements of $ran(g)$ are now positive!

The $ \xcO $ in above formula choses some such function $g,$ but we do not
know which.
Let $g[x_{i,j}]$ denote the chosen one.
 \xDH
Note that both $C[x_{i,j}]$ and $g[x_{i,j}]$ are undefined if there are no
$x_{i',j' }>x_{i,j}.$
 \xEj

 \xDH

\label{Case DOA}

The conflicts will be between the ``vertical'' (negative) columns, and
``horizontal'' (positive) lines of the $g.$ It is a very graphical
construction. If we start with a positive point, the
negative columns correspond to the $ \xcA $ in Yablo's
construction, the $g$ to the $ \xcE,$ if we start with a negative one, it
is
the other way round.

More precisely:

 \xEh
 \xDH Let $x_{i,j}+,$ consider $C[x_{i,j}]$ (the chosen column) -
recall all elements of $C[x_{i,j}]$ are negative.

If there is $x_{i',j' } \xbe C[x_{i,j}]$ which is not maximal in
$C[x_{i,j}],$ then
$g[x_{i',j' }]$ intersects $C[x_{i,j}],$ a contradiction, as all elements
of $ran(g[x_{i',j' }])$
are positive.

Of course, such non-maximal $x_{i',j' }$ need not exist. In that case, we
have to
try again with the negative element $x_{i',j' }$ and Case (5.2).

 \xDH
Let $x_{i,j}-,$ consider $g[x_{i,j}]$ (the chosen function) - recall, all
elements of $ran(g[x_{i,j}])$ are positive.

If there is $x_{i' j' } \xbe ran(g[x_{i,j}])$ such that $C[x_{i',j' }]
\xcs ran(g[x_{i,j}]) \xEd \xCQ,$
we have a contradiction.
So $x_{i,j} \xcP x_{i',j' } \xcP x_{i'',j'' } \xbe C[x_{i',j' }].$
(This case is impossible in Yablo's original construction.)

Otherwise, we have to try to work with the elements in $ran(g[x_{i,j}])$
and Case (5.1) above, etc.

Note that we can chose a suitable $x_{i',j' } \xbe ran(g[x_{i,j}]),$ in
particular
one which is not maximal (if such exist), but we have no control over
the choice of $C[x_{i',j' }],$ so we might have to exhaust all finite
columns, until the choice is only from a set of infinite columns.

 \xEj

See Diagram 
\ref{Diagram Or-And-1} (page 
\pageref{Diagram Or-And-1}), upper part, and
Example \ref{Example Or-And} (page \pageref{Example Or-And}),
Case \ref{Case OA6} (page \pageref{Case OA6}).

 \xEj

\ed

\be

$\hspace{0.01em}$


\label{Example Or-And}

 \xEh

 \xDH Case 1

\label{Case OA1}

Consider the structure $C:=\{x_{i,j}:i< \xba,j< \xbo \},$ with columns
$C_{i}:=\{x_{i,j}:j< \xbo \}.$

Take a standard enumeration $f$ of $C,$ e.g. $f(0):=x_{0,0},$ then
enumerate the $x_{i,j}$
s.t. $max\{i,j\}=1,$ then $max\{i,j\}=2,$ etc. As $f$ is bijective,
$f^{-1}(x_{i,j})$ is
defined.

(More precisely, let $m:=max\{i,j\},$ and e.g. we go first horizontally
from left to
right over the columns up to column $m-1,$ then in column $m$ upwards,
i.e.
$x_{0,m}, \Xl,x_{m-1,m},x_{m,0}, \Xl,x_{m,m}$.)

Define the relation $x_{i,j} \xcP x_{i',j' }$ iff
$f^{-1}(x_{i,j})<f^{-1}(x_{i',j' }).$ Obviously,
$ \xcP $ is transitive and free from cycles.
$C \xex k:=\{x_{i,j} \xbe C:f^{-1}(x_{i,j}) \xCe k\}$ etc. are defined.

We now show that the structure has no truth values.

Suppose $x_{i,j}+.$

Consider $C[x_{i,j}],$ let $i':=i(x_{i,j})$ and chose $x_{i',j' } \xbe
C[x_{i,j}].$
$x_{i',j' }$ is false, $ran(g[x_{i',j' }])$ intersects $C_{i' }$ above
$x_{i',j' },$ so we have a contradiction.

In particular, $x_{0,0}+$ is impossible.

Suppose $x_{0,0}-,$ then $ran(g[x_{0,0}]) \xEd \xCQ,$ chose $x_{i,j} \xbe
ran(g[x_{0,0}]),$ so $x_{i,j}+,$
but we saw that this is impossible.

Note: for $ \xba =1,$ we have Yablo's construction.

 \xDH Case 2

\label{Case OA2}

We now show that the same construction with columns of height 2
does not work, it has an escape path.

Set $x_{i,0}+,$ $x_{i,1}-$ for all $i.$ $C_{i(x_{i,0})}$ might be $C_{i},$
so $C[x_{i,0}]=\{x_{i,1}\},$ which
is possible. $ran(g[x_{i,1}])$ might be $\{x_{j,0}:j>i\},$ which is
possible again.

 \xDH
Case 3

\label{Case OA3}

We change the order in
Case \ref{Case OA2} (page \pageref{Case OA2}),
to a (horizontally) ranked order:

$x_{i,j}<x_{i',j' }$ iff $i<i' $ (thus, between columns), we now have a
contradiction:

Consider $x_{i,0}-$ and $g[x_{i,0}].$ $g[x_{i,0}](i)$ is undefined. Let
$x_{i+1,j' }:=g[x_{i,0}](i+1),$ thus $x_{i+1,j' }+.$
Consider $C[x_{i+1,j' }].$ This must be some $C_{i' }$ for $i' >i+1,$
so $C[x_{i+1,j' }] \xcs ran(g[x_{i,0}]) \xEd \xCQ,$ and we have a
contradiction.

For $x_{0,0}+,$ consider $C[x_{0,0}],$ this must be some $C_{i},$ $i>0,$
take $x_{i,0} \xbe C_{i},$ $x_{i,0}-,$ and
continue as above.

 \xDH
Case 4

\label{Case OA4}

Modify
Case \ref{Case OA1} (page \pageref{Case OA1}),
to a ranked order, but this time vertically:
$x_{i,j}<x_{i',j' }$ iff $j<j'.$

Take $x_{i,j}+,$ consider $C[x_{i,j}],$ this may be (part of) any $C_{i'
}$ (beginning
at $j+1).$ Take e.g. $x_{i',j+1}- \xbe C[x_{i,j}],$ consider $g[x_{i'
,j+1}],$
a choice function in $\{C_{k} \xex x_{i',j+1}:k, \xbo \}$
this will intersect $C[x_{i,j}].$
In particular, $x_{0,0}+$ is impossible.

Suppose $x_{0,0}-,$ take $x_{i,j}+ \xbe ran(g[x_{0,0}]),$ and continue as
above.

Note that the case with just one $C_{i}$ is the original Yablo
construction.

 \xDH
Case 5

\label{Case OA5}

We modify
Case \ref{Case OA1} (page \pageref{Case OA1})
again:
Consider the ranked order by $x_{i,j} \xcP x_{i',j' }$ iff
$max\{i,j\}<max\{i',j' \}.$

This is left to the reader as an exercise.

 \xDH Case 6

\label{Case OA6}

The general argument is as follows (and applies to general partial
orders, too):
 \xEh
 \xDH
We show that $x_{i,j}+$ leads to a contradiction.
(In our terminology, $x_{i,j}$ is the head.)
 \xEh
 \xDH
we find $C[x_{i,j}]$ (negative elements)
above $x_{i,j}$ - but we have no control over the choice of $C[x_{i,j}],$
 \xDH we chose $x_{i',j' } \xbe C[x_{i,j}]$ - if there is a minimal such,
chose this one,
it must not be a maximal element in $C[x_{i,j}]$ $(x_{i',j' }$ is the
knee),
 \xDH we find $g[x_{i',j' }]$ (positive elements)
above $x_{i',j' }$ - again we have no control over the
choice of this $g.$ But, as $C[x_{i,j}] \xex x_{i',j' }$ is not empty,
$ran(g[x_{i',j' }]) \xcs C[x_{i,j}] \xEd \xCQ,$ so we have a
contradiction
$(x_{i'',j'' } \xbe ran(g[x_{i',j' }]) \xcs C[x_{i,j}]$ is the foot).
 \xDH
We apply the reasoning to $x_{0,0}.$
 \xEj
 \xDH
We show that $x_{0,0}-$ leads to a contradiction.
 \xEh
 \xDH We have $g[x_{0,0}]$ (positive elements) -
but no control over the choice of $g.$ In particular,
it may be arbitrarily high up.
 \xDH we chose $x_{i,j} \xbe ran(g[x_{0,0}])$ with enough room above it
for the
argument about $x_{i,j}+$ in (6.1).
 \xEj

 \xDH Case 
\ref{Case OA3} (page 
\pageref{Case OA3})  is similar, except that we work
horizontally, not vertically.

 \xEj

 \xEj

\ee

\br

$\hspace{0.01em}$


\label{Remark Condition}

Note that we may use here the ideas of
Condition 
\ref{Condition Yablo} (page 
\pageref{Condition Yablo}), based on the discussion
in Section 
\ref{Section Saw-Blades-Discussion} (page 
\pageref{Section Saw-Blades-Discussion}),
to construct and check constructions for arbitrary formulas in
disjunctive normal form.

\er

$ \xCO $

$ \xCO $

$ \xCO $

\clearpage

\begin{diagram}

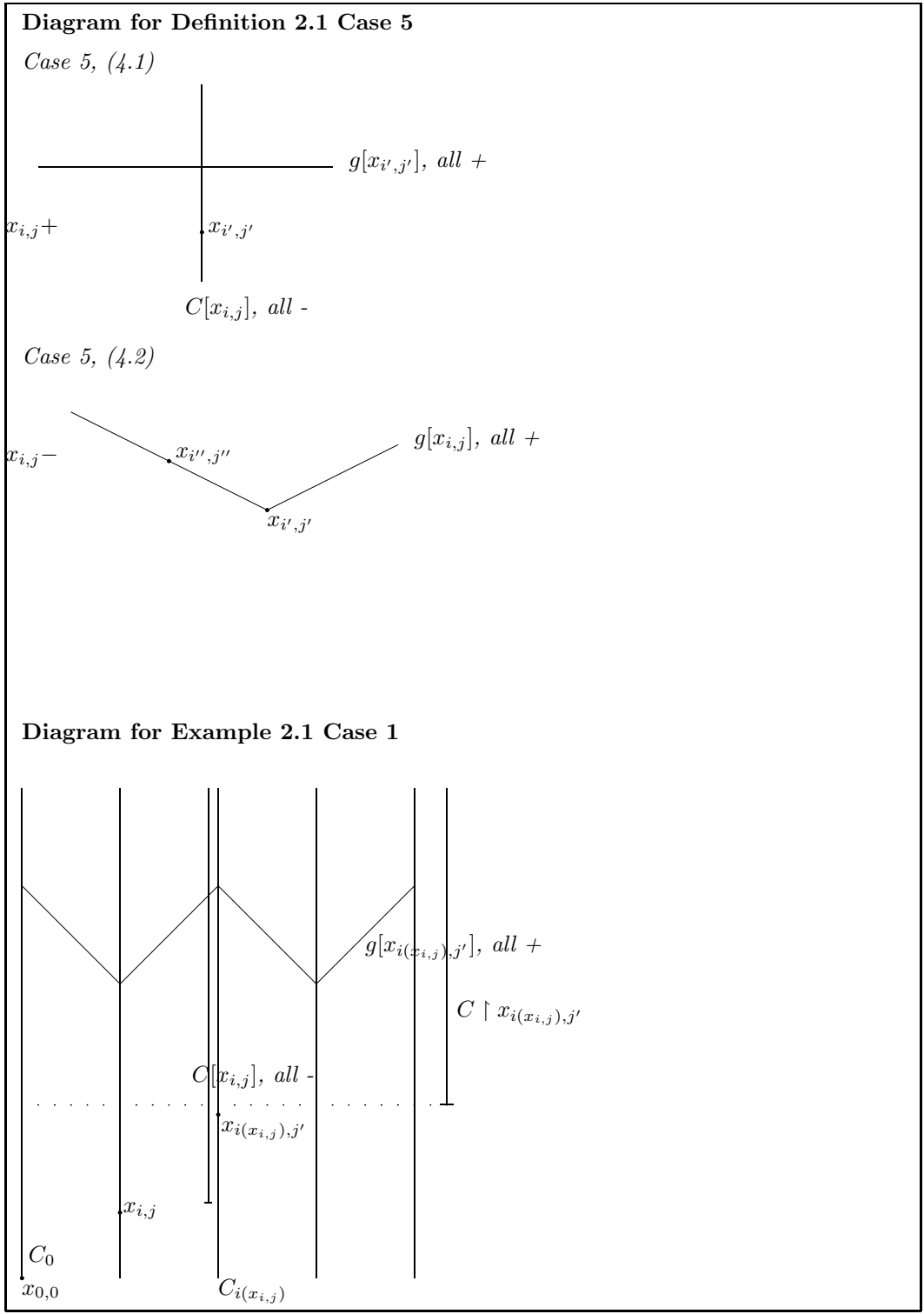
Contradictions for Disjunctive Normal Forms

\label{Diagram Or-And-1}
\index{Diagram Or-And-1}

\unitlength0.5mm
\begin{picture}(280,400)(0,0)

\put(0,0){\line(1,0){280}}
\put(0,0){\line(0,1){400}}
\put(280,400){\line(-1,0){280}}
\put(280,400){\line(0,-1){400}}

\put(5,392){{\rm\bf Diagram for Definition
    \ref{Definition Or-And} Case \ref{Case DOA}}}

\put(5,380){Case \ref{Case DOA}, (4.1)}

\put(60,315){\line(0,1){60}}
\put(10,350){\line(1,0){90}}
\put(60,330){\circle*{1}}
\put(62,330){$x_{i',j'}$}
\put(55,305){$C[x_{i,j}]$, all -}
\put(0,330){$x_{i,j}+$}
\put(105,350){$g[x_{i',j'}]$, all +}

\put(5,290){Case \ref{Case DOA}, (4.2)}

\put(80,245){\circle*{1}}
\put(80,240){$x_{i',j'}$}
\put(80,245){\line(-2,1){60}}
\put(80,245){\line(2,1){40}}
\put(0,260){$x_{i,j}-$}
\put(125,265){$g[x_{i,j}]$, all +}
\put(50,260){\circle*{1}}
\put(52,261){$x_{i'',j''}$}

\put(5,175){{\rm\bf Diagram for Example
    \ref{Example Or-And} Case \ref{Case OA1}}}

\put(5,10){\line(0,1){150}}
\put(35,10){\line(0,1){150}}
\put(65,10){\line(0,1){150}}
\put(95,10){\line(0,1){150}}
\put(125,10){\line(0,1){150}}

\put(5,10){\circle*{1}}
\put(7,15){$C_0$}
\put(5,5){$x_{0,0}$}
\put(65,5){$C_{i(x_{i,j})}$}
\put(36,30){$x_{i,j}$}
\put(35,30){\circle*{1}}

\put(66,55){$x_{i(x_{i,j}),j'}$}
\put(65,60){\circle*{1}}
\put(57,70){$C[x_{i,j}]$, all -}
\put(62,33){\line(0,1){127}}
\put(61,33){\line(1,0){2}}

\put(135,63){\line(0,1){97}}
\multiput(135,63)(-5,0){26}{\circle*{0.2}}
\put(133,63){\line(1,0){4}}
\put(138,90){$C \xex x_{i(x_{i,j}),j'}$}

\put(5,130){\line(1,-1){30}}
\put(65,130){\line(1,-1){30}}
\put(35,100){\line(1,1){30}}
\put(95,100){\line(1,1){30}}
\put(110,110){$g[x_{i(x_{i,j}),j'}]$, all +}

\end{picture}

\end{diagram}

\vspace{4mm}

\clearpage

\clearpage
\section{
Saw Blades
}

\label{Section Saw-Blades-All}
\subsection{
Introduction
}

Yablo works with contradictions in the form of $x_{0}= \xCN x_{1} \xcu
\xCN x_{2},$ $x_{1}= \xCN x_{2},$
graphically $x_{0} \xcP x_{1} \xcP x_{2},$ $x_{0} \xcP x_{2}.$ They are
combined in a formally simple total
order, which, however, blurs conceptual differences.

We discuss here different, conceptually very clear and simple, examples
of a Yablo-like construction.

In particular, we emphasize the difference between $x_{1}$ and $x_{2}$ in
the
Yablo contradictions. The contradiction is
finished in $x_{2},$ but not in $x_{1},$ requiring the barring of escape
routes in $x_{1}.$
We repair the possible escape routes by constructing new contradictions
for the same origin. This is equivalent to closing under transitivity in
the
individual ``saw blades'' - see below.

The main difference to Yablo's construction is that we first construct a
contradiction
$x_{0} \xcP x_{1} \xcP y_{0},$ $x_{0} \xcP y_{0},$ and later $x_{1} \xcP
x_{2},$ and then $x_{0} \xcP x_{2},$ like Yablo, so
our construction is first more liberal, but we later see that we have
to continue like Yablo.

So we use the same ``cells'' as Yablo does for the contradictions, but
analyse
the way they are put together.

We use the full strength of the conceptual difference between $x_{1}$ and
$x_{2}$
(in above notation)
only in Section 
\ref{Section Simplifications} (page 
\pageref{Section Simplifications}), where we show that
preventing
$x_{2}$ from being TRUE is sufficient, whereas we need $x_{1}$ to be
contradictory, see also
Section 
\ref{Section Saw-Blades-Discussion} (page 
\pageref{Section Saw-Blades-Discussion}).
Thus, we obtain a minimal, i.e. necessary and sufficient, construction
for combining Yablo cells in this way.

We analyse the constructions (Yablo's original construction and our
variants)
in detail, and see how they follow from the prerequisites $(x_{0}+$ is
contradictory, and any arrow $x_{0} \xcP z$ leads again to a contradiction
in case
$z+$ (and thus $x_{0}-)$).
\subsection{
Saw Blades
}

\label{Section Saw-Blades}

First, we show the escape route problem.

\be

$\hspace{0.01em}$


\label{Example Escape}

Consider Construction 
\ref{Construction Saw-Blade} (page 
\pageref{Construction Saw-Blade})
without closing under transitivity, i.e. the only arrows originating
in $x_{ \xbs,0}$ will be $x_{ \xbs,0} \xcP x_{ \xbs,1}$ and $x_{ \xbs
,0} \xcP y_{ \xbs,0},$ etc.

Let $x_{ \xbs,0}=TRUE,$ then $x_{ \xbs,1}$ is an $ \xco,$ and we
pursue the path $x_{ \xbs,0} \xcP x_{ \xbs,1} \xcP x_{ \xbs,2},$ this
has no contradiction so far, and we
continue with $x_{ \xbs,2}=TRUE,$ $x_{ \xbs,3}$ is $ \xco $ again,
we continue and have $x_{ \xbs,0} \xcP x_{ \xbs,1} \xcP x_{ \xbs,2}
\xcP x_{ \xbs,3} \xcP x_{ \xbs,4},$ etc.,
never meeting a contradiction, so we have an escape path.

\ee

$ \xCO $

\bcs

$\hspace{0.01em}$


\label{Construction Saw-Blade}

We construct a saw blade $ \xbs,$ $SB_{ \xbs }.$

 \xEh

 \xDH
``Saw Blades''
 \xEh
 \xDH
Let $x_{ \xbs,0} \xcP x_{ \xbs,1} \xcP x_{ \xbs,2} \xcP x_{ \xbs,3}
\xcP x_{ \xbs,4},$  \Xl.

$x_{ \xbs,0} \xcP y_{ \xbs,0},$ $x_{ \xbs,1} \xcP y_{ \xbs,0},$ $x_{
\xbs,1} \xcP y_{ \xbs,1},$ $x_{ \xbs,2} \xcP y_{ \xbs,1},$
$x_{ \xbs,2} \xcP y_{ \xbs,2},$ $x_{ \xbs,3} \xcP y_{ \xbs,2},$ $x_{
\xbs,3} \xcP y_{ \xbs,3},$ $x_{ \xbs,4} \xcP y_{ \xbs,3},$  \Xl.

we call the construction
a ``saw blade'', with ``teeth'' $y_{ \xbs,0},$ $y_{ \xbs,1},$ $y_{ \xbs
,2},$  \Xl.
and ``back'' $x_{ \xbs,0},$ $x_{ \xbs,1},$ $x_{ \xbs,2},$  \Xl.

We call $x_{ \xbs,0}$ the start of the blade.

See Diagram \ref{Diagram Sawblade-2} (page \pageref{Diagram Sawblade-2}),
upper part.

 \xDH
Add (against escape), e.g. first $x_{ \xbs,0} \xcP x_{ \xbs,2},$ $x_{
\xbs,0} \xcP y_{ \xbs,1},$
then $x_{ \xbs,1} \xcP x_{ \xbs,3},$ $x_{ \xbs,1} \xcP y_{ \xbs,2},$
now
we have to add $x_{ \xbs,0} \xcP x_{ \xbs,3},$ $x_{ \xbs,0} \xcP y_{
\xbs,2},$ etc, recursively.
This is equivalent to
closing the saw blade under transitivity with negative arrows $ \xcP.$
This is easily seen.

 \xDH
We set $x_{ \xbs,i}= \xcU \xCN z_{ \xbs,j},$ for all $x_{ \xbs,j}$
such that $x_{ \xbs,i} \xcP z_{ \xbs,j},$
as in the original Yablo construction.

 \xEj

 \xDH
Composition of saw blades
 \xEh
 \xDH
Add for the teeth $y_{ \xbs,0},$ $y_{ \xbs,1},$ $y_{ \xbs,2}$  \Xl.
their own saw blades, i.e.
start at $y_{ \xbs,0}$ a new saw blade $SB_{ \xbs,0}$ with $y_{ \xbs
,0}=x_{ \xbs,0,0},$ at
$y_{ \xbs,1}$ a new saw blade $SB_{ \xbs,1}$ with $y_{ \xbs,1}=x_{ \xbs
,1,0},$ etc.

 \xDH
Do this recursively.

I.e., at every tooth of every saw blade start a new saw blade.
See Diagram \ref{Diagram Sawblade-3} (page \pageref{Diagram Sawblade-3}).
 \xEj
 \xEj

Note:

It is NOT necessary to close the whole structure (the individual saw
blades
together) under transitivity.

\ecs

\bfa

$\hspace{0.01em}$


\label{Fact Saw-Blade}

All $z_{ \xbs,i}$ in all saw blades so constructed are contradictory,
i.e.
assigning them a truth value leads to a contradiction.

\efa

\subparagraph{
Proof
}

$\hspace{0.01em}$


The argument is almost the same as for Yablo's construction.

Fix some saw blade $SB_{ \xbs }$ in the construction.

 \xEh

 \xDH
Take any $z_{ \xbs,i}$ with $z_{ \xbs,i}+,$ i.e. $z_{ \xbs,i}=TRUE.$ We
show that this is contradictory.

 \xEh

 \xDH
Case 1: $z_{ \xbs,i}$ is one of the $x_{ \xbs,i}$ i.e. it is in the back
of the blade.

Take any $x_{ \xbs,i' }$ in the back such that there is an arrow $x_{
\xbs,i} \xcP x_{ \xbs,i' }$
$(i':=i+1$ suffices).
Then $x_{ \xbs,i' }=FALSE,$ and we have an $ \xco $ at $x_{ \xbs,i' }.$
Take any $z_{ \xbs,j}$ such that $x_{ \xbs,i' } \xcP z_{ \xbs,j},$ by
transitivity, $x_{ \xbs,i} \xcP z_{ \xbs,j},$ so $z_{ \xbs,j}=FALSE,$
but as $x_{ \xbs,i' }=FALSE,$ $z_{ \xbs,j}=TRUE,$
contradiction.

 \xDH
Case 2: $z_{ \xbs,i}$ is one of the $y_{ \xbs,i},$ i.e. a tooth of the
blade.

Then $y_{ \xbs,i}$ is the start of the new blade starting at $y_{ \xbs
,i},$ and we
argue as above in Case 1.

 \xEj

 \xDH
Take any $z_{ \xbs,i}$ with $z_{ \xbs,i}-,$ i.e. $z_{ \xbs,i}=FALSE,$
and we have an $ \xco $ at $z_{ \xbs,i},$ and
one of the successors of $z_{ \xbs,i},$ say $z_{ \xbs,j},$ has to be
TRUE.
We just saw that this is impossible.

(For the intuition:
If $z_{ \xbs,i}$ is in the back of the blade, all of its successors are
in the
same blade.
If $z_{ \xbs,i}$ is one of the teeth of the blade, all of its successors
are
in the new blade, starting at $z_{ \xbs,i}.$
In both cases, $z_{ \xbs,j}=TRUE$ leads to a contradiction, as we saw
above.)

 \xEj

\br

$\hspace{0.01em}$


\label{Remark Saw-Blade}

Note that all $y_{ \xbs,i}$ are contradictory, too, not only the $x_{
\xbs,i}.$
We will see in
Section \ref{Section Simplifications} (page \pageref{Section Simplifications})
that we can achieve this by simpler means, as we need to consider here
the case $x_{ \xbs,i} \xco $ only, the contradiction for the case $x_{
\xbs,i} \xcu $ is already
treated.

Thus, we seemingly did not fully use here the conceptual clarity of
difference
between $x_{1}$ and $x_{2}$ alluded to in the beginning of
Section \ref{Section Saw-Blades} (page \pageref{Section Saw-Blades}).
See, however, the discussion in
Section 
\ref{Section Saw-Blades-Discussion} (page 
\pageref{Section Saw-Blades-Discussion}).

\er

$ \xCO $

$ \xCO $

\clearpage

\begin{diagram}

Saw Blade

\label{Diagram Sawblade-2}
\index{Diagram Sawblade-2}

\unitlength0.6mm
\begin{picture}(230,210)(0,0)

\put(0,0){\line(1,0){230}}
\put(0,0){\line(0,1){210}}
\put(230,210){\line(-1,0){230}}
\put(230,210){\line(0,-1){210}}

\put(5,197){Horizontal lines stand for negative arrows from left to right,}
\put(5,190){the other lines for negative arrows from top to bottom}

\put(5,182){{\rm\bf Diagram Single Saw Blade }}
\put(5,175){Start of the saw blade $\sigma$ beginning at $x_{\sigma,0}$,}
\put(5,168){before closing under transitivity }

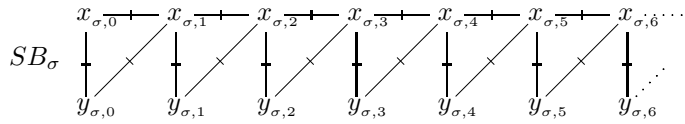
\put(5,150){$SB_{\sigma}$}


\put(20,160){$x_{_{\sigma,0}}$}
\put(20,140){$y_{_{\sigma,0}}$}
\put(26,161){\line(1,0){12}}
\put(32,160){\line(0,1){2}}
\put(22,143){\line(0,1){14}}
\put(21,150){\line(1,0){2}}
\put(24,143){\line(1,1){16}}
\put(31,151.5){\line(1,-1){1.5}}


\put(40,160){$x_{_{\sigma,1}}$}
\put(40,140){$y_{_{\sigma,1}}$}
\put(46,161){\line(1,0){12}}
\put(52,160){\line(0,1){2}}
\put(42,143){\line(0,1){14}}
\put(41,150){\line(1,0){2}}
\put(44,143){\line(1,1){16}}
\put(51,151.5){\line(1,-1){1.5}}


\put(60,160){$x_{_{\sigma,2}}$}
\put(60,140){$y_{_{\sigma,2}}$}
\put(66,161){\line(1,0){12}}
\put(72,160){\line(0,1){2}}
\put(62,143){\line(0,1){14}}
\put(61,150){\line(1,0){2}}
\put(64,143){\line(1,1){16}}
\put(71,151.5){\line(1,-1){1.5}}


\put(80,160){$x_{_{\sigma,3}}$}
\put(80,140){$y_{_{\sigma,3}}$}
\put(86,161){\line(1,0){12}}
\put(92,160){\line(0,1){2}}
\put(82,143){\line(0,1){14}}
\put(81,150){\line(1,0){2}}
\put(84,143){\line(1,1){16}}
\put(91,151.5){\line(1,-1){1.5}}


\put(100,160){$x_{_{\sigma,4}}$}
\put(100,140){$y_{_{\sigma,4}}$}
\put(106,161){\line(1,0){12}}
\put(112,160){\line(0,1){2}}
\put(102,143){\line(0,1){14}}
\put(101,150){\line(1,0){2}}
\put(104,143){\line(1,1){16}}
\put(111,151.5){\line(1,-1){1.5}}


\put(120,160){$x_{_{\sigma,5}}$}
\put(120,140){$y_{_{\sigma,5}}$}
\put(126,161){\line(1,0){12}}
\put(132,160){\line(0,1){2}}
\put(122,143){\line(0,1){14}}
\put(121,150){\line(1,0){2}}
\put(124,143){\line(1,1){16}}
\put(131,151.5){\line(1,-1){1.5}}


\put(140,160){$x_{_{\sigma,6}}$}
\put(140,140){$y_{_{\sigma,6}}$}
\multiput(146,161)(2,0){5}{\circle*{0.2}}
\put(142,143){\line(0,1){14}}
\put(141,150){\line(1,0){2}}
\multiput(144,143)(1.5,1.5){5}{\circle*{0.2}}

\put(5,125){Read $x_{\xbs,0} \xcP x_{\xbs,1} \xcP y_{\xbs,0}$,
$x_{\xbs,0} \xcP y_{\xbs,0}$, etc., more precisely:}

\put(5,118){$x_{\xbs,0} = \xCN x_{\xbs,1} \xcu \xCN y_{\xbs,0}$,
$x_{\xbs,1} = \xCN y_{\xbs,0} \xcu \xCN x_{\xbs,2} \xcu \xCN y_{\xbs,1}$,
etc.}

\put(5,82){{\rm\bf Diagram Simplified Saw Blade }}
\put(5,75){Start of the saw blade before closing}
\put(5,68){the blade (without ``decoration'') under transitivity }


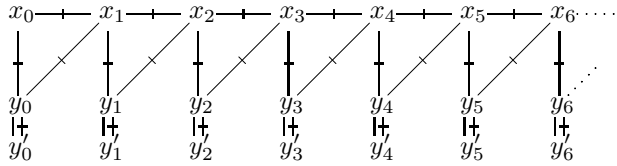
\put(20,60){$x_0$}
\put(20,40){$y_0$}
\put(20,30){$y'_0$}
\put(26,61){\line(1,0){12}}
\put(32,60){\line(0,1){2}}
\put(22,43){\line(0,1){14}}
\put(21,34){\line(0,1){4}}
\put(23,34){\line(0,1){4}}
\put(22,36){\line(1,0){2}}
\put(21,50){\line(1,0){2}}
\put(24,43){\line(1,1){16}}
\put(31,51.5){\line(1,-1){1.5}}


\put(40,60){$x_1$}
\put(40,40){$y_1$}
\put(40,30){$y'_1$}
\put(46,61){\line(1,0){12}}
\put(52,60){\line(0,1){2}}
\put(42,43){\line(0,1){14}}
\put(41,34){\line(0,1){4}}
\put(43,34){\line(0,1){4}}
\put(42,36){\line(1,0){2}}
\put(41,50){\line(1,0){2}}
\put(44,43){\line(1,1){16}}
\put(51,51.5){\line(1,-1){1.5}}


\put(60,60){$x_2$}
\put(60,40){$y_2$}
\put(60,30){$y'_2$}
\put(66,61){\line(1,0){12}}
\put(72,60){\line(0,1){2}}
\put(62,43){\line(0,1){14}}
\put(61,34){\line(0,1){4}}
\put(63,34){\line(0,1){4}}
\put(62,36){\line(1,0){2}}
\put(61,50){\line(1,0){2}}
\put(64,43){\line(1,1){16}}
\put(71,51.5){\line(1,-1){1.5}}


\put(80,60){$x_3$}
\put(80,40){$y_3$}
\put(80,30){$y'_3$}
\put(86,61){\line(1,0){12}}
\put(92,60){\line(0,1){2}}
\put(82,43){\line(0,1){14}}
\put(81,34){\line(0,1){4}}
\put(83,34){\line(0,1){4}}
\put(82,36){\line(1,0){2}}
\put(81,50){\line(1,0){2}}
\put(84,43){\line(1,1){16}}
\put(91,51.5){\line(1,-1){1.5}}


\put(100,60){$x_4$}
\put(100,40){$y_4$}
\put(100,30){$y'_4$}
\put(106,61){\line(1,0){12}}
\put(112,60){\line(0,1){2}}
\put(102,43){\line(0,1){14}}
\put(101,34){\line(0,1){4}}
\put(103,34){\line(0,1){4}}
\put(102,36){\line(1,0){2}}
\put(101,50){\line(1,0){2}}
\put(104,43){\line(1,1){16}}
\put(111,51.5){\line(1,-1){1.5}}


\put(120,60){$x_5$}
\put(120,40){$y_5$}
\put(120,30){$y'_5$}
\put(126,61){\line(1,0){12}}
\put(132,60){\line(0,1){2}}
\put(122,43){\line(0,1){14}}
\put(121,34){\line(0,1){4}}
\put(123,34){\line(0,1){4}}
\put(122,36){\line(1,0){2}}
\put(121,50){\line(1,0){2}}
\put(124,43){\line(1,1){16}}
\put(131,51.5){\line(1,-1){1.5}}


\put(140,60){$x_6$}
\put(140,40){$y_6$}
\put(140,30){$y'_6$}
\multiput(146,61)(2,0){5}{\circle*{0.2}}
\put(142,43){\line(0,1){14}}
\put(141,34){\line(0,1){4}}
\put(143,34){\line(0,1){4}}
\put(142,36){\line(1,0){2}}
\put(141,50){\line(1,0){2}}
\multiput(144,43)(1.5,1.5){5}{\circle*{0.2}}

\put(5,15){Read $y_{0} =  y'_{0} \xcu \xCN y'_{0}$ etc.}

\end{picture}

\end{diagram}

\vspace{10mm}

\clearpage

\clearpage

\begin{diagram}

Composed Saw Blades

\label{Diagram Sawblade-3}
\index{Diagram Sawblade-3}

\unitlength0.8mm
\begin{picture}(170,210)(0,0)

\put(0,0){\line(1,0){170}}
\put(0,0){\line(0,1){210}}
\put(170,210){\line(-1,0){170}}
\put(170,210){\line(0,-1){210}}

\put(5,197){Horizontal lines stand for negative arrows from left to right,}
\put(5,190){the other lines for negative arrows from top to bottom}

\put(5,182){{\rm\bf Diagram Composition of Saw Blades }}
\put(5,175){Composition of saw blades (without additional arrows) }
\put(5,168){The fat dots indicate identity, e.g. $y_{0,0}=x_{0,0,0}$ }

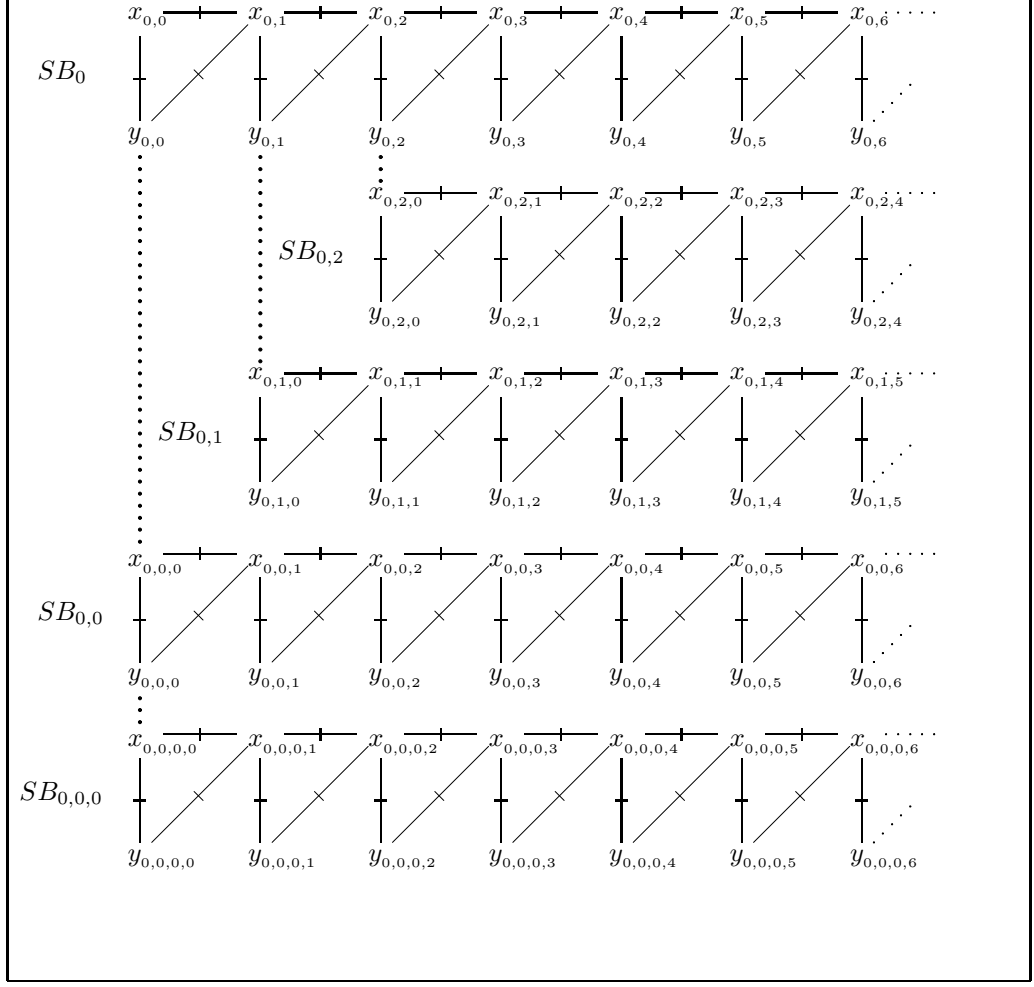
\put(5,150){$SB_{0}$}


\put(20,160){$x_{_{0,0}}$}
\put(20,140){$y_{_{0,0}}$}
\put(26,161){\line(1,0){12}}
\put(32,160){\line(0,1){2}}
\put(22,143){\line(0,1){14}}
\put(21,150){\line(1,0){2}}
\put(24,143){\line(1,1){16}}
\put(31,151.5){\line(1,-1){1.5}}


\put(40,160){$x_{_{0,1}}$}
\put(40,140){$y_{_{0,1}}$}
\put(46,161){\line(1,0){12}}
\put(52,160){\line(0,1){2}}
\put(42,143){\line(0,1){14}}
\put(41,150){\line(1,0){2}}
\put(44,143){\line(1,1){16}}
\put(51,151.5){\line(1,-1){1.5}}


\put(60,160){$x_{_{0,2}}$}
\put(60,140){$y_{_{0,2}}$}
\put(66,161){\line(1,0){12}}
\put(72,160){\line(0,1){2}}
\put(62,143){\line(0,1){14}}
\put(61,150){\line(1,0){2}}
\put(64,143){\line(1,1){16}}
\put(71,151.5){\line(1,-1){1.5}}


\put(80,160){$x_{_{0,3}}$}
\put(80,140){$y_{_{0,3}}$}
\put(86,161){\line(1,0){12}}
\put(92,160){\line(0,1){2}}
\put(82,143){\line(0,1){14}}
\put(81,150){\line(1,0){2}}
\put(84,143){\line(1,1){16}}
\put(91,151.5){\line(1,-1){1.5}}


\put(100,160){$x_{_{0,4}}$}
\put(100,140){$y_{_{0,4}}$}
\put(106,161){\line(1,0){12}}
\put(112,160){\line(0,1){2}}
\put(102,143){\line(0,1){14}}
\put(101,150){\line(1,0){2}}
\put(104,143){\line(1,1){16}}
\put(111,151.5){\line(1,-1){1.5}}


\put(120,160){$x_{_{0,5}}$}
\put(120,140){$y_{_{0,5}}$}
\put(126,161){\line(1,0){12}}
\put(132,160){\line(0,1){2}}
\put(122,143){\line(0,1){14}}
\put(121,150){\line(1,0){2}}
\put(124,143){\line(1,1){16}}
\put(131,151.5){\line(1,-1){1.5}}


\put(140,160){$x_{_{0,6}}$}
\put(140,140){$y_{_{0,6}}$}
\multiput(146,161)(2,0){5}{\circle*{0.2}}
\put(142,143){\line(0,1){14}}
\put(141,150){\line(1,0){2}}
\multiput(144,143)(1.5,1.5){5}{\circle*{0.2}}

\put(45,120){$SB_{0,2}$}

\multiput(62,133)(0,2){3}{\circle*{0.5}}


\put(60,130){$x_{_{0,2,0}}$}
\put(60,110){$y_{_{0,2,0}}$}
\put(66,131){\line(1,0){12}}
\put(72,130){\line(0,1){2}}
\put(62,113){\line(0,1){14}}
\put(61,120){\line(1,0){2}}
\put(64,113){\line(1,1){16}}
\put(71,121.5){\line(1,-1){1.5}}


\put(80,130){$x_{_{0,2,1}}$}
\put(80,110){$y_{_{0,2,1}}$}
\put(86,131){\line(1,0){12}}
\put(92,130){\line(0,1){2}}
\put(82,113){\line(0,1){14}}
\put(81,120){\line(1,0){2}}
\put(84,113){\line(1,1){16}}
\put(91,121.5){\line(1,-1){1.5}}


\put(100,130){$x_{_{0,2,2}}$}
\put(100,110){$y_{_{0,2,2}}$}
\put(106,131){\line(1,0){12}}
\put(112,130){\line(0,1){2}}
\put(102,113){\line(0,1){14}}
\put(101,120){\line(1,0){2}}
\put(104,113){\line(1,1){16}}
\put(111,121.5){\line(1,-1){1.5}}


\put(120,130){$x_{_{0,2,3}}$}
\put(120,110){$y_{_{0,2,3}}$}
\put(126,131){\line(1,0){12}}
\put(132,130){\line(0,1){2}}
\put(122,113){\line(0,1){14}}
\put(121,120){\line(1,0){2}}
\put(124,113){\line(1,1){16}}
\put(131,121.5){\line(1,-1){1.5}}


\put(140,130){$x_{_{0,2,4}}$}
\put(140,110){$y_{_{0,2,4}}$}
\multiput(146,131)(2,0){5}{\circle*{0.2}}
\put(142,113){\line(0,1){14}}
\put(141,120){\line(1,0){2}}
\multiput(144,113)(1.5,1.5){5}{\circle*{0.2}}

\put(25,90){$SB_{0,1}$}

\multiput(42,103)(0,2){18}{\circle*{0.5}}


\put(40,100){$x_{_{0,1,0}}$}
\put(40,80){$y_{_{0,1,0}}$}
\put(46,101){\line(1,0){12}}
\put(52,100){\line(0,1){2}}
\put(42,83){\line(0,1){14}}
\put(41,90){\line(1,0){2}}
\put(44,83){\line(1,1){16}}
\put(51,91.5){\line(1,-1){1.5}}


\put(60,100){$x_{_{0,1,1}}$}
\put(60,80){$y_{_{0,1,1}}$}
\put(66,101){\line(1,0){12}}
\put(72,100){\line(0,1){2}}
\put(62,83){\line(0,1){14}}
\put(61,90){\line(1,0){2}}
\put(64,83){\line(1,1){16}}
\put(71,91.5){\line(1,-1){1.5}}


\put(80,100){$x_{_{0,1,2}}$}
\put(80,80){$y_{_{0,1,2}}$}
\put(86,101){\line(1,0){12}}
\put(92,100){\line(0,1){2}}
\put(82,83){\line(0,1){14}}
\put(81,90){\line(1,0){2}}
\put(84,83){\line(1,1){16}}
\put(91,91.5){\line(1,-1){1.5}}


\put(100,100){$x_{_{0,1,3}}$}
\put(100,80){$y_{_{0,1,3}}$}
\put(106,101){\line(1,0){12}}
\put(112,100){\line(0,1){2}}
\put(102,83){\line(0,1){14}}
\put(101,90){\line(1,0){2}}
\put(104,83){\line(1,1){16}}
\put(111,91.5){\line(1,-1){1.5}}


\put(120,100){$x_{_{0,1,4}}$}
\put(120,80){$y_{_{0,1,4}}$}
\put(126,101){\line(1,0){12}}
\put(132,100){\line(0,1){2}}
\put(122,83){\line(0,1){14}}
\put(121,90){\line(1,0){2}}
\put(124,83){\line(1,1){16}}
\put(131,91.5){\line(1,-1){1.5}}


\put(140,100){$x_{_{0,1,5}}$}
\put(140,80){$y_{_{0,1,5}}$}
\multiput(146,101)(2,0){5}{\circle*{0.2}}
\put(142,83){\line(0,1){14}}
\put(141,90){\line(1,0){2}}
\multiput(144,83)(1.5,1.5){5}{\circle*{0.2}}

\put(5,60){$SB_{0,0}$}

\multiput(22,73)(0,2){33}{\circle*{0.5}}


\put(20,69){$x_{_{0,0,0}}$}
\put(20,50){$y_{_{0,0,0}}$}
\put(26,71){\line(1,0){12}}
\put(32,70){\line(0,1){2}}
\put(22,53){\line(0,1){14}}
\put(21,60){\line(1,0){2}}
\put(24,53){\line(1,1){16}}
\put(31,61.5){\line(1,-1){1.5}}


\put(40,69){$x_{_{0,0,1}}$}
\put(40,50){$y_{_{0,0,1}}$}
\put(46,71){\line(1,0){12}}
\put(52,70){\line(0,1){2}}
\put(42,53){\line(0,1){14}}
\put(41,60){\line(1,0){2}}
\put(44,53){\line(1,1){16}}
\put(51,61.5){\line(1,-1){1.5}}


\put(60,69){$x_{_{0,0,2}}$}
\put(60,50){$y_{_{0,0,2}}$}
\put(66,71){\line(1,0){12}}
\put(72,70){\line(0,1){2}}
\put(62,53){\line(0,1){14}}
\put(61,60){\line(1,0){2}}
\put(64,53){\line(1,1){16}}
\put(71,61.5){\line(1,-1){1.5}}


\put(80,69){$x_{_{0,0,3}}$}
\put(80,50){$y_{_{0,0,3}}$}
\put(86,71){\line(1,0){12}}
\put(92,70){\line(0,1){2}}
\put(82,53){\line(0,1){14}}
\put(81,60){\line(1,0){2}}
\put(84,53){\line(1,1){16}}
\put(91,61.5){\line(1,-1){1.5}}


\put(100,69){$x_{_{0,0,4}}$}
\put(100,50){$y_{_{0,0,4}}$}
\put(106,71){\line(1,0){12}}
\put(112,70){\line(0,1){2}}
\put(102,53){\line(0,1){14}}
\put(101,60){\line(1,0){2}}
\put(104,53){\line(1,1){16}}
\put(111,61.5){\line(1,-1){1.5}}


\put(120,69){$x_{_{0,0,5}}$}
\put(120,50){$y_{_{0,0,5}}$}
\put(126,71){\line(1,0){12}}
\put(132,70){\line(0,1){2}}
\put(122,53){\line(0,1){14}}
\put(121,60){\line(1,0){2}}
\put(124,53){\line(1,1){16}}
\put(131,61.5){\line(1,-1){1.5}}


\put(140,69){$x_{_{0,0,6}}$}
\put(140,50){$y_{_{0,0,6}}$}
\multiput(146,71)(2,0){5}{\circle*{0.2}}
\put(142,53){\line(0,1){14}}
\put(141,60){\line(1,0){2}}
\multiput(144,53)(1.5,1.5){5}{\circle*{0.2}}

\put(2,30){$SB_{0,0,0}$}

\multiput(22,43)(0,2){3}{\circle*{0.5}}


\put(20,39){$x_{_{0,0,0,0}}$}
\put(20,20){$y_{_{0,0,0,0}}$}
\put(26,41){\line(1,0){12}}
\put(32,40){\line(0,1){2}}
\put(22,23){\line(0,1){14}}
\put(21,30){\line(1,0){2}}
\put(24,23){\line(1,1){16}}
\put(31,31.5){\line(1,-1){1.5}}


\put(40,39){$x_{_{0,0,0,1}}$}
\put(40,20){$y_{_{0,0,0,1}}$}
\put(46,41){\line(1,0){12}}
\put(52,40){\line(0,1){2}}
\put(42,23){\line(0,1){14}}
\put(41,30){\line(1,0){2}}
\put(44,23){\line(1,1){16}}
\put(51,31.5){\line(1,-1){1.5}}


\put(60,39){$x_{_{0,0,0,2}}$}
\put(60,20){$y_{_{0,0,0,2}}$}
\put(66,41){\line(1,0){12}}
\put(72,40){\line(0,1){2}}
\put(62,23){\line(0,1){14}}
\put(61,30){\line(1,0){2}}
\put(64,23){\line(1,1){16}}
\put(71,31.5){\line(1,-1){1.5}}


\put(80,39){$x_{_{0,0,0,3}}$}
\put(80,20){$y_{_{0,0,0,3}}$}
\put(86,41){\line(1,0){12}}
\put(92,40){\line(0,1){2}}
\put(82,23){\line(0,1){14}}
\put(81,30){\line(1,0){2}}
\put(84,23){\line(1,1){16}}
\put(91,31.5){\line(1,-1){1.5}}


\put(100,39){$x_{_{0,0,0,4}}$}
\put(100,20){$y_{_{0,0,0,4}}$}
\put(106,41){\line(1,0){12}}
\put(112,40){\line(0,1){2}}
\put(102,23){\line(0,1){14}}
\put(101,30){\line(1,0){2}}
\put(104,23){\line(1,1){16}}
\put(111,31.5){\line(1,-1){1.5}}


\put(120,39){$x_{_{0,0,0,5}}$}
\put(120,20){$y_{_{0,0,0,5}}$}
\put(126,41){\line(1,0){12}}
\put(132,40){\line(0,1){2}}
\put(122,23){\line(0,1){14}}
\put(121,30){\line(1,0){2}}
\put(124,23){\line(1,1){16}}
\put(131,31.5){\line(1,-1){1.5}}


\put(140,39){$x_{_{0,0,0,6}}$}
\put(140,20){$y_{_{0,0,0,6}}$}
\multiput(146,41)(2,0){5}{\circle*{0.2}}
\put(142,23){\line(0,1){14}}
\put(141,30){\line(1,0){2}}
\multiput(144,23)(1.5,1.5){5}{\circle*{0.2}}

\end{picture}

\end{diagram}

\vspace{4mm}

\clearpage

\clearpage
\subsection{
Discussion of Saw Blades
}

\label{Section Saw-Blades-Discussion}

We analyse here the construction of saw blades from the prerequisites
and show that Yablo's construction is a necessary substructure.
The analysis applies as well to Yablo's original structure - it suffices
to omit the $y$'s.

It will be illustrated in
Diagram \ref{Diagram X1} (page \pageref{Diagram X1})  and
Diagram \ref{Diagram X2} (page \pageref{Diagram X2}).
They will also show that the negative arrows are necessary - though this
is not important here.
\subsubsection{
Preliminaries
}

Yablo's structure is sufficient, it shows that $x_{0}$ there does not have
a consistent truth value.

We constructed
in Section \ref{Section Saw-Blades} (page \pageref{Section Saw-Blades})
a structure from the prerequisites which is more general,
and contains Yablo's structure as a substructure. This shows that
Yablo's structure is necessary (and minimal) under our prerequisites.

As we are mainly interested in constructing Yablo's structure from the
prerequisites, we will stop analysing the rest of the saw blade structure
after
some initial steps, in order to simplify things.

One of the consequences of the construction is that we will have
a constructive proof of the nesessity of infinite depth and width.

We will work as usual with conjunctions of negations.

Recall:

\bd

$\hspace{0.01em}$


\label{Definition Terminology}

 \xEh
 \xDH
A Yablo cell (YC) has the form $x \xcP x' \xcP y,$ $x \xcP y.$

$x$ is the head, $x' $ the knee, $y$ the foot.

Without any branching at $x',$ $x+$ is impossible, as it is
contradictory.

If $x+,$ then $x' - \xco.$ Thus, if we add some $x' \xcP y' $ (or $x'
\xcp y'),$ we have
$x' -=y \xco y',$ so $x+$ is not contradictory any more.
Any $x' \xcP y' $ opens a new possibility for $x+$ to be consistent,
although $x$
is the head of a YC.

 \xDH
If, for all $x' \xcP y',$ we also have $x \xcP y',$ $x+$ is again
impossible.

We call such a system of $x \xcP x',$ and for all $x' \xcP y_{i}$ also $x
\xcP y_{i}$ a
Yablo Cell System (YCS) - of course, with head $x$ and knee $x'.$

Thus, if we want an contradiction at $x,$ a YC with head $x$ does not
suffice,
we need a YCS with head $x.$

 \xDH
If we discuss general contradictory cells, we will write CC instead of YC,
and CCS instead of YCS.

 \xEj
\subsubsection{
Prerequisites
}

\label{Section Prereq}

\ed

We use the following prerequisites:

 \xEh
 \xDH
$x_{0}$ is contradictory
 \xDH
If there is an arrow $x_{0} \xcP x,$ then $x$ is contradictory.
 \xEj

Translated to YC's and $YCS$'s, this means:

 \xEh
 \xDH
$x_{0}$ is the head of a YCS (for $x_{0}+)$
 \xDH
For any $x_{0} \xcP x,$ and $x_{0}-,$ thus $x+ \xcU,$ $x+$ is the head of
a YCS.
 \xEj
\subsubsection{
Strategy and Details of the Construction
}

\label{Section Strat}

\label{Section DetailsConstr}

During the inductive construction, we will append new $YC$'s, say $z \xcP
z' \xcP z'',$
$z \xcP z'',$ to the construction done so far, but we will see that
prequisite (2) in Section 
\ref{Section Prereq} (page 
\pageref{Section Prereq})  will force us to add new
arrows to $z',$ say $z' \xcP z''',$ so the YCS property at $z$ is
destroyed,
and we will have to repair the YCS by adding a new arrow $z \xcP z'''.$

The final construction is the union of infinitely many steps,
so all YCS's which were destroyed at some step will be repaired later
again.

\bcs

$\hspace{0.01em}$


\label{Construction Details}

 \xEh

 \xDH
Start:

We begin with $x_{0} \xcP x_{1} \xcP y_{0},$ $x_{0} \xcP y_{0}.$

See Diagram \ref{Diagram X1} (page \pageref{Diagram X1}), Top.

This is a YCS with head $x_{0},$ thus $x_{0}$ is contradictory, and
prerequisite (1)
in Section \ref{Section Prereq} (page \pageref{Section Prereq})  is satisfied.

 \xDH
Prerequisite (2) is not satisfied, and we add two new YC's to $x_{1}$ and
$y_{0}$
to satisfy prerequisite (2):

$x_{1} \xcP x_{2} \xcP y_{1},$ $x_{1} \xcP y_{1},$ and

$y_{0} \xcP z_{0} \xcP z_{0}',$ $y_{0} \xcP z_{0}'.$

See Diagram \ref{Diagram X1} (page \pageref{Diagram X1}), Center, only
$x_{1} \xcP x_{2} \xcP y_{1},$ $x_{1} \xcP y_{1}$ shown.

 \xDH
The new arrows $x_{1} \xcP x_{2}$ and $x_{1} \xcP y_{1}$ destroy the YCS
at $x_{0}$ with knee $x_{1}.$

We could, of course, begin an new YC at $x_{0}$ with different knee, say
$x_{0} \xcP x_{1}' \xcP y_{0}',$ which is a YCS, but this would take us
back to step (1),
step (2) would have been a dead end, and we would be in the same
situation again. Thus, this does not lead anywhere.

Thus, we have to repair the situation of step (2) by adding
$x_{0} \xcP x_{2}$ and $x_{0} \xcP y_{1},$ and we have a YCS again.

See Diagram \ref{Diagram X1} (page \pageref{Diagram X1}), Bottom.

 \xDH
Prerequisite (2) is again not satisfied, and we add two new $YC$'s,
to $x_{2}$ and $y_{1},$ say

$x_{2} \xcP x_{3} \xcP y_{2},$ $x_{2} \xcP y_{2},$ and

$y_{1} \xcP z_{1} \xcP z_{1}',$ $y_{1} \xcP z_{1}'.$

See Diagram \ref{Diagram X2} (page \pageref{Diagram X2}), Top, only
$x_{2} \xcP x_{3} \xcP y_{2},$ $x_{2} \xcP y_{2}$ shown.

Note:

It might be possible to find or add arrows $x_{2} \xcP u$ and $x_{2} \xcP
u' $ and $u \xcP u' $
to already
constructed $u,$ $u',$ without introducing cycles. We continue with $u,$
$u',$
adding arrows $x_{0} \xcP u,$ $x_{0} \xcP u',$ etc., but as the
construction already
achieved is finite, in the end, we will have to add new points and arrows.
More precisely, such $u,$ $u' $ cannot be some $x_{i}$ already
constructed, this would
create a cycle, it has to be some $y$ or $z.$ We will neglect
the y's and $z' $ anyway, as our main
interest here lies in the $x_{i}$ and the arrows between them.

 \xDH
We now want to construct arrows $x_{0} \xcP x_{3}$ and $x_{0} \xcP y_{2}$
using a broken YCS
property at $x_{0},$ but we have here a new knee, $x_{2}$ and not $x_{1}$
(which is
the knee for $x_{0}).$ So, we have to proceed indirectly.

 \xDH
The new arrows $x_{2} \xcP y_{2}$ and $x_{2} \xcP x_{3}$ destroy the YCS
at $x_{1}$ with knee $x_{2},$
and we add $x_{1} \xcP x_{3}$ and $x_{1} \xcP y_{2}$ to repair the YCS.

See Diagram \ref{Diagram X2} (page \pageref{Diagram X2}), Center.

 \xDH
But, now, we go from $x_{0}$ via the knee, have new arrows
$x_{1} \xcP x_{3}$ and $x_{1} \xcP y_{2},$ and have to repair the YCS of
$x_{0}$ by adding
$x_{0} \xcP x_{3}$ and $x_{0} \xcP y_{2}.$

See Diagram \ref{Diagram X2} (page \pageref{Diagram X2}), Bottom.

 \xDH
We use here in (6) and (7) the general property of the construction
that $x_{i+1}$ is the knee of the YCS with head $x_{i}.$

Consequently, if we add an arrow $x_{i+1} \xcp x,$ then we have to add an
arrow $x_{i} \xcP x$
to repair the YCS with head $x_{i}.$

The argument goes downward until $x_{0}$ (as $x_{i}$ is the knee of the
YCS with head
$x_{i-1}$ etc.), so, if we add an arrow $x_{i+1} \xcP x,$ we have to add
new arrows
$x_{k} \xcP x$ for all $0 \xck k \xck x_{i}.$

 \xDH
 \xEh
 \xDH
Thus, there are arrows $x_{i} \xcP x_{j}$ for all $i,j,$ $i<j,$ and the
construction is
transitive for the $x_{i}' s.$

 \xDH
Every $x_{i}$ is head of a YCS with knee $x_{i+1}.$

 \xDH
Thus, every arrow from any $x_{i}$ to any $x_{j}$ goes to the head of a
YCS, and
not only the arrows from $x_{0}.$

This property is ``accidental'', and due to the fact that for any arrow
$x_{i} \xcP x_{j},$ there is also an arrow $x_{0} \xcP x_{j},$ and
property (2) holds for $x_{0}$
by prerequisite.

 \xEj

 \xDH
The construction has infinite depth and branching.

 \xEj

\clearpage

\begin{diagram}

Inductive Construction 1

\label{Diagram X1}
\index{Diagram X1}

\unitlength0.8mm
\begin{picture}(170,210)(0,0)

\put(0,0){\line(1,0){170}}
\put(0,0){\line(0,1){210}}
\put(170,210){\line(-1,0){170}}
\put(170,210){\line(0,-1){210}}

\put(5,197){Horizontal lines stand for negative arrows from left to right,}
\put(5,190){the other lines for negative arrows from top to bottom}

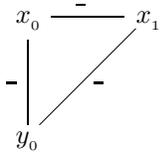
\put(80,160){Contradiction at $x_0+$}
\put(80,154){(If $x_0-$, all arrows from $x_0$ have to lead to a}
\put(80,150){contradiction, so $x_0 \xcP y_0$ and $x_0 \xcP x_1$ }
\put(80,146){have to be negative arrows.}
\put(80,142){But if $x_0+$, we need a contradiction in }
\put(80,138){$ \{ x_0,x_1,y_0 \}$, so $x_1 \xcP y_0$ has to be}
\put(80,134){negative, too.)}


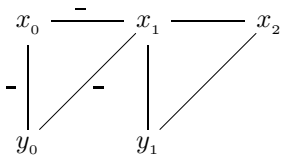
\put(10,160){$x_{_{0}}$}
\put(10,140){$y_{_{0}}$}
\put(16,161){\line(1,0){12}}
\put(12,143){\line(0,1){14}}
\put(14,143){\line(1,1){16}}

\put(20,163){\line(1,0){1.5}}
\put(8.5,150){\line(1,0){1.5}}
\put(23,150){\line(1,0){1.5}}


\put(30,160){$x_{_{1}}$}


\put(80,110){Contradiction for $x_0-, x_0.. x_1$ has to lead to}
\put(80,105){new Yablo Cell}


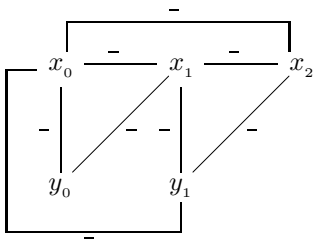
\put(10,110){$x_{_{0}}$}
\put(10,90){$y_{_{0}}$}
\put(16,111){\line(1,0){12}}
\put(12,93){\line(0,1){14}}
\put(14,93){\line(1,1){16}}

\put(20,113){\line(1,0){1.5}}
\put(8.5,100){\line(1,0){1.5}}
\put(23,100){\line(1,0){1.5}}


\put(30,110){$x_{_{1}}$}
\put(30,90){$y_{_{1}}$}
\put(36,111){\line(1,0){12}}
\put(32,93){\line(0,1){14}}
\put(34,93){\line(1,1){16}}


\put(50,110){$x_{_{2}}$}


\put(80,60){Repair contradiction at $x_0+$}
\put(80,56){(due to new branches at $x_1-)$}
\put(80,52){In future, we will not draw the new arrows }
\put(80,48){to the $y_i's$}

\put(80,40){(Again, $x_0 \xcP x_2$ and $x_0 \xcP y_1$ have to be }
\put(80,36){negative. But, as $\{x_0,x_1,y_1\}$ and }
\put(80,32){$\{x_0,x_1,x_2\}$ have to be contradictory, $x_1 \xcP y_1$ }
\put(80,28){and $x_1 \xcP x_2$ have to be negative, too.)}


\put(10,60){$x_{_{0}}$}
\put(10,40){$y_{_{0}}$}
\put(16,61){\line(1,0){12}}
\put(12,43){\line(0,1){14}}
\put(14,43){\line(1,1){16}}

\put(13,63){\line(0,1){5}}
\put(50,63){\line(0,1){5}}
\put(13,68){\line(1,0){37}}

\put(8,60){\line(-1,0){5}}
\put(3,60){\line(0,-1){27}}
\put(3,33){\line(1,0){29}}
\put(32,38){\line(0,-1){5}}


\put(30,60){$x_{_{1}}$}
\put(30,40){$y_{_{1}}$}
\put(36,61){\line(1,0){12}}
\put(32,43){\line(0,1){14}}
\put(34,43){\line(1,1){16}}


\put(50,60){$x_{_{2}}$}

\put(20,63){\line(1,0){1.5}}
\put(8.5,50){\line(1,0){1.5}}
\put(23,50){\line(1,0){1.5}}

\put(40,63){\line(1,0){1.5}}
\put(28.5,50){\line(1,0){1.5}}
\put(43,50){\line(1,0){1.5}}

\put(30,70){\line(1,0){1.5}}
\put(16,32){\line(1,0){1.5}}


\end{picture}

\end{diagram}

\vspace{10mm}

\clearpage

\clearpage

\begin{diagram}

Inductive Construction 2

\label{Diagram X2}
\index{Diagram X2}

\unitlength0.7mm
\begin{picture}(200,210)(0,0)

\put(0,0){\line(1,0){200}}
\put(0,0){\line(0,1){210}}
\put(200,210){\line(-1,0){200}}
\put(200,210){\line(0,-1){210}}

\put(5,197){Horizontal lines stand for negative arrows from left to right,}
\put(5,190){the other lines for negative arrows from top to bottom}

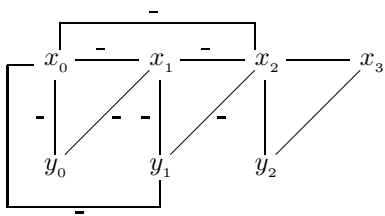
\put(80,160){Repair contradiction at $x_0-$: new branch $x_0.. x_2$ has }
\put(80,156){to lead to a new Yablo Cell}


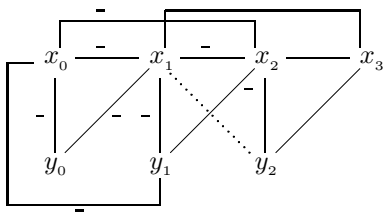
\put(10,160){$x_{_{0}}$}
\put(10,140){$y_{_{0}}$}
\put(16,161){\line(1,0){12}}
\put(12,143){\line(0,1){14}}
\put(14,143){\line(1,1){16}}

\put(13,163){\line(0,1){5}}
\put(50,163){\line(0,1){5}}
\put(13,168){\line(1,0){37}}


\put(30,160){$x_{_{1}}$}
\put(30,140){$y_{_{1}}$}
\put(36,161){\line(1,0){12}}
\put(32,143){\line(0,1){14}}
\put(34,143){\line(1,1){16}}


\put(50,160){$x_{_{2}}$}
\put(50,140){$y_{_{2}}$}
\put(56,161){\line(1,0){12}}
\put(52,143){\line(0,1){14}}
\put(54,143){\line(1,1){16}}


\put(70,160){$x_{_{3}}$}

\put(20,163){\line(1,0){1.5}}
\put(8.5,150){\line(1,0){1.5}}
\put(23,150){\line(1,0){1.5}}

\put(40,163){\line(1,0){1.5}}
\put(28.5,150){\line(1,0){1.5}}
\put(43,150){\line(1,0){1.5}}

\put(30,170){\line(1,0){1.5}}
\put(16,132){\line(1,0){1.5}}

\put(8,160){\line(-1,0){5}}
\put(3,160){\line(0,-1){27}}
\put(3,133){\line(1,0){29}}
\put(32,138){\line(0,-1){5}}


\put(80,110){Repair contradiction at $x_1+$ due to new branch $x_2.. x_3$}
\put(80,104){As $x_1$ is the head of a Yablo Cell, we have the dotted}
\put(80,100){line from $x_1$ to $y_2$}


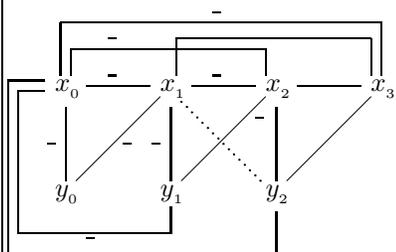
\put(10,110){$x_{_{0}}$}
\put(10,90){$y_{_{0}}$}
\put(16,111){\line(1,0){12}}
\put(12,93){\line(0,1){14}}
\put(14,93){\line(1,1){16}}

\put(13,113){\line(0,1){5}}
\put(50,113){\line(0,1){5}}
\put(13,118){\line(1,0){37}}


\put(30,110){$x_{_{1}}$}
\put(30,90){$y_{_{1}}$}
\put(36,111){\line(1,0){12}}
\put(32,93){\line(0,1){14}}
\put(34,93){\line(1,1){16}}

\put(33,113){\line(0,1){7}}
\put(70,113){\line(0,1){7}}
\put(33,120){\line(1,0){37}}

\multiput(34,108)(1,-1){16}{\circle*{0.2}}


\put(50,110){$x_{_{2}}$}
\put(50,90){$y_{_{2}}$}
\put(56,111){\line(1,0){12}}
\put(52,93){\line(0,1){14}}
\put(54,93){\line(1,1){16}}


\put(70,110){$x_{_{3}}$}

\put(20,113){\line(1,0){1.5}}
\put(8.5,100){\line(1,0){1.5}}
\put(23,100){\line(1,0){1.5}}

\put(40,113){\line(1,0){1.5}}
\put(28.5,100){\line(1,0){1.5}}
\put(48,105){\line(1,0){1.5}}

\put(20,120){\line(1,0){1.5}}
\put(16,82){\line(1,0){1.5}}

\put(8,110){\line(-1,0){5}}
\put(3,110){\line(0,-1){27}}
\put(3,83){\line(1,0){29}}
\put(32,88){\line(0,-1){5}}


\put(80,60){Repair contradiction at $x_0+$ due to new branch $x_1.. x_3$}

\put(80,54){(Because of the dotted line $x_1$ to $y_2$, we also need a }

\put(80,50){line $x_0$ to $y_2$. It has to be negative (for the case $x_0-$).}

\put(80,44){As $\{x_0,x_1,y_2\}$ is contradictory, $x_1 \xcP y_2$ is negative.}

\put(80,38){As $\{x_0,x_1,x_3\}$ is contradictory, $x_1 \xcP x_3$ is negative.}

\put(80,32){As $\{x_1,x_2,x_3\}$ is contradictory, $x_2 \xcP x_3$ is negative.}

\put(80,26){But as $x_2$ is the head of a Yablo Cell,}

\put(80,22){and $x_2 \xcP y_2$, $x_2 \xcP x_3$, so $x_3 \xcP y_2$, too.}

\put(80,14){Thus, ALL lines in the diagram are negative.)}


\put(10,60){$x_{_{0}}$}
\put(10,40){$y_{_{0}}$}
\put(16,61){\line(1,0){12}}
\put(12,43){\line(0,1){14}}
\put(14,43){\line(1,1){16}}

\put(13,63){\line(0,1){5}}
\put(50,63){\line(0,1){5}}
\put(13,68){\line(1,0){37}}

\put(11,63){\line(0,1){10}}
\put(72,63){\line(0,1){10}}
\put(11,73){\line(1,0){61}}

\put(40,75){\line(1,0){1.5}}


\put(30,60){$x_{_{1}}$}
\put(30,40){$y_{_{1}}$}
\put(36,61){\line(1,0){12}}
\put(32,43){\line(0,1){14}}
\put(34,43){\line(1,1){16}}

\put(33,63){\line(0,1){7}}
\put(70,63){\line(0,1){7}}
\put(33,70){\line(1,0){37}}

\multiput(34,58)(1,-1){16}{\circle*{0.2}}


\put(50,60){$x_{_{2}}$}
\put(50,40){$y_{_{2}}$}
\put(56,61){\line(1,0){12}}
\put(52,43){\line(0,1){14}}
\put(54,43){\line(1,1){16}}


\put(70,60){$x_{_{3}}$}

\put(20,63){\line(1,0){1.5}}
\put(8.5,50){\line(1,0){1.5}}
\put(23,50){\line(1,0){1.5}}

\put(40,63){\line(1,0){1.5}}
\put(28.5,50){\line(1,0){1.5}}

\put(48,55){\line(1,0){1.5}}

\put(20,70){\line(1,0){1.5}}
\put(16,32){\line(1,0){1.5}}

\put(8,60){\line(-1,0){5}}
\put(3,60){\line(0,-1){27}}
\put(3,33){\line(1,0){29}}
\put(32,38){\line(0,-1){5}}

\put(8,62){\line(-1,0){7}}
\put(1,62){\line(0,-1){33}}
\put(1,29){\line(1,0){51}}
\put(52,29){\line(0,1){8}}

\put(26,27){\line(1,0){1.5}}


\end{picture}

\end{diagram}

\clearpage

\clearpage
\subsection{
Simplifications of the Saw Blade Construction
}

\label{Section Simplifications}

\ecs

We show here that it is not necessary to make the $y_{ \xbs,i}$
contradictory
in a recursive construction, as in
Construction 
\ref{Construction Saw-Blade} (page 
\pageref{Construction Saw-Blade}).
It suffices to prevent them to be true.

We discuss three, much simplified, Saw Blade constructions.

Thus, we fully use here the conceptual difference of $x_{1}$ and $x_{2},$
as alluded to at the beginning of
Section \ref{Section Saw-Blades} (page \pageref{Section Saw-Blades}).

Note, however, that the back of each saw blade ``hides'' a Yablo
construction. The separate treatment of the teeth illustrates the
conceptual difference, but it cannot escape blurring it again in
the back of the blade.

First, we discuss some ``false'' simplifications which do not work.
\subsubsection{
``Simplifications'' that Will Not Work
}

\label{Section NotWork}

We try to simplify here the Saw Blade construction.
Throughout, we consider formulas of pure conjunctions.

We start with a Yablo Cell, but try to continue otherwise.

So we have $x_{0} \xcP x_{1} \xcP x_{2},$ $x_{0} \xcP x_{2}.$ So $x_{0}+$
is impossible. We now try to treat $x_{0}-.$
We see in Construction 
\ref{Construction Simple-1} (page 
\pageref{Construction Simple-1})
that appending $x_{2} \xch_{ \xCL }y_{2}$ may take care of the necessary
contradiction at $x_{2},$
see Diagram \ref{Diagram Sawblade-2} (page \pageref{Diagram Sawblade-2}),
lower part.
When we try to do the same at $x_{1},$ i.e. some $x_{1} \xch_{ \xCL
}x_{3},$ we solve again
the necessary contradiction at $x_{1},$ but run into a problem with
$x_{0}+,$
as $x_{1}$ is an $ \xco.$ So $x_{3}$ has to be contradictory. If we
continue
$x_{3} \xch_{ \xCL }x_{4} \xch_{ \xCL }x_{5}$ etc., this will not work, as
we may set all such $x_{i}-,$
and have a model. In abstract terms, we only procrastinate the same
problem
without solving anything.
Of course, we could append after some time new
Yablo Cells, as in the saw blade construction, but this is cheating,
as the ``true'' construction begins only later.
This shows the difference between ``knee'' and ``foot'', it works with
a foot, but not a knee. (Recall that all $x_{i},i>1$ are both foot and
knee.)

Suppose we add not only $x_{1} \xch_{ \xCL }x_{3},$ but also $x_{0} \xch_{
\xCL }x_{3},$ then we solve
$x_{0}+,$ but $x_{0}-$ is not solved.

Working with cells of the type (2.1) in
Example \ref{Example Simple-Cells} (page \pageref{Example Simple-Cells})
will lead to similar problems.

Consequently, any attempt to use a ``pipeline'', avoiding infinite
branching, is doomed:

Instead of
$x_{0} \xcP x_{1},$ $x_{0} \xcP x_{2},$  \Xl. etc. we construct a
``pipeline'' of $x_{i}',$ with
$x_{0} \xcP x'_{1},$ $x_{1} \xcP x'_{2},$ etc, $x'_{1} \xcp x'_{2} \xcp
x'_{3} \Xl.,$ and $x'_{1} \xcp x_{1},$ $x'_{2} \xcp x_{2},$ etc.
or similarly, to have infinitely many contradictions for paths
from $x_{0}.$

As this is a set of classical formulas, this cannot achieve inconsistency,
see Fact \ref{Fact Infinite} (page \pageref{Fact Infinite}).
\subsubsection{
Real Simplifications
}

\bcs

$\hspace{0.01em}$


\label{Construction Simple-1}

 \xEh

 \xDH

Take ONE saw blade $ \xbs,$ and attach (after closing under transitivity)
at all $y_{ \xbs,i}$ a SINGLE Yablo Cell $y_{ \xbs,i} \xcP u_{ \xbs,i}
\xcP v_{ \xbs,i},$ $y_{ \xbs,i} \xcP v_{ \xbs,i}.$
We call this the decoration, it is not involved in closure under
transitivity.

 \xEh

 \xDH

Any node $z$ in the saw blade (back or tooth) cannot be $z+,$ this leads
to a contradiction:

If $z=x_{i}$ (in the back):

Let $x_{i}+:$
Take $x_{i} \xcP x_{i+1}$ (any $x_{j},$ $i<j$ would do), if $x_{i+1} \xcP
r,$ then by transitivity,
$x_{i} \xcP r,$ so we have a contradiction.

If $z=y_{i}$ (a tooth):

$y_{i}+$ is contradictory by the ``decoration'' appended to $y_{i}.$

 \xDH

Any $x_{i}-$ $(x_{i}$ in the back, as a matter of fact, $x_{0}-$ would
suffice)
is impossible:

Consider any $x_{i} \xcP r,$ then $r+$ is impossible, as we just saw.

Note: there are no arrows from the back of the blade to the decoration.

 \xEj

 \xDH

We can simplify even further. The only thing we need about the $y_{i}$ is
that
they cannot be $+.$ Instead of decorating them with a Yablo Cell, any
contradiction will do, the simplest one is $y_{i}=y'_{i} \xcu \xCN
y'_{i}.$ Even
just one $y' $ s.t. $y_{i}=y' \xcu \xCN y' $ for all $y_{i}$ would do. (Or
a constant
FALSE.)

See Diagram \ref{Diagram Sawblade-2} (page \pageref{Diagram Sawblade-2}),
lower part.

Formally, we set

$x_{0}:= \xcU \{ \xCN x_{i}:i>0\} \xcu \xcU \{ \xCN y_{i}:i \xcg 0\},$

for $j>0:$

$x_{j}:= \xcU \{ \xCN x_{i}:i>j\} \xcu \xcU \{ \xCN y_{i}:i \xcg j-1\},$

and

$y_{j}:=y'_{j} \xcu \xCN y'_{j}.$

 \xDH

In a further step, we see that the $y_{i}$ (and thus the $y'_{i})$ need
not
be different from each other, one $y$ and one $y' $ suffice.

Thus, we set $x_{j}:= \xcU \{ \xCN x_{i}:i>j\} \xcu \xCN y,$ $y:=y' \xcu
\xCN y'.$

(Intuitively, the cells are arranged in a circle, with $y$ at the center,
and
$y' $ ``sticking out''. We might call this a ``curled saw blade''.

 \xDH
When we throw away the $y_{j}$ altogether, we have Yablo's construction.
this works, as we have the essential part in the $x_{i}$'s, and used the
$y_{j}' s$
only as a sort of scaffolding.

 \xEj

\ecs

\br

$\hspace{0.01em}$


\label{Remark Simple-2}

It seems difficult to conceptually simplify even further, as
Section 
\ref{Section Saw-Blades-Discussion} (page 
\pageref{Section Saw-Blades-Discussion})
shows basically the need for the
construction of the single Saw Blades.
We have to do something about the teeth, and above
Construction 
\ref{Construction Simple-1} (page 
\pageref{Construction Simple-1}), in particular cases (2) and (3) are
simple solutions.

\er

The construction is robust, as the following easy remarks show
(see also Section 
\ref{Section Procrastination} (page 
\pageref{Section Procrastination})):

 \xEh
 \xDH
Suppose we have ``gaps'' in the closure under transitivity, so, e.g. not
all $x_{0} \xcP x_{i}$ exist, they always exist only for $i>n.$ (And all
other $x_{k} \xcP x_{l}$
exist.) Then $x_{0}$ is still contradictory. Proof: Suppose $x_{0}+,$ then
we have the contradiction $x_{0} \xcP x_{n} \xcP x_{n+1}$ and $x_{0} \xcP
x_{n+1}.$ Suppose $x_{0}-,$ let $x_{0} \xcP x_{i}.$
As $x_{i}$ is unaffected, $x_{i}+$ is impossible.

 \xDH
Not only $x_{0}$ has gaps, but other $x_{i},$ too. Let again $x_{n}$ be an
upper bound for
the gaps. As above, we see that $x_{0}+,$ but also all $x_{i}+$ are
impossible.
If $x_{0} \xcP x_{i},$ as $x_{i}+$ is impossible, $x_{0}-$ is impossible.

 \xDH
$x_{0}$ has unboundedly often gaps, the other $x_{i}$ are not affected.
Thus, for $i \xEd 0,$ $x_{i}+$ and $x_{i}-$ are impossible. Thus, $x_{0}-$
is impossible, as
all $x_{i}+$ are, and $x_{0}+$ is, as all $x_{i}-$ are.

 \xEj

See also Section 
\ref{Section Procrastination} (page 
\pageref{Section Procrastination}).

$ \xCO $

$ \xCO $


\clearpage
\section{
The ``Right'' Level of Abstraction
}

\label{Section AbstractionLevel}
\subsection{
Introduction
}

The basic elements in Yablo's construction are negative arrows, from which
Yablo Cells are built. We show here that we can build arbitrarily complex
structures equivalent to negative arrows, and as a matter of fact to any
propositional logical operator. This suggests that the right level of
abstraction to consider more general Yablo-like structures is not the
level of single arrows, but rather of paths.

Note that  \cite{BS17} and  \cite{Wal23} also work with paths
in graphs.
\subsection{
Expressing Logical Operators by Combinations of Yablo Cells
}

\label{Section NegationCells}

\br

$\hspace{0.01em}$


\label{Remark Negation}

In the diagrams in
Diagram \ref{Diagram Negation} (page \pageref{Diagram Negation})  the
inner path $x-y-z$ is barred by $x-z$, and the outer path contains 3 or 2
negations.

Note that in all diagrams
Diagram \ref{Diagram Negation} (page \pageref{Diagram Negation}), and
Diagram \ref{Diagram TRUE} (page \pageref{Diagram TRUE}), upper part,
$y$ will always be FALSE, and as $x= \xCN y \xcu \xCN z,$ $ \xCN y$ is
TRUE and will not be
considered, the value of $x$ depends only on the branch through $z.$

E.g., in
Diagram \ref{Diagram Negation} (page \pageref{Diagram Negation}), upper part,
the path $z-x$ has uneven length, so the diagram describes negation, in
Diagram \ref{Diagram Negation} (page \pageref{Diagram Negation}), lower part,
the path $z' -z-x$ has even length, so the diagram describes identity.

\er

\bd

$\hspace{0.01em}$


\label{Definition Negation-Type}

We define negations and their type.

See Diagram \ref{Diagram Negation} (page \pageref{Diagram Negation}).

A negation (diagram) is a diagram all of whose arrows are of the $ \xcP $
kind.

We define negation diagrams, or, simply, negations, and their type.

 \xEI

 \xDH
An arrow $x \xcP x' $ is a negation of type 0.

A negation of the type
$x \xcH z,$ $x \xcH y \xcH z,$ $y \xcH y' \xcH z$
(see Diagram 
\ref{Diagram Negation} (page 
\pageref{Diagram Negation}), upper left)
is of type $n+1$ iff

all negations inside are of type $ \xck n,$ and at least one
negation inside is of type $n.$

 \xEJ

Similarly, we may define the type of an identity via the type of the
negations it is composed of - see
Diagram \ref{Diagram Negation} (page \pageref{Diagram Negation}),
lower part.

This all shows that we may blur the basic structure almost ad libitum,
making a characterisation difficult.
(See here Diagram \ref{Diagram TRUE} (page \pageref{Diagram TRUE}), lower part,
etc. too.)

\ed

\br

$\hspace{0.01em}$


\label{Remark Bedeutung}

Thus, it seems very difficult to describe Yablo type diagrams on the
level of single arrows. It seems to be the wrong level of
abstraction.

\er

We now give some examples.
\subsection{
Some Examples
}

\be

$\hspace{0.01em}$


\label{Example Logic}

See Diagram \ref{Diagram Negation} (page \pageref{Diagram Negation}), and
Diagram \ref{Diagram TRUE} (page \pageref{Diagram TRUE}).

 \xEh
 \xDH Negation (1)

$x= \xCN y \xcu \xCN z,$ $y= \xCN z \xcu \xCN y',$ $y' = \xCN z,$ so $
\xCN y=z \xco y',$ and
$x=(z \xco y') \xcu \xCN z$ $=$ $(z \xco \xCN z) \xcu \xCN z$ $=$ $ \xCN
z.$

 \xDH Negation (2)

$x= \xCN y \xcu \xCN z,$ $y= \xCN y'' \xcu \xCN y',$ $y' = \xCN z,$ $y''
=z,$ so $y= \xCN z \xcu \xCN y',$ and
continue as above.

 \xDH Identity

$x= \xCN y \xcu \xCN z,$ $y= \xCN z \xcu \xCN z',$ $z= \xCN z',$ so $
\xCN y=z \xco z',$ and $x=(\xCN z' \xco z') \xcu z' =z'.$

 \xDH TRUE

$x= \xCN y \xcu \xCN z,$ $y= \xCN y' \xcu \xCN z,$ $y' = \xCN z,$ $z= \xCN
z' \xcu \xCN z'',$ $z' = \xCN z''.$
Thus, $z=z'' \xcu \xCN z'',$ $ \xCN z=z'' \xco \xCN z'' =TRUE,$ $y'
=TRUE,$ $y=FALSE \xcu TRUE=FALSE,$ and
$x=TRUE \xcu TRUE.$

 \xDH $z' \xcu \xCN u$

$x= \xCN z \xcu \xCN y \xcu \xCN u,$ $y= \xCN z \xcu \xCN u' \xcu \xCN u,$
$z= \xCN z',$ $u' = \xCN u,$ so
$y=z' \xcu u \xcu \xCN u,$ $ \xCN y= \xCN z' \xco \xCN u \xco u=TRUE,$ and
$x=z' \xcu TRUE \xcu \xCN u=z' \xcu \xCN u.$

 \xEj

\ee

.

\clearpage

\begin{diagram}

Propositional Formulas 1

\label{Diagram Negation}
\index{Diagram Negation}

\unitlength0.7mm
\begin{picture}(200,210)(0,0)

\put(0,0){\line(1,0){200}}
\put(0,0){\line(0,1){210}}
\put(200,210){\line(-1,0){200}}
\put(200,210){\line(0,-1){210}}

\put(5,205){Lines represent negated downward arrows, except for the }
\put(5,199){line y"-z, which stands for a positive downward arrow. }

\put(5,190){{\rm\bf Diagram Negation }}

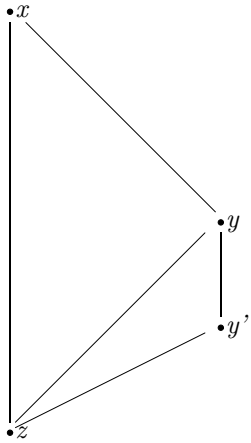
\put(10,180){\circle*{1}}
\put(50,140){\circle*{1}}
\put(50,120){\circle*{1}}
\put(10,100){\circle*{1}}

\put(11,179){x}
\put(51,139){y}
\put(51,119){y'}
\put(11,99){z}

\put(10,102){\line(0,1){76}}
\put(50,122){\line(0,1){16}}
\put(11,102){\line(1,1){36}}
\put(11,101){\line(2,1){36}}
\put(49,142){\line(-1,1){36}}

\put(85,180){We can achieve this using a diamond, too, see}
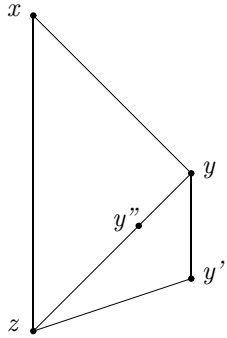
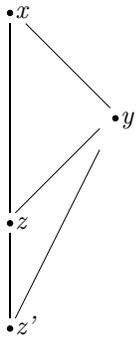
\put(85,175){Section \ref{Section Diamonds}}

\put(90,155){\circle*{1}}
\put(120,125){\circle*{1}}
\put(110,115){\circle*{1}}
\put(120,105){\circle*{1}}
\put(90,95){\circle*{1}}

\put(85,155){x}
\put(122,125){y}
\put(122,105){y'}
\put(105,115){y"}
\put(85,95){z}

\put(90,95){\line(0,1){60}}
\put(90,95){\line(1,1){20}}
\put(90,95){\line(3,1){30}}
\put(120,105){\line(0,1){20}}
\put(110,115){\line(1,1){10}}
\put(120,125){\line(-1,1){30}}

\put(5,75){{\rm\bf Diagram Identity }}

\put(10,65){\circle*{1}}
\put(30,45){\circle*{1}}
\put(10,25){\circle*{1}}
\put(10,5){\circle*{1}}

\put(11,64){x}
\put(31,44){y}
\put(11,24){z}
\put(11,4){z'}

\put(10,27){\line(0,1){36}}
\put(11,27){\line(1,1){16}}
\put(29,47){\line(-1,1){16}}
\put(10,7){\line(0,1){16}}
\put(11,7){\line(1,2){16}}

\end{picture}

\end{diagram}

\clearpage

\clearpage

\begin{diagram}

Propositional Formulas 2

\label{Diagram TRUE}
\index{Diagram TRUE}

\unitlength0.7mm
\begin{picture}(200,210)(0,0)

\put(0,0){\line(1,0){200}}
\put(0,0){\line(0,1){210}}
\put(200,210){\line(-1,0){200}}
\put(200,210){\line(0,-1){210}}

\put(50,200){{\rm\bf Diagram TRUE }}
\put(50,190){{(Lines represent negated downward }}
\put(50,186){{arrows) }}

\put(10,200){\circle*{1}}
\put(30,180){\circle*{1}}
\put(30,170){\circle*{1}}
\put(10,160){\circle*{1}}
\put(30,140){\circle*{1}}
\put(10,120){\circle*{1}}

\put(15,200){x}
\put(35,180){y}
\put(35,170){y'}
\put(17,160){z}
\put(35,140){z'}
\put(15,120){z"}

\put(10,120){\line(0,1){80}}
\put(30,170){\line(0,1){10}}
\put(30,180){\line(-1,-1){20}}
\put(30,180){\line(-1,1){20}}
\put(30,140){\line(-1,-1){20}}
\put(30,140){\line(-1,1){20}}
\put(30,170){\line(-2,-1){20}}

\put(100,100){{\rm\bf Diagram $z' \xcu \xCN u$ }}
\put(100,90){{(Lines represent negated downward arrows) }}

\put(50,100){\circle*{1}}
\put(50,80){\circle*{1}}
\put(50,60){\circle*{1}}
\put(90,40){\circle*{1}}
\put(10,40){\circle*{1}}
\put(10,20){\circle*{1}}

\put(50,100){\line(0,-1){40}}
\put(50,100){\line(2,-3){40}}
\put(50,100){\line(-2,-3){40}}

\put(10,40){\line(0,-1){20}}
\put(10,40){\line(1,1){40}}

\put(90,40){\line(-1,1){40}}
\put(90,40){\line(-2,1){40}}

\put(10,20){\line(2,3){40}}

\put(55,100){x}
\put(95,40){u}
\put(55,80){y}

\put(55,60){u'}

\put(4,40){z}
\put(4,20){z'}

\end{picture}

\end{diagram}

\clearpage

\subsection{
Paths Instead of Arrows
}

\label{Section Paths-Arrows}

We concluded above that the level of arrows might be the wrong level of
abstraction.
We change perspective, and consider paths, i.e. sequences of arrows,
sometimes neglecting branching points.

We assume all formulas attached to nodes are as in Yablo's paper, i.e.
pure conjunctions of positive or negative
nodes.

\bcd

$\hspace{0.01em}$


\label{Condition Yablo}

Note that we speak only about special cases here: The formulas
are of the type $ \xcU \xCN,$ we use Yablo cells for contradictions,
and every arrow (negative path below) leads directly to a new Yablo cell,
see Section \ref{Section Prereq} (page \pageref{Section Prereq}).

(We use an extension of the language beyond $ \xcU \xCN $ for
simplifications
(``cutting unwanted branches''),
so, strictly speaking,
our result is NOT an equivalence result. In addition to $ \xCN x_{i}$ as
components of $ \xcU,$ we also admit formulas of the type $(x_{i} \xco
\xCN x_{i})$ or the
constant TRUE.)

Suppose a graph $ \xbG $ contains an infinite set of points $x_{i}:$ $i
\xbe \xbo,$
which can be connected by negative paths $ \xbs_{i,j}:x_{i} \Xl x_{j}$
for all $i<j,$
then we can interpret this graph so that $x_{0}$ (as a matter of fact, all
$x_{i})$
cannot have a truth value.

Thus, the graph is directed, the arrows are not necessarily labelled
(labels $+$ or -), but labelling can be done such that for each pair
of nodes $x,y,$ where $y$ is a successor of $x,$ there is a negative path
from $x$ to $y.$
(See Definition 
\ref{Definition Path-Value} (page 
\pageref{Definition Path-Value})).

Remark: The condition ``negative paths'' is not trivial. Suppose that
$ \xbs:x \Xl y$ and $ \xbs':y \Xl z$ are both negative, and $ \xbs \xDO
\xbs' $ is the only path
from $x$ to $z,$ then the condition is obviously false. The author does
not
know how to characterize graphs which satisfy the condition. Some kind
of ``richness'' will probably have to hold.

\ecd

\subparagraph{
Proof
}

$\hspace{0.01em}$


We saw in Section 
\ref{Section Saw-Blades-Discussion} (page 
\pageref{Section Saw-Blades-Discussion})
that we can generate
from the prerequisites a structure isomorphic to Yablo's structure.

Conversely:
Let $y$ be any node in $ \xbG,$ if $y \xcP y' $ is not part of the $
\xbs_{i,j},$ then
interpret $y \xcP y' $ by $(y' \xco \xCN y').$ Thus, $y \xcP y' $ has no
influence on the
truth of $y.$ In particular, the branching points inside the $ \xbs_{i,j}$
disappear, and the $ \xbs_{i,j}$ become trivial. The $x_{i}$ are not
influenced
by anything apart from the $x_{j},$ $j>i.$

Thus, the construction of the
$x_{i},$ $x_{j},$ $ \xbs_{i,j}$ is equivalent to Yablo's structure.

\bd

$\hspace{0.01em}$


\label{Definition Sloppy}

Our argument above was a bit sloppy. We can make it more precise.

If $x$ is contradictory, i.e. both $x+$ and x- lead to a contradiction,
we have $M(x)= \xCQ,$ and $M(\xCN x)= \xCQ,$
with $M(.)$ the set of models.

Still, the following definition seems reasonable:

$M(x \xcu \xbf):=M(x) \xcs M(\xbf),$ and $M(x \xco \xbf):=M(x) \xcv M(
\xbf).$

In above proof, we used $M(x \xcu TRUE)=M(x).$
\clearpage
\section{
Possible Contradiction Cells
}

\label{Section ContradictionCells}
\subsection{
Introduction
}

\ed

We consider alternative basic contradiction cells, similar to Yablo cells.
We saw in Example 
\ref{Example Simple-Cells} (page 
\pageref{Example Simple-Cells})
that Yablo cells are the simplest contradictions suitable
to a contradictory construction (i.e. without possible truth values).
We now examine more complicated contradiction cells, the only
interesting one is the ``diamond'', but this cell fails to fit into our
construction principle - at least in our logical framework, see
Section \ref{Section Nested-Diamonds} (page \pageref{Section Nested-Diamonds}).

The Yablo cell has the form $x \xcP y \xcP z,$ $x \xcP z,$ whereas the
diamond
has the form $x \xcP y \xcP z,$ $x \xcP y' \xcp z.$
The diamond is more complicated, as it contains two branching points, $y$
and $y',$
before the contradiction is complete, and we have to ``synchronize''
what happens at those branching points.
\subsection{
Prerequisites
}

In Yablo's construction, and our Saw Blades, the basic contradiction had
the form $x \xcP y \xcP z,$ $x \xcP z$ (Yablo Cell, YC).
We generalize this now, as illustrated in
Diagram 
\ref{Diagram YC2} (page 
\pageref{Diagram YC2}), $x_{0}$ corresponds to $x,$ $a_{2}$ to $y,$
$b_{2}$ to $z.$
The arrow $x \xcP y$ is replaced by a perhaps more complicated path via
$a_{1},$ $x \xcP z$
by a perhaps more complicated path via $b_{2},$ but we still want $a_{2}$
and $b_{2}$ origins
of new contradictions in the case $x_{0}-.$ Thus, both $ \xbs:x_{0} \Xl
a_{2}$ and $ \xbs':x_{0} \Xl b_{2}$
have to be negative paths.

The most interesting modification is
of the type of Case (2.4.2)
in Remark 
\ref{Remark General-YC} (page 
\pageref{Remark General-YC}), where the contradiction originating
at $x_{0}$
is via an additional point, $z.$ We call this a ``Diamond'', and examine
this case in more detail in Section 
\ref{Section Diamonds} (page 
\pageref{Section Diamonds}), and
Section \ref{Section Nested-Diamonds} (page \pageref{Section Nested-Diamonds}).

In the Diamond case, we need a positive arrow (equivalently two negative
arrows without branching in between), either from $a_{2}$ to $z,$ or from
$b_{2}$ to $z.$

In all cases, $ \xbs $ from $x_{0}$ to $a_{2},$ and $ \xbs' $ from
$x_{0}$ to $b_{2}$ has to be
negative, so $x_{0}-$ results in $a_{2}+$ and $b_{2}+.$ Omitting $a_{1}$
and $b_{1},$
and adding a direct arrow $a_{2} \xcP b_{2}$ results in the original Yablo
Cell.

We continue to work with pure conjunctions, now of positive or negative
propositional variables. (We will be forced to look at more complicated
formulas in Remark 
\ref{Remark Rhombus-Basic} (page 
\pageref{Remark Rhombus-Basic}).)
\subsection{
Various Types
}

\label{Section VariousTypes}

\br

$\hspace{0.01em}$


\label{Remark General-YC}

See Diagram \ref{Diagram YC2} (page \pageref{Diagram YC2}).

We examine here the different modifications of Yablo's basic construction.

Note that e.g. $x_{0} \Xl a_{1}$ need not be an arrow, it may be a longer
path.

The contradiction need not be e.g. $a_{2} \Xl b_{2},$ but they might take
a ``detour''
$a_{2} \Xl z$ and $b_{2} \Xl z.$ (The cases like $a_{2} \Xl z \Xl b_{2}$
are already covered by $a_{2}..b_{2}.)$
So $z$ will be the ``point of conflict''. We will call these cases
``detour to $z$'' (and the construction ``diamonds'').

We do not break e.g. $x_{0} \Xl a_{1}$ further down, this
suffices for our analysis.

The situation:
 \xEh
 \xDH $x_{0}+$ must be contradictory, the contradiction formed with
(perhaps part of)
$ \xbs $ and $ \xbs',$

 \xDH both $a_{2}+$ and $b_{2}+$ must be contradictory if $x_{0}-.$
 \xEj

Recall that the paths $x_{0} \Xl a_{1} \Xl a_{2}$ and $x_{0} \Xl b_{1} \Xl
b_{2}$ are both negative
(as $a_{2}$ and $b_{2}$ must be positive, if $x_{0}$ is negative).

We look at the following cases of building the contradiction at $x_{0},$
e.g. Case (1.1) means that we have the contradiction between
$x_{0} \Xl a_{1} \Xl b_{1}$ and $x_{0} \Xl b_{1},$
Case (1.4.1) means that we have a contradiction between
$x_{0} \Xl a_{1} \Xl z$ and $x_{0} \Xl b_{1}..z$.

 \xEh

 \xDH $a_{1}$ as additional branching point
 \xEh
 \xDH from $a_{1}$ to $b_{1}$
 \xDH from $a_{1}$ to $b_{2}$
 \xDH from $a_{1}$ to $b_{3}$
 \xDH from $a_{1}$ to $z$ and
 \xEh
 \xDH from $b_{1}$ to $z,$ $b_{1}$ as additional branching point
 \xDH from $b_{2}$ to $z,$ $b_{2}$ as additional branching point
 \xDH from $b_{3}$ to $z,$ $b_{3}$ as additional branching point
 \xEj
 \xEj

 \xDH $a_{2}$ (again $b_{i}$ as additional branching points)
 \xEh
 \xDH from $a_{2}$ to $b_{1}$
 \xDH from $a_{2}$ to $b_{2}$
 \xDH from $a_{2}$ to $b_{3}$
 \xDH from $a_{2}$ to $z$ and
 \xEh
 \xDH from $b_{1}$ to $z$
 \xDH from $b_{2}$ to $z$
 \xDH from $b_{3}$ to $z$
 \xEj
 \xEj

 \xDH $a_{3}$ (again $a_{3}$ and $b_{i}$ as additional branching points)
 \xEh
 \xDH from $a_{3}$ to $b_{1}$
 \xDH from $a_{3}$ to $b_{2}$
 \xDH from $a_{3}$ to $b_{3}$
 \xDH from $a_{3}$ to $z$ and
 \xEh
 \xDH from $b_{1}$ to $z$
 \xDH from $b_{2}$ to $z$
 \xDH from $b_{3}$ to $z$
 \xEj

 \xEj

 \xEj

We will examine the cases now.

 \xEh

 \xDH Case (1):

In all cases, we have to branch at $a_{1}$ with $a_{1}+ \xcu,$
as we need to have a contradiction for $x_{0}+,$ so we have $x_{0} \xcp
a_{1}.$
But then, for $x_{0}-,$ we branch at $a_{1}- \xco.$ Now, one branch will
lead to $a_{2},$
the other to $b_{1},b_{2},b_{3},$ or $z.$

Cases (1.1) and (1.2):

Suppose $x_{0}-,$ and $a_{1}$ is chosen at $x_{0}.$ As $a_{1}-,$ suppose
$b_{1}$ is chosen at $a_{1}.$
By prerequisite (we need a contradiction),
$x_{0} \Xl a_{1} \Xl b_{1}$ contradicts $x_{0} \Xl b_{1},$ moreover $x_{0}
\Xl b_{1}..b_{2}$ is
negative, so $x_{0} \Xl a_{1} \Xl b_{1} \Xl b_{2}$ is positive, and the
chosen path misses
$a_{2}$ and makes $b_{2}$ negative, contradicting (2) above in ``The
situation''.
(Thus, appending the same type of construction at $b_{2}$ again,
etc. will result in an escape path $x_{0}-,$ $b_{2}-,$ etc.)

Case (1.3) and (1.4):

Consider $x_{0}-,$ chose at $x_{0}$ $a_{1}$ and chose at $a_{1}$ $b_{3}$
(or $z).$ Then
$x \Xl a_{1} \Xl b_{3}$ or $x \Xl a_{1} \Xl z$ will not meet $a_{2}$ nor
$b_{2}.$

 \xDH Case (2):

Cases (2.1) and (2.2) are (equivalent to) a Yablo Cell.

Cases (2.3), (2.4.2), and (2.4.3) are equivalent, and discussed below
(``Diamond''),
see Section \ref{Section Diamonds} (page \pageref{Section Diamonds}).
Note that, e.g. in case (2.4.2), to have a contradiction, we need
$val(\xbr) \xEd val(\xbr').$

Cases (2.4.1) and (1.4.2) are symmetrical.

 \xDH Case (3):

Similar considerations as for Case (2) apply.

 \xEj

\clearpage

\begin{diagram}

General Contradiction

\label{Diagram YC2}
\index{Diagram YC2}

\unitlength1mm
\begin{picture}(140,190)(0,0)

\put(0,0){\line(1,0){140}}
\put(0,0){\line(0,1){190}}
\put(140,190){\line(-1,0){140}}
\put(140,190){\line(0,-1){190}}

\put(20,180){See Remark \ref{Remark General-YC}}

\put(15,80){$x_{0}+$}

\put(40,90){$a_1$}
\put(60,100){$a_2$}
\put(80,110){$a_3$}

\put(100,80){$z$}

\put(20,160){Recall that $val(x_0..a_2)=val(x_0..b_2)=-$}
\put(20,156){but $val(\xbr) \xEd val(\xbr')$}

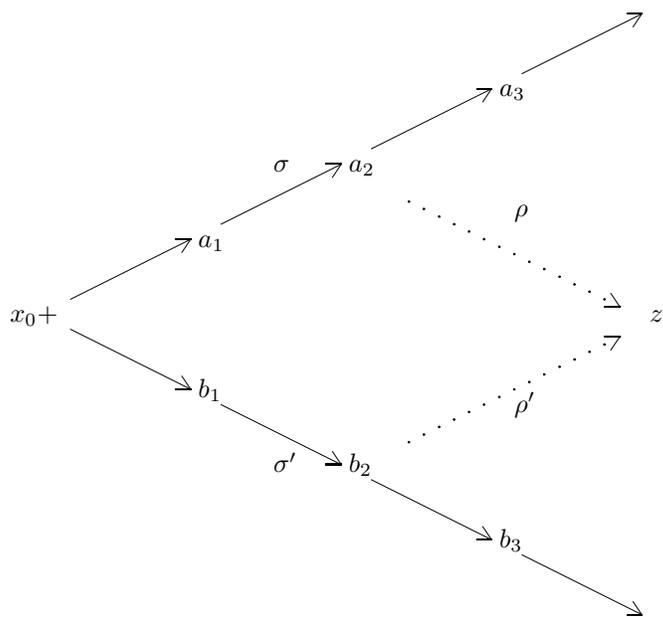
\put(23,83){\line(2,1){16}}
\put(39,91){\line(-1,0){2}}
\put(39,91){\line(-2,-3){1.2}}

\put(43,93){\line(2,1){16}}
\put(59,101){\line(-1,0){2}}
\put(59,101){\line(-2,-3){1.2}}

\put(63,103){\line(2,1){16}}
\put(79,111){\line(-1,0){2}}
\put(79,111){\line(-2,-3){1.2}}

\put(83,113){\line(2,1){16}}
\put(99,121){\line(-1,0){2}}
\put(99,121){\line(-2,-3){1.2}}

\put(50,100){$\xbs $}

\put(40,70){$b_1$}
\put(60,60){$b_2$}
\put(80,50){$b_3$}

\put(23,79){\line(2,-1){16}}
\put(39,71){\line(-1,0){2}}
\put(39,71){\line(-2,3){1.2}}

\put(43,69){\line(2,-1){16}}
\put(59,61){\line(-1,0){2}}
\put(59,61){\line(-2,3){1.2}}

\put(63,59){\line(2,-1){16}}
\put(79,51){\line(-1,0){2}}
\put(79,51){\line(-2,3){1.2}}

\put(83,49){\line(2,-1){16}}
\put(99,41){\line(-1,0){2}}
\put(99,41){\line(-2,3){1.2}}

\put(82,94){$\xbr $}
\put(82,68){$\xbr'$}

\put(50,60){$\xbs'$}

\multiput(68,64)(2,1){15}{\circle*{0.2}}
\put(96,78){\line(-1,0){2}}
\put(96,78){\line(-2,-3){1.2}}

\multiput(68,96)(2,-1){15}{\circle*{0.2}}
\put(96,82){\line(-1,0){2}}
\put(96,82){\line(-2,3){1.2}}

\put(20,24){The solid arrows stand for positive or negative paths}
\put(20,20){For the dotted arrows see Remark 
\ref{Remark General-YC}, Case (2.4.2) -}
\put(20,16){they form part of the Diamond contradiction: }
\put(20,12){$x_0..b_1..b_2..z$, $x_0..a_1..a_2..z$.}

\end{picture}

\end{diagram}

\vspace{10mm}

\clearpage

\subsection{
Diamonds
}

\label{Section Diamonds}

\er

We first give some simple examples as
``warming up exercises'', before we turn in
Section 
\ref{Section Nested-Diamonds} (page 
\pageref{Section Nested-Diamonds})  to a proof that - at least in our
setting - Diamonds cannot replace Yablo Cells, see
Remark \ref{Remark Rhombus-Basic} (page \pageref{Remark Rhombus-Basic}).

Recall that in all examples below
$val(x_{0},x_{1})=val(x_{0},x_{2})=-.$

\br

$\hspace{0.01em}$


\label{Remark Rhombus-3a}

We may, of course, imitate a diamond $a \xcP b \xcP d,$ $a \xcP c \xcp d$
by a triangle
$a \xcP b \xcP d,$ $a \xcP c \xcp d,$ where the ``detour'' via $c$ is
integrated in the line
$a \xcP d.$ But this is ``cheating'', and would neglect the essential
property,
that each of $b$ and $c$ should have a new diamond attached to it.

\er

\br

$\hspace{0.01em}$


\label{Remark Explanation}

Case (1), see Diagram 
\ref{Diagram Rhombus-1.1} (page 
\pageref{Diagram Rhombus-1.1}).

$x_{0}$ and $x_{1}$ should be heads of contradiction cells.

Take e.g. Case (1.1). Here, $x_{1}$ is a conjunction of positive or
negative
occurences of $y_{1},$ $y_{2},$ $y,$ Thus, if $x_{0}+,$ $ \xCN x_{1}$ is a
disjunction of
(the negatives) of above occurences. So, it might e.g. be $y_{1} \xco \xCN
y_{2} \xco y,$
and we cannot be sure of going to $y,$ and thus we are not sure to have
a contradiction for the case $x_{0}+.$
In more detail:

If $x_{0}+,$ then $x_{1}-,$ and let $x_{1}-$ be expressed by some logical
formula
$ \xbf = \xbf_{1} \xco  \Xl  \xco \xbf_{n},$ with each $ \xbf_{i}$ a
conjunction of letters.
We know that, for the contradiction to work for $x_{0},$ $val(x_{1},y)
\xEd val(x_{2},y)$
has to hold. Suppose $val(x_{2},y)=-,$ so val $(x_{1},y)=+.$
Then, in each conjunction, we must have $y.$ Take now the case $x_{0}-,$
then $x_{1}+$
and consequently $val(x_{1},y_{1}) \xEd val(x_{1},y_{2},y_{1})$ has to
hold. Fix $val(y_{2},y_{1}).$
The negation of $ \xbf $ is the disjunction of the (negated) conjunctions
of all choice functions in the conjunctions for $ \xbf.$ But one of
those conjunctions is just $ \xCN y \xcu  \Xl  \xcu \xCN y= \xCN y,$ it
says nothing about
$y_{1}$ or $y_{2},$ so we cannot force any coherence between
$val(x_{1},y_{1})$ and
$val(x_{1},y_{2},y_{1}).$

Note that this argument works for arbitrary infinite formulas, too.

The argument for Case (1.2) is similar.

Case (3) is more interesting, see
Diagram \ref{Diagram Rhombus-3.1} (page \pageref{Diagram Rhombus-3.1}).

We want $x_{0},$ $x_{1},$ and $x_{2}$ to be heads of contradiction cells.
We discuss Case (3.1), (3.2) is similar.

Suppose $val(x_{1},y)=-,$ $val(x_{2},y)=+$

We distinguish two possibilities:

Case A: $val(y,z)=+.$

Then $val(x_{1},z)=+,$ $val(x_{2},z)=-.$
If $x_{0}+,$ then $x_{1}- \xcO,$ $x_{2}- \xcO.$ Suppose $x_{1}$ chooses
$z,$ and $x_{2}$ chooses $y,$ then
we have no contradiction.

Case $B:$ $val(y,z)=-.$

Then $val(x_{1},z)=-,$ $val(x_{2},z)=+.$
Again, if $x_{1}$ chooses $z,$ and $x_{2}$ chooses $y,$ we have no
contradiction.

\clearpage

\begin{diagram}

Simple Diamonds 1

\label{Diagram Rhombus-1.1}
\index{Diagram Rhombus-1.1}

\unitlength1mm
\begin{picture}(140,190)(0,0)

\put(0,0){\line(1,0){140}}
\put(0,0){\line(0,1){190}}
\put(140,190){\line(-1,0){140}}
\put(140,190){\line(0,-1){190}}

\put(10,180){See Remark \ref{Remark Explanation}}

\put(10,175){Unmarked lines denote positive or negative upward pointing arrows.}


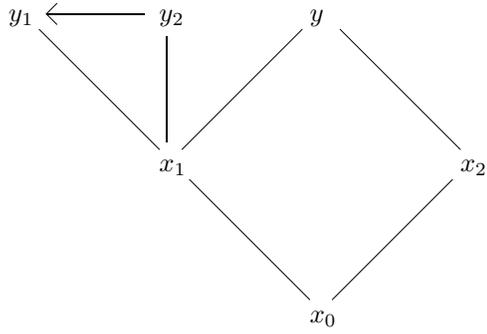
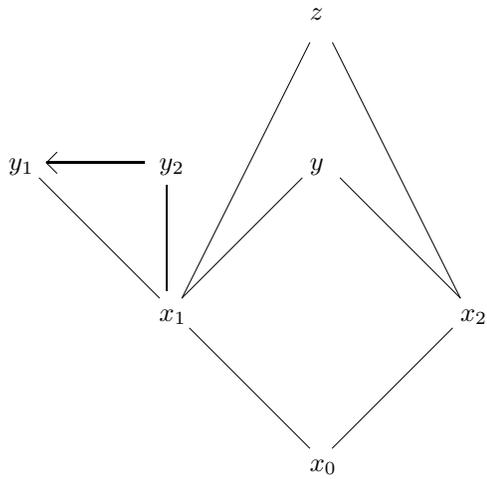
\put(100,140){Case (1.1)}
\put(60,120){$x_{0}$}
\put(40,140){$x_{1}$}
\put(80,140){$x_{2}$}
\put(20,160){$y_{1}$}
\put(60,160){$y$}
\put(40,160){$y_2$}

\put(60,123){\line(-1,1){16}}
\put(63,123){\line(1,1){16}}

\put(40,143){\line(-1,1){16}}
\put(43,143){\line(1,1){16}}

\put(80,143){\line(-1,1){16}}

\put(38,161){\line(-1,0){13}}
\put(41,158){\line(0,-1){14}}
\put(25,161){\line(1,1){1.4}}
\put(25,161){\line(1,-1){1.4}}


\put(100,40){Case (1.2)}

\put(60,20){$x_{0}$}
\put(40,40){$x_{1}$}
\put(80,40){$x_{2}$}
\put(20,60){$y_{1}$}
\put(60,60){$y$}
\put(60,80){$z$}
\put(40,60){$y_2$}

\put(60,23){\line(-1,1){16}}
\put(63,23){\line(1,1){16}}

\put(40,43){\line(-1,1){16}}
\put(43,43){\line(1,1){16}}

\put(80,43){\line(-1,1){16}}

\put(43,43){\line(1,2){17}}
\put(80,43){\line(-1,2){17}}

\put(38,61){\line(-1,0){13}}
\put(41,58){\line(0,-1){14}}
\put(25,61){\line(1,1){1.4}}
\put(25,61){\line(1,-1){1.4}}

\end{picture}

\end{diagram}

\vspace{10mm}

\clearpage

\clearpage

\begin{diagram}

Simple Diamonds 2

\label{Diagram Rhombus-3.1}
\index{Diagram Rhombus-3.1}

\unitlength1mm
\begin{picture}(140,190)(0,0)

\put(0,0){\line(1,0){140}}
\put(0,0){\line(0,1){190}}
\put(140,190){\line(-1,0){140}}
\put(140,190){\line(0,-1){190}}

\put(10,180){See Remark \ref{Remark Explanation} }
\put(10,175){Unmarked lines denote positive or negative upward pointing arrows.}


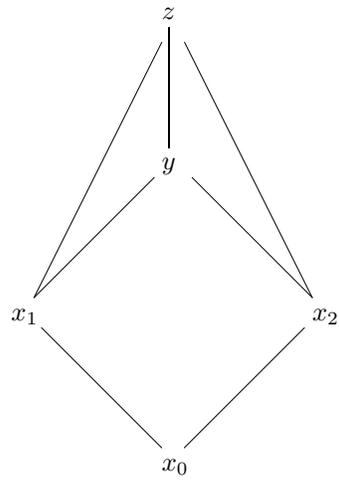
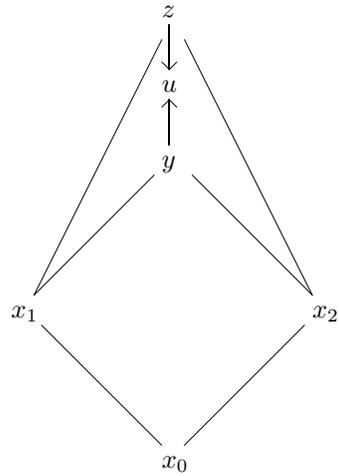
\put(100,30){Case (3.2)}

\put(60,10){$x_{0}$}
\put(40,30){$x_{1}$}
\put(80,30){$x_{2}$}
\put(60,50){$y$}
\put(60,70){$z$}
\put(60,60){$u$}

\put(60,13){\line(-1,1){16}}
\put(63,13){\line(1,1){16}}

\put(43,33){\line(1,1){16}}

\put(80,33){\line(-1,1){16}}

\put(43,33){\line(1,2){17}}
\put(80,33){\line(-1,2){17}}

\put(61,53){\line(0,1){6}}
\put(61,63){\line(0,1){6}}

\put(61,59){\line(1,-1){1}}
\put(61,59){\line(-1,-1){1}}

\put(61,63){\line(1,1){1}}
\put(61,63){\line(-1,1){1}}


\put(100,120){Case (3.1)}

\put(60,100){$x_{0}$}
\put(40,120){$x_{1}$}
\put(80,120){$x_{2}$}
\put(60,140){$y$}
\put(60,160){$z$}

\put(60,103){\line(-1,1){16}}
\put(63,103){\line(1,1){16}}

\put(43,123){\line(1,1){16}}
\put(80,123){\line(-1,1){16}}

\put(43,123){\line(1,2){17}}
\put(80,123){\line(-1,2){17}}

\put(61,143){\line(0,1){16}}

\end{picture}

\end{diagram}

\vspace{10mm}

\clearpage

\clearpage
\subsubsection{
Nested Diamonds
}

\label{Section Nested-Diamonds}

\er

We present here some examples and problems of constructions with nested
diamonds.

\br

$\hspace{0.01em}$


\label{Remark Versuch}

We first try to imitate the Yablo construction using diamonds instead of
triangles.
See Diagram \ref{Diagram Versuch} (page \pageref{Diagram Versuch}),
for the moment both sides.

The basic construction is $x_{0}-x_{1,1}-x_{2},$ $x_{0}-x_{1,2}-x_{2}.$

To make the construction recursive, we have to add, among other things,
a new diamond at $x_{1,1},$ $x_{1,1}-x_{3,1}-x_{4}$ and
$x_{1,1}-x_{3,2}-x_{4},$ etc.

 \xEh
 \xDH
Consider first the left hand side of the diagram.

$x_{0}-x_{1,2}-x_{2}$ contradicts $x_{0}-x_{1,1}-x_{2},$ and
$x_{0}-x_{1,2}-x_{2}-x_{3,2}$ contradicts $x_{0}-x_{1,1}-x_{3,2}.$

$x_{0}-x_{1,1}-x_{3,1}-x_{4}$ and $x_{0}-x_{1,2}-x_{2}-x_{3,2}-x_{4}$ are
not contradictory, just as
$x_{0}-x_{1,1}-x_{3,1}-x_{5,2}$ and
$x_{0}-x_{1,2}-x_{2}-x_{3,2}-x_{4}-x_{5,2}$ are not contradictory.

More generally, the paths via $x_{i,1}$ to $x_{i+1}$ and via $x_{i,1}$ to
$x_{i+2,2}$ behave
the same way. It is also easy to see that $x_{i,1}-x_{i+2,2}-x_{i+3}$ and
$x_{i,1}-x_{i+2,1}-x_{i+3}$ behave differently with respect to
contradicting
the simple vertical path, e.g.
$x_{0}-x_{1,2}-x_{2}-x_{3,2}$ contradicts $x_{0}-x_{1,1}-x_{3,2},$ but
$x_{0}-x_{1,1}-x_{3,1}-x_{4}$ and $x_{0}-x_{1,2}-x_{2}-x_{3,2}-x_{4}$ are
not contradictory.

Consider now $x_{0} \xDL,$ then $x_{1,1} \xDN \xcO,$ and we have to make
all paths
to contradict the corresponding simple vertical path.
This works for $x_{1,1}-x_{2}$ and $x_{1,1}-x_{3,2},$ but not for
$x_{1,1}-x_{3,1}-x_{4}$ and
$x_{1,1}-x_{3,1}-x_{5,2}.$
As $x_{3,1} \xDL \xcU,$ we have a ``second chance'' to continue via
$x_{0}-x_{1,1}-x_{3,1}-x_{5,1},$ but $x_{5,1}$ is $ \xDN \xcO,$ so
$x_{5,1}-x_{6}$ and $x_{5,1}-x_{7,2}$ are ok,
but $x_{5,1}-x_{7,1}$ is not, and we have an escape path.

 \xDH
The construction on the right hand side of the diagram avoids this
problem,
as for $x_{0} \xDL,$ $x_{3,1} \xDL \xcU,$ and
$x_{0}-x_{1,1}-x_{3,1}-x_{6,2}$ contradicts
$x_{0}-x_{1,2}-x_{2}-x_{3,2}-x_{4}-x_{5}-x_{6,2}.$

But we created a different problem for the added diamond:
If $x_{1,1} \xDL,$ then $x_{3,1} \xDN,$ so not only
$x_{1,1}-x_{3,1}-x_{4}$ has to contradict $x_{1,1}-x_{3,2}-x_{4},$ but
also
$x_{1,1}-x_{3,1}-x_{6,2}$ has to contradict
$x_{1,1}-x_{3,2}-x_{4}-x_{5}-x_{6,2},$
which does not work.

In addition, we have a triangle contradiction in $x_{3,1} \xcp x_{4} \xcP
x_{5} \xcP x_{6,2},$ $x_{3,1} \xcP x_{6,2}.$

In summary: both attempts fail, but, of course, a different approach
might still work.

 \xDH
Note that the roles of $ \xcU $ and $ \xcO $ are exchanged when we go from
$x_{0}$ to
$x_{1,1}$ and $x_{1,2}.$

In Yablo's construction, this presents no problem, as the structure seen
from
$x_{0}$ and any other $x_{i}$ is the same.

This is not the case in our example.

 \xEj

\clearpage

\begin{diagram}

Nested Diamonds 1

\label{Diagram Versuch}
\index{Diagram Versuch}

\unitlength0.7mm
\begin{picture}(200,220)(0,0)

\put(0,0){\line(1,0){200}}
\put(0,0){\line(0,1){220}}
\put(200,220){\line(-1,0){200}}
\put(200,220){\line(0,-1){220}}

\multiput(40,10)(0,20){11}{\circle*{1.5}}

\put(40,10){\line(0,1){200}}

\put(42,10){$x_0$}
\put(42,50){$x_2$}
\put(42,90){$x_4$}
\put(42,130){$x_6$}
\put(42,170){$x_8$}
\put(42,210){$x_{10}$}

\put(42,30){$x_{1,2}$}
\put(42,70){$x_{3,2}$}
\put(42,110){$x_{5,2}$}
\put(42,150){$x_{7,2}$}
\put(42,190){$x_{9,2}$}

\multiput(42,20)(0,20){10}{$-$}

\multiput(20,30)(0,40){5}{\circle*{1.5}}

\put(10,30){$x_{1,1}$}
\put(10,70){$x_{3,1}$}
\put(10,110){$x_{5,1}$}
\put(10,150){$x_{7,1}$}
\put(10,190){$x_{9,1}$}

\put(20,30){\line(0,1){160}}

\multiput(15,50)(0,40){4}{$-$}

\put(40,10){\line(-1,1){20}}

\put(30,15){$-$}

\multiput(20,30)(0,40){5}{\line(1,1){20}}
\multiput(20,30)(0,40){4}{\line(1,2){20}}

\multiput(30,37)(0,40){5}{$+$}
\multiput(30,57)(0,40){4}{$-$}


\multiput(100,10)(0,20){8}{\circle*{1.5}}
\put(100,10){\line(0,1){140}}

\put(102,10){$x_0$}
\put(102,50){$x_2$}
\put(102,90){$x_4$}
\put(102,130){$x_{6,2}$}

\put(102,30){$x_{1,2}$}
\put(102,70){$x_{3,2}$}
\put(102,110){$x_{5}$}
\put(102,150){$x_{7}$}

\multiput(102,20)(0,20){7}{$-$}

\put(80,30){\circle*{1.5}}
\put(80,70){\circle*{1.5}}
\put(80,130){\circle*{1.5}}

\put(80,30){\line(0,1){100}}

\put(70,30){$x_{1,1}$}
\put(70,70){$x_{3,1}$}
\put(70,130){$x_{6,1}$}

\put(100,10){\line(-1,1){20}}

\put(90,15){$-$}

\multiput(80,30)(0,40){2}{\line(1,1){20}}
\put(80,130){\line(1,1){20}}

\put(80,30){\line(1,2){20}}
\put(80,70){\line(1,3){20}}

\put(90,37){$+$}
\put(90,77){$+$}
\put(90,137){$+$}
\put(90,57){$-$}
\put(90,112){$-$}

\put(75,50){$-$}
\put(75,100){$-$}

\end{picture}

\end{diagram}

\vspace{10mm}

\clearpage

\er

\br

$\hspace{0.01em}$


\label{Remark Rhombus-Basic}

See Diagram 
\ref{Diagram Rhombus-Basic} (page 
\pageref{Diagram Rhombus-Basic}), and
in more detail Diagram 
\ref{Diagram Matrix-4} (page 
\pageref{Diagram Matrix-4}).

We neglect the triangles like $x_{1,1}-x_{3,2}-x_{2},$ $x_{1,1}-x_{2}.$
(See (5) below for
a comment.)

The diagrams show the basic recursive construction, but additional
Diamonds
still have to be added on top of the existing Diamonds.

To make the Diagram 
\ref{Diagram Rhombus-Basic} (page 
\pageref{Diagram Rhombus-Basic})
more readable, we noted some points several
times $(x_{2},$ $x_{4,1},$ $x_{4,2}),$ they are connected by vertical fat
lines.

We have 7 diamonds, $(x_{0},$ $x_{1,1},$ $x_{2},$ $x_{1,2}),$ $(x_{1,1},$
$x_{3,1},$ $x_{4,1},$ $x_{3,2}),$
$(x_{1,2},$ $x_{3,3},$ $x_{4,2},$ $x_{3,4}),$
$(x_{3,1},$ $x_{5,1},$ $x_{6,1},$ $x_{5,2}),$ $(x_{3,2},$ $x_{5,3},$
$x_{6,2},$ $x_{5,4}),$
$(x_{3,3},$ $x_{5,5},$ $x_{6,3},$ $x_{5,6}),$ $(x_{3,4},$ $x_{5,7},$
$x_{6,4},$ $x_{5,8}).$
They are drawn using thick lines.

In all diamonds, the bottom and the upper right lines are supposed to
be negative, the upper left line is meant positive for a contradiction.

This diagram, or some modification of it seems the right way to
to use diamonds instead of triangles for the contradictions recursively,
as the triangles are used recursively in Yablo's construction.

We will discuss here problems with this construction.

 \xEh

 \xDH ``Synchronization''

 \xEh
 \xDH We have a conflict between the diamonds starting at $x_{1,1}$ and
$x_{1,2}$ and
the diamond starting at $x_{0}.$

If $x_{0}+ \xcU,$ then $x_{1,1}- \xcO $ and $x_{1,2}- \xcO.$ As the
choices at $x_{1,1}$ and $x_{1,2}$ are
independent, any branch $x_{0}-x_{1,1}-x_{3,1}-x_{2},$
$x_{0}-x_{1,1}-x_{3,2}-x_{2},$ $x_{0}-x_{1,1}-x_{2}$
combined with any branch $x_{0}-x_{1,2}-x_{3,3}-x_{2},$
$x_{0}-x_{1,2}-x_{3,4}-x_{2},$ $x_{0}-x_{1,2}-x_{2}$
must be conflicting, thus, given $x_{0}-x_{1,1}-x_{2}$ is negative, all
branches
on the left must be negative, likewise, all branches on the right must
positive.

However, if $x_{1,1}$ is positive, the diamond $x_{1,1}-x_{3,1}-x_{2},$
$x_{1,1}-x_{3,2}-x_{2}$
has to be contradictory, so not both branches may be negative.

 \xDH A solution is to ``synchronise'' the choices at $x_{1,1}$ and
$x_{1,2}$ which
can be done e.g. by the formula

$x_{0}$ $=$ $ \xCN x_{1,1} \xcu \xCN x_{1,2}$ $ \xcu $ $[(x_{3,1} \xcu
x_{3,4}) \xco (x_{3,2} \xcu x_{3,3}) \xco (x_{2} \xcu \xCN x_{2})],$ and

$ \xCN x_{0}$ $=$ $x_{1,1} \xco x_{1,2}$ $ \xco $ $[ \xCN (x_{3,1} \xcu
x_{3,4}) \xcu \xCN (x_{3,2} \xcu x_{3,3}) \xcu \xCN (x_{2} \xcu \xCN
x_{2})]$ $=$
$x_{1,1} \xco x_{1,2}$ $ \xco $ $[(\xCN x_{3,1} \xco \xCN x_{3,4}) \xcu (
\xCN x_{3,2} \xco \xCN x_{3,3}) \xcu (\xCN x_{2} \xco x_{2})].$

(This poses a problem, as validity of $ \xCN x_{0}$ might involve validity
of $ \xCN x_{2},$
so we construct an Escape path. This can be transformed into a semantical
argument, so we cannot avoid considering $ \xCN x_{2}.)$

This formula is of a different type than $ \xcU \xCN.$ In addition,
arrows from $x_{0}$
to the $x_{3,i}$ are missing.

Even if we simplify further and consider only the paths
$x_{0}-x_{1,i}-x_{3,j}-x_{2}$
(i.e. do not consider the direct paths $x_{0}-x_{1,i}-x_{2}),$
we have the formula

$x_{0}$ $=$ $ \xCN x_{1,1} \xcu \xCN x_{1,2}$ $ \xcu $ $[(x_{3,1} \xcu
x_{3,4}) \xco (x_{3,2} \xcu x_{3,3})],$ and

$ \xCN x_{0}$ $=$ $x_{1,1} \xco x_{1,2}$ $ \xco $ $[ \xCN (x_{3,1} \xcu
x_{3,4}) \xcu \xCN (x_{3,2} \xcu x_{3,3})]$ $=$
$x_{1,1} \xco x_{1,2}$ $ \xco $ $[(\xCN x_{3,1} \xco \xCN x_{3,4}) \xcu (
\xCN x_{3,2} \xco \xCN x_{3,3})],$
so we can make $x_{1,1}$ and $x_{1,2}$ false, and chose $(\xCN x_{3,1}$
and $ \xCN x_{3,3})$ or
$(\xCN x_{3,4}$ and $ \xCN x_{3,2}),$ and have $x_{2}-$ or $x_{2}+,$ the
first possibility
results again in an escape path.

Note: Adding the information about $x_{1,1}$ and $x_{1,2}$ does not help:
$ \xCN x_{1,1}=x_{3,1} \xco x_{3,2},$ $ \xCN x_{1,2}=x_{3,3} \xco
x_{3,4},$ so $ \xCN x_{1,1} \xcu \xCN x_{1,2}$ $=$
$(x_{3,1} \xcu x_{3,3})$ $ \xco $ $(x_{3,1} \xcu x_{3,4})$ $ \xco $
$(x_{3,2} \xcu x_{3,3})$ $ \xco $ $(x_{3,2} \xcu x_{3,4})$ still
leaves the possibility $(\xCN x_{3,1}$ and $ \xCN x_{3,3})$ by making the
last disjunct
true.

 \xDH Summary:

This concerns, of course, the type of constructions we analysed here:
direct nesting of diamonds, and direct contradictions.

The use of diamonds, together with synchronisation problems, seems to
present
unsolvable problems, not only for the simple type of formulas considered
here, but also to arbitrary formulas.

There is no synchronisation problem in Yablo's construction, as there is
only one choice there, and ``meeting'' is avoided by the universal
quantifier. The latter is impossible with diamonds as we need a
contradiction
on both sides, thus both sides cannot be uniform.

 \xEj

 \xDH
If we modify the construction, and consider for $x_{0}+$ not the
direct lines $x_{1,1}-x_{3,1}$ and $x_{1,1}-x_{3,2},$ but go higher up
$(x_{3,i}$ are all $+,$ $ \xcu)$
e.g. lines $x_{4,1}-x_{2},$ then this will not work either, as one way
to go to $x_{4,1}$ is positive, the other is negative, so the paths
$x_{1,1}-x_{3,1}-x_{4,1}-x_{2}$ and $x_{1,1}-x_{3,2}-x_{4,1}-x_{2}$ cannot
both be positive
as required.

 \xDH
If we try to go higher up, e.g. to $x_{5,3}$ and $x_{5,4}$ we have again
an $ \xco,$
and need a contradiction in the diamond $(x_{3,2},$ $x_{5,3},$ $x_{6,2},$
$x_{5,4}),$
so we are in the same situation as with the diamond $(x_{1,1},$ $x_{3,1},$
$x_{4,1},$ $x_{3,2}).$

 \xDH
(If we create after $x_{3,1}$ or $x_{3,2}$ new branching points, these
points
will be $- \xcO $ either for $x_{0}+$ or for $x_{1,1}+,$ so one case will
fail.)

 \xDH
Back to the triangles like $x_{1,1}-x_{3,2}-x_{2},$ $x_{1,1}-x_{2}.$

The probably simplest way would be to create an intermediate point, say
$x_{1,1,a}$
between $x_{1,1}$ and $x_{2}$ (here with $x_{1,1} \xcP x_{1,1,a} \xcP
x_{2},$ on the right hand side
of the diagram with $x_{1,2} \xcP x_{1,2,a} \xcp x_{2})$ to create a
situation similar
to the diamonds discussed here.

We did not investigate this any further, as it is outside our framework.

 \xDH
We give a tentative, abstract appoach to the problem, allowing
a mixture of different types of contradictions, in
Remark \ref{Remark Abstraction} (page \pageref{Remark Abstraction}).

 \xEj

\er

\br

$\hspace{0.01em}$


\label{Remark Abstraction}

In a general construction, we may use a mixture of different contradictory
cells, Yablo triangles, diamonds (if we can solve the associate problems),
etc.

Thus, it is probably not the right level of abstraction to consider
full CCS's (see Definition 
\ref{Definition Terminology} (page 
\pageref{Definition Terminology})),
as they might be composed of CC's of different types.

We will, however, keep the following:
 \xEh
 \xDH the formulas will be basically of the form $ \xcU \xCN,$
 \xDH $x_{0}$ is the head of a contradiction, so $x_{0}+$ is impossible,
 \xDH every path from $x_{0}$ will lead to (the head of) a new
contradiction, and this
path will be negative, so we make $x_{0}-$ impossible, too.
 \xEj

Consequently:

 \xEI
 \xDH if $x_{0}+,$ and a negative path $ \xbs $ leads from
$x_{0}$ to the head $x_{1}$ of a new contradiction,
then $x_{1}- \xcO,$ so we have to form new negative paths $ \xbs' $ from
$x_{0}$ to contradict
all $x_{1}' $ reached from $x_{1}$ to build contradictions starting at
$x_{1}.$
These new $x_{1}' $ will have the same properties as $x_{1},$
leading to an infinite (in depth and width) construction as detailed in
Construction \ref{Construction Details} (page \pageref{Construction Details}),
Diagram \ref{Diagram X1} (page \pageref{Diagram X1}), and
Diagram \ref{Diagram X2} (page \pageref{Diagram X2}).

Of course, we might have to add new properties, like ``synchronisation''.

 \xEJ

In Section 
\ref{Section Procrastination} (page 
\pageref{Section Procrastination})  we discuss additional problems
which arise if
we do not construct contradictions as soon as possible, preventing local
solutions.

\clearpage

\begin{diagram}

Nested Diamonds 2

\label{Diagram Rhombus-Basic}
\index{Diagram Rhombus-Basic}

\unitlength0.7mm
\begin{picture}(200,210)(0,0)

\put(0,0){\line(1,0){200}}
\put(0,0){\line(0,1){210}}
\put(200,210){\line(-1,0){200}}
\put(200,210){\line(0,-1){210}}

\put(20,190){See Remark \ref{Remark Rhombus-Basic}}

\thicklines

\put(80,10){\circle*{1}}
\put(78,7){\footnotesize $x_0$}
\put(80,10){\line(-2,3){40}}
\put(80,10){\line(2,3){40}}

\put(40,70){\circle*{1}}
\put(33,70){\footnotesize $x_{1,1}$}
\put(120,70){\circle*{1}}
\put(122,70){\footnotesize $x_{1,2}$}
\put(40,70){\line(-2,3){20}}
\put(40,70){\line(2,3){20}}
\put(40,70){\line(2,1){40}}
\put(120,70){\line(-2,3){20}}
\put(120,70){\line(2,3){20}}
\put(120,70){\line(-2,1){40}}

\put(80,90){\circle*{1}}
\put(78,86){\footnotesize $x_{2}$}

\put(20,100){\circle*{1}}
\put(13,100){\footnotesize $x_{3,1}$}
\put(60,100){\circle*{1}}
\put(62,100){\footnotesize $x_{3,2}$}

\put(20,100){\line(-1,4){10}}
\put(20,100){\line(0,1){40}}
\put(20,100){\line(1,1){20}}
\put(60,100){\line(1,4){10}}
\put(60,100){\line(0,1){40}}
\put(60,100){\line(-1,1){20}}
\put(100,100){\line(-1,4){10}}
\put(100,100){\line(0,1){40}}
\put(100,100){\line(1,1){20}}
\put(140,100){\line(1,4){10}}
\put(140,100){\line(0,1){40}}
\put(140,100){\line(-1,1){20}}

\put(100,100){\circle*{1}}
\put(93,100){\footnotesize $x_{3,3}$}
\put(140,100){\circle*{1}}
\put(142,100){\footnotesize $x_{3,4}$}

\put(80,120){\circle*{1}}
\put(82,122){\footnotesize $x_{2}$}
\put(80,140){\circle*{1}}
\put(78,142){\footnotesize $x_{2}$}

\put(40,120){\circle*{1}}
\put(33,120){\footnotesize $x_{4,1}$}
\put(120,120){\circle*{1}}
\put(122,120){\footnotesize $x_{4,2}$}

\put(10,140){\circle*{1}}
\put(3,140){\footnotesize $x_{5,1}$}
\put(20,140){\circle*{1}}
\put(22,140){\footnotesize $x_{5,2}$}
\put(60,140){\circle*{1}}
\put(53,140){\footnotesize $x_{5,3}$}
\put(70,140){\circle*{1}}
\put(64,139){\footnotesize $x_{5,4}$}

\put(90,140){\circle*{1}}
\put(84,140){\footnotesize $x_{5,5}$}
\put(100,140){\circle*{1}}
\put(102,140){\footnotesize $x_{5,6}$}
\put(140,140){\circle*{1}}
\put(133,140){\footnotesize $x_{5,7}$}
\put(150,140){\circle*{1}}
\put(144,139){\footnotesize $x_{5,8}$}

\put(10,140){\line(1,2){5}}
\put(20,140){\line(-1,2){5}}
\put(60,140){\line(1,2){5}}
\put(70,140){\line(-1,2){5}}

\put(90,140){\line(1,2){5}}
\put(100,140){\line(-1,2){5}}
\put(140,140){\line(1,2){5}}
\put(150,140){\line(-1,2){5}}

\put(15,150){\circle*{1}}
\put(13,153){\footnotesize $x_{6,1}$}
\put(65,150){\circle*{1}}
\put(63,153){\footnotesize $x_{6,2}$}
\put(95,150){\circle*{1}}
\put(93,153){\footnotesize $x_{6,3}$}
\put(145,150){\circle*{1}}
\put(143,153){\footnotesize $x_{6,4}$}

\put(40,160){\circle*{1}}
\put(37,163){\footnotesize $x_{4,1}$}
\put(120,160){\circle*{1}}
\put(117,163){\footnotesize $x_{4,2}$}

\put(39.9,120){\line(0,1){40}}
\put(40.1,120){\line(0,1){40}}
\put(119.9,120){\line(0,1){40}}
\put(120.1,120){\line(0,1){40}}

\put(79.9,90){\line(0,1){50}}
\put(80.1,90){\line(0,1){50}}

\thinlines

\put(20,100){\line(3,1){60}}
\put(60,100){\line(1,1){20}}
\put(100,100){\line(-1,1){20}}
\put(140,100){\line(-3,1){60}}

\put(40,120){\line(2,1){40}}
\put(120,120){\line(-2,1){40}}

\put(10,140){\line(3,2){30}}
\put(20,140){\line(1,1){20}}
\put(60,140){\line(-1,1){20}}
\put(70,140){\line(-3,2){30}}

\put(90,140){\line(3,2){30}}
\put(100,140){\line(1,1){20}}
\put(140,140){\line(-1,1){20}}
\put(150,140){\line(-3,2){30}}

\put(15,150){\line(5,2){25}}
\put(65,150){\line(-5,2){25}}
\put(95,150){\line(5,2){25}}
\put(145,150){\line(-5,2){25}}

\end{picture}

\end{diagram}

\vspace{10mm}

\clearpage

\clearpage

\begin{diagram}

Nested Diamonds 2, Details

\label{Diagram Matrix-4}
\index{Diagram Matrix-4}

\unitlength0.5mm
\begin{picture}(280,220)(0,0)

\put(0,0){\line(1,0){280}}
\put(0,0){\line(0,1){220}}
\put(280,220){\line(-1,0){280}}
\put(280,220){\line(0,-1){220}}

\put(20,205){Example for synchronisation, see Remark \ref{Remark Rhombus-Basic}}

\put(80,20){\circle*{1}}
\put(78,15){\footnotesize $x_0$}

\put(80,20){\line(-1,2){40}}
\put(80,20){\line(1,2){40}}

\put(62,60){$-$}
\put(95,60){$-$}

\put(40,100){\circle*{1}}
\put(120,100){\circle*{1}}

\put(30,100){\footnotesize $x_{1,1}$}
\put(122,100){\footnotesize $x_{1,2}$}

\put(40,100){\line(-1,1){20}}
\put(40,100){\line(0,1){20}}
\put(40,100){\line(1,2){40}}
\put(120,100){\line(1,1){20}}
\put(120,100){\line(0,1){20}}
\put(120,100){\line(-1,2){40}}

\put(25,110){$-$}
\put(36,110){$-$}
\put(122,110){$-$}
\put(133,110){$-$}

\put(20,120){\circle*{1}}
\put(40,120){\circle*{1}}
\put(120,120){\circle*{1}}
\put(140,120){\circle*{1}}

\put(10,120){\footnotesize $x_{3,1}$}
\put(30,120){\footnotesize $x_{3,2}$}
\put(121,120){\footnotesize $x_{3,3}$}
\put(142,120){\footnotesize $x_{3,4}$}

\put(20,120){\line(1,1){60}}
\put(40,120){\line(2,3){40}}
\put(140,120){\line(-1,1){60}}
\put(120,120){\line(-2,3){40}}

\multiput(20,120)(-1,1){13}{\circle*{0.3}}
\multiput(40,120)(-1,1){13}{\circle*{0.3}}
\multiput(120,120)(1,1){13}{\circle*{0.3}}
\multiput(140,120)(1,1){13}{\circle*{0.3}}

\put(46,150){$+$}
\put(57,150){$-$}
\put(69,150){$+$}
\put(89,150){$-$}
\put(101,150){$+$}
\put(112,150){$-$}

\put(80,180){\circle*{1}}
\put(78,185){\footnotesize $x_2$}


\put(180,100){\circle*{1}}
\put(182,100){\footnotesize $x_{1,3}$}

\put(80,20){\line(5,4){100}}
\put(180,100){\line(-5,4){100}}
\multiput(180,100)(1,1){13}{\circle*{0.3}}
\multiput(180,100)(0,1){13}{\circle*{0.3}}

\put(123,60){$-$}
\put(125,150){$+$}

\put(10,5){Adding $x_{1,3}$ etc. does not help, as we also add the same problems
again.}

\end{picture}

\end{diagram}

\vspace{10mm}

\clearpage

\clearpage
\section{
Illustration of (Finite) Procrastination
}

\label{Section Procrastination}
\subsection{
Introduction
}

\er

This section is about the original Yablo construction and its
modifications, postponing contradictions.
Throughout, we work with Yablo's order (of the natural numbers).

In Yablo's construction, full transitivity guarantees contradictions.
If full transitivity is absent, we may still have contradictions,
provided we have ``enough'' transitivity. We discuss this here.

The most important problem seems to be to find a notation which is easy to
read
and write, expressing the alternation between $ \xcA $ (or $ \xcU)$ and $
\xcE $ (or $ \xcO,$
arbitrary choice). We do this as follows:
 \xEh
 \xDH vertical lines beginning at some $x$ illustrate all $y>x,$ they
express
$ \xcA y.y>x.$
 \xDH diagonal lines to the left starting at $x$ and going to $x' $
express the choice
of $x' $ for $x,$ i.e. $ \xcE x'.x' >x.$
 \xEj

First, we consider the meaning of nested quantifiers in our context.
This is obvious, but writing it down helps.

We did not elaborate all details, but hope that these notes suffice
to illustrate the important cases and considerations.

Note that the situation also suggests - in hindsight - a game theoretic
approach, as taken e.g. in  \cite{Wal23}.

\br

$\hspace{0.01em}$


\label{Remark Koeln}

(Trivial)

Consider the order in the graph as in Yablo's construction.
Quantifiers are mostly restricted, and range only over all elements
bigger than some reference point: $ \xcA y.y>x$  \Xl., etc.

In the case of mixed quantifiers, essentially the last one decides the
meaning. If the chain of quantifiers ends by $ \xcA,$ then the property
holds
for an end segment of the graph, if it ends by $ \xcE,$ then we only know
that
there is still some element where it holds. In the case $ \xcA \xcE,$
this means
that there are cofinally many elements where is holds.

Some examples for illustration:

 \xEh
 \xDH
$ \xcE x \xcA y:$ all $y$ from $x+1$ onward have a certain property, this
is a full
end segment of the natural numbers.

In particular, this set has a non-empty
intersection with any cofinal, i.e. infinite, sequence.

 \xDH
$ \xcA x \xcE y:$ this might e.g. be the next prime after $x,$ the next
prime $+1,$ etc.

Thus, the intersection with a cofinal sequence might be empty, but not
the intersection with an end segment.

 \xDH
$ \xcE \xcA \xcE:$ this is like $ \xcA \xcE,$ starting at a certain
value, again, the
intersection with a cofinal sequence may be empty.

 \xDH
$ \xcA x \xcE y \xcA z:$ this is like $ \xcE y \xcA y,$ for all $x$ there
is a start, and from
then onward, all $z$ are concerned, i.e. we have a full end segment.

 \xDH
Thus, $ \xcA \xcE $ $ \xcs $ $ \xcE \xcA $ $ \xEd \xCQ $

 \xDH
Thus, $ \xcA \xcE $ $ \xcs $ $ \xcA \xcE $ $= \xCQ $ is possible, see
above the primes and primes $+$ 1
cases.

 \xDH
$ \xcA \xcE:$ it helps to interpret this as $f(a)$ for all a in $ \xcA,$
similarly $ \xcA \xcE \xcA \xcE $ as $f_{1}(a),$ $f_{2}(b)$ etc.
Chose first a, then $f_{1}(a),$ then for all $b>f_{1}(a)$ $f_{2}(b).$

 \xDH
In the following diagrams, the diagonal lines describe the choice
functions
$f(x)$ corresponding to $ \xcE.$

 \xEj

\er

In Yablo like diagrams it sometimes helps to consider ``local transitivity
'',
which we introduce now.

\bd

$\hspace{0.01em}$


\label{Definition Local}

$ \xBc x_{n},x_{n' } \xBe $ is locally transitive iff

(1) $x_{n} \xcP x_{n' }$ and

(2) $ \xcA x_{n'' }$ s.t. $x_{n' } \xcP x_{n'' },$ $x_{n} \xcP x_{n'' }$
too.

\ed

\bfa

$\hspace{0.01em}$


\label{Fact Local}

 \xEh
 \xDH If $ \xBc x_{n},x_{n' } \xBe $ is locally transitive for some $x_{n' },$
then
$x_{n}+$ is
impossible.

 \xDH Let $m<n,$ and

$ \xcA $ $x_{n}$ s.t. $x_{m} \xcP x_{n}$ $ \xcE $ $x_{n' }$ s.t.
$ \xBc x_{n},x_{n' } \xBe $ is locally transitive,

then $x_{m}-$ is impossible.

 \xDH Thus, we need (1) and (2) for $x_{0}.$
 \xEj
\subsection{
Discussion of More Complicated Cases
}

\efa

\be

$\hspace{0.01em}$


\label{Example Escape-2}

See Diagram \ref{Diagram Escape} (page \pageref{Diagram Escape}).

We consider the case $x_{i}+.$

 \xEh
 \xDH
Notation:

The nodes are ordered by the natural order of $ \xbo.$
We have here partial transitivity of $ \xcP $ only, with $x \xcP x+1$ as
basis.

$x \xDc y$ means: $x<y$ and $x \xcP y.$

 \xDH
The vertical lines show the nodes above the bottom node.

A full circle at a node expresses that there is an arrow $ \xcP $ from the
bottom
node to this node, an empty node expresses that there is no such arrow.
No circle expresses that the case is left open (but usually there will be
an
arrow) - see the text below.

The dotted slanted lines
indicate choices. E.g., if $x_{1}-,$ then there is some $x_{2},$ s.t.
$x_{1} \xDc x_{2}$ (and thus
$x_{2}+).$

 \xDH
We want to show that $x_{0}+$ is contradictory.
 \xEh

 \xDH
In Yablo's construction, we find a contradiction going through $x_{0}+1,$
as a matter of fact, through all $x_{1}>x_{0}.$ Thus, chose any
$x_{1}>x_{0},$ then
$x_{1}- \xcO,$ so we have to find for all $x_{1} \xcP x_{2}$ an arrow
$x_{0} \xcP x_{2},$ which,
of course, exists.

 \xDH
In our diagram, we try to find a contradiction through $x_{1}.$

Suppose $x_{0} \xDc x_{1},$ so for $x_{0}+,$ $x_{1}- \xcO,$ and for all
$x_{2}$ s.t. $x_{1} \xDc x_{2}$ we look
for an arrow $x_{0} \xcP x_{2},$ i.e. for $x_{0} \xDc x_{2},$ transitivity
for $ \xDc.$

This does not hold, see the empty circle at $x_{2}$ in the $x_{0}$ column.
So $x_{1}$ cannot be a suitable knee in the contradiction for $x_{0}.$

 \xDH
We try to ``mend'' $x_{2}$ now, and make it contradictory itself.

Suppose $x_{2} \xDc x_{4} \xDc x_{5},$ then we want $x_{2} \xDc x_{5},$
but $x_{5}$ has an empty circle in
the $x_{2}$ column, so this does not hold.

Suppose $x_{2} \xDc x_{6} \xDc x_{7}.$ As noted by the full circle at
$x_{7}$ in the $x_{2}$ column,
$x_{2} \xDc x_{7}.$

So, for this one choice of $x_{7}$ at $x_{6},$ we have a contradiction at
$x_{2}.$
To have a contradiction for $x_{2}$ in all choices at $x_{6},$ we need: if
$x_{6} \xDc y,$ then $x_{2} \xDc y$ has to hold.
This is, of course transitivity for $x_{2} \xDc x_{6} \xDc y.$

(As indicated by the full circle for $x_{7}$ in the $x_{0}$ column, we
also
have a contradiction between $x_{0} \xDc x_{7}$ and the path of length 4
$x_{0} \xDc x_{1} \xDc x_{2} \xDc x_{6} \xDc x_{7}$ - for this path, and
neither for all $y$ such that
$x_{6} \xDc y,$ nor for all $y$ such that $x_{1} \xDc y.)$

Note that this is not necessary for all $x_{6}$ with $x_{2} \xDc x_{6},$
one such $x_{6}$ suffices
to show that $x_{2}+$ is contradictory.

 \xDH
Suppose we do not have some $x_{1}$ with $x_{0} \xDc x_{1}$ s.t. for all
$x_{2}$ with $x_{1} \xDc x_{2}$
we also have $x_{0} \xDc x_{2},$ then we can look for another, better,
$x_{1},$ or try
to mend the old $x_{1}.$

Suppose we have $x_{0} \xDc x_{1} \xDc x_{2},$ but not $x_{0} \xDc x_{2}.$
We make now $x_{2}$ itself contradictory,
as described above in (3.3).

This repair possibility holds, of course, recursively.

If the procedure fails (for all $x_{1},$ $x_{6}$ etc.) we construct an
escape path
which shows that $x_{0}$ is not contradictory.

 \xDH
Apart from the remark in (3.3) above,
we do not discuss here more complicated contradictions
with paths longer than 1, as any contradicting path, e.g. $x_{0} \xDc
x_{1}' \xDc x_{1}'' \xDc x_{2}$
would involve a new ``OR'', here at $x_{1}',$ so we have new branchings,
and no control over meeting $x_{2}$ - unless the new branchings also meet
$x_{2}.$

 \xEj

 \xEj

\clearpage

\begin{diagram}

Procrastination

\label{Diagram Escape}
\index{Diagram Escape}

\unitlength1mm
\begin{picture}(140,190)(0,0)

\put(0,0){\line(1,0){140}}
\put(0,0){\line(0,1){190}}
\put(140,190){\line(-1,0){140}}
\put(140,190){\line(0,-1){190}}

\put(12,170){See the discussion in Example \ref{Example Escape-2} }

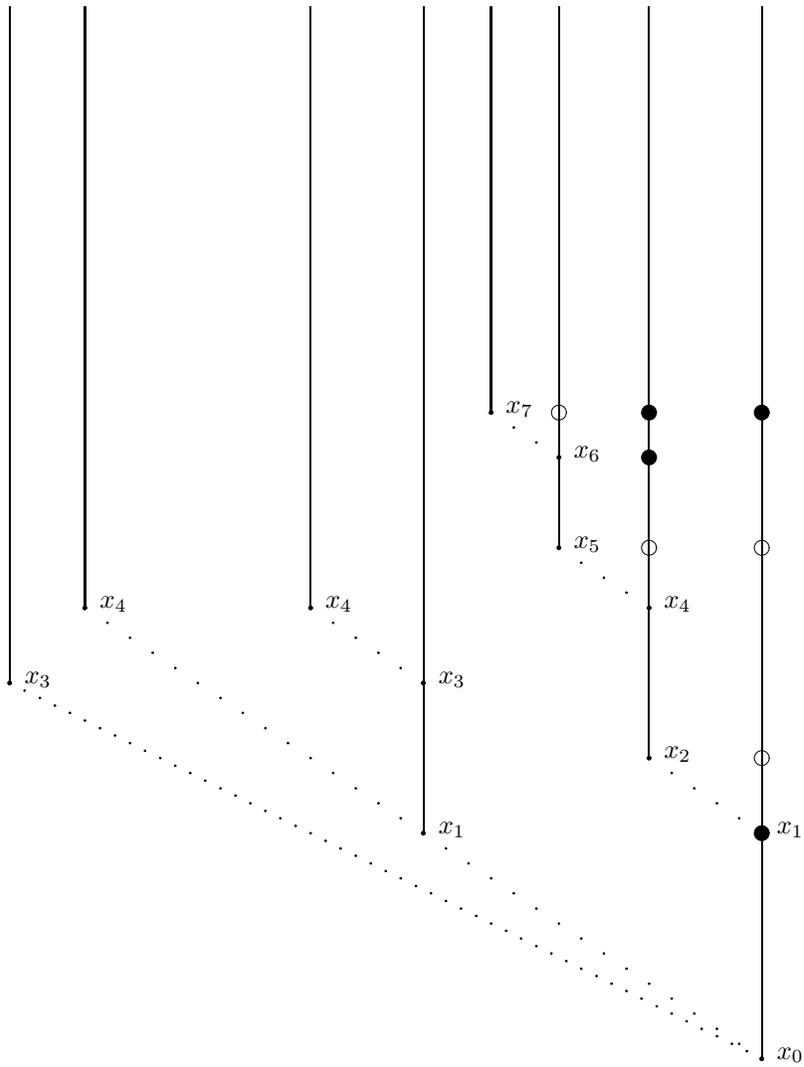
\put(10,70){\circle*{0.5}}
\put(12,70){$x_3$}
\put(10,70){\line(0,1){90}}

\put(20,80){\circle*{0.5}}
\put(22,80){$x_4$}
\put(20,80){\line(0,1){80}}

\put(50,80){\circle*{0.5}}
\put(52,80){$x_4$}
\put(50,80){\line(0,1){80}}

\put(65,50){\circle*{0.5}}
\put(67,50){$x_1$}
\put(65,50){\line(0,1){110}}

\put(65,70){\circle*{0.5}}
\put(67,70){$x_3$}
\put(74,106){\circle*{0.5}}
\put(76,106){$x_7$}
\put(74,106){\line(0,1){54}}

\put(83,88){\circle*{0.5}}
\put(85,88){$x_5$}
\put(83,88){\line(0,1){72}}

\put(83,100){\circle*{0.5}}
\put(85,100){$x_6$}
\put(95,60){\circle*{0.5}}
\put(97,60){$x_2$}
\put(95,60){\line(0,1){100}}

\put(95,80){\circle*{0.5}}
\put(97,80){$x_4$}
\put(110,20){\circle*{0.5}}
\put(112,20){$x_0$}
\put(110,20){\line(0,1){140}}

\put(110,50){\circle*{0.5}}
\put(112,50){$x_1$}

\multiput(110,20)(-2,1){50}{\circle*{0.2}}
\multiput(110,20)(-3,2){16}{\circle*{0.2}}

\multiput(65,50)(-3,2){16}{\circle*{0.2}}

\multiput(110,50)(-3,2){5}{\circle*{0.2}}
\multiput(65,70)(-3,2){5}{\circle*{0.2}}

\multiput(95,80)(-3,2){4}{\circle*{0.2}}

\multiput(83,100)(-3,2){3}{\circle*{0.2}}

\put(110,50){\circle*{2}}
\put(110,60){\circle{2}}
\put(110,88){\circle{2}}
\put(110,106){\circle*{2}}

\put(95,88){\circle{2}}
\put(95,100){\circle*{2}}
\put(95,106){\circle*{2}}

\put(83,106){\circle{2}}

\end{picture}

\end{diagram}

\vspace{10mm}

\clearpage

\clearpage
\section{
Some Additional Comments
}

\label{Section Additional}
\subsection{
Introduction
}

\ee

We make the truth value ``impossible'', noted here (X), formal, and
argue that the approach to nested diamonds
in Section 
\ref{Section Nested-Diamonds} (page 
\pageref{Section Nested-Diamonds})
is excessive, we can do with one half. But one half is essentially a Yablo
Cell,
and the other half does not work.

We summarize the results in
Section \ref{Section Branch} (page \pageref{Section Branch}).
\subsection{
The New Truth Value in Detail
}

We introduce a new truth value (X), with the meaning that, in the present
structure, $x=(X)$ $x$ can neither be true, nor false, all such
assumptions lead to a contradiction.

We define the operators $ \xCN,$ $ \xcu,$ $ \xco $ for (X) as follows:

$ \xCN (X)=(X)$

$(X) \xcu (X)=(X)$

$(X) \xco (X)=(X)$

$(X) \xcu TRUE=(X)$

$(X) \xcu FALSE=FALSE$

$(X) \xco TRUE=TRUE$

$(X) \xco FALSE=(X)$

Thus, in particular, $(X) \xcu \xCN (X)$ $=$ $(X) \xco \xCN (X)$ $=$ (X).

This generalizes to infinitary versions of $ \xcu $ and $ \xco,$ e.g., in
$ \xcU x_{i},$
we may eliminate $x_{i}=TRUE,$ and $x_{i}=FALSE$ makes $ \xcU x_{i}$
FALSE.
\subsection{
How Much Recursion Do We Need
}

\label{Section How-Much}

Thus, when we branch $x=x' \xcu x'',$ and want to make $x=(X)$ in the
final
construction, we need to make only $x' =(X)$ or $x'' =(X),$ so we have to
make
a recursive construction only in one branch. This was done already in
Construction 
\ref{Construction Simple-1} (page 
\pageref{Construction Simple-1}).

On the other hand, we need infinite depth for $x_{i}$ with $x_{i}=(X),$
as, otherwise,
we could fill in (classical) truth values from below.

Conversely, such an infinite sequence can be used to construct a Yablo
graph.
\subsubsection{
Back to Nested Diamonds
}

In our discussion of nested diamonds, we tried to do too much:
Consider, e.g., Diagram 
\ref{Diagram Rhombus-Basic} (page 
\pageref{Diagram Rhombus-Basic}).
We tried to append recursively new diamonds at both sides of the diamond,
e.g. at $x_{1,1}$ and $x_{1,2},$ but as
shown above, one recursive construction would have been sufficient.
We were blinded by the symmetry, instead of looking at what is needed.

We now modify Yablo's construction with a suitable diamond:

Instead of $x_{i} \xcP x_{i+1} \xcP x_{i+2}$ and $x_{i} \xcP x_{i+2},$ or,
more generally,
$x_{i} \xcP x_{i' } \xcP x_{j}$ and $x_{i} \xcP x_{j},$ $j \xcg i+2,$ we
 \xEh
 \xDH introduce new points $x_{i,j}$ and $y_{i,j}$

 \xDH set
$x_{i}= \xCN x_{i+1} \xcu \xcU \xCN x_{i,j}:j>i+1$

 \xDH set
$x_{i,j}=x_{j} \xco (y_{i,j} \xcu \xCN y_{i,j})=x_{j}$

 \xDH
Thus, we have

$ \xCN x_{i}=x_{i+1} \xco \xcO x_{i,j}:j>i+1$

$ \xCN x_{i,j}$ $=$ $ \xCN x_{j}$ $ \xcu $ $TRUE= \xCN x_{j}$

$x_{i}= \xCN x_{i+1} \xcu \xcU \xCN x_{j}:j>i+1$

$ \xCN x_{i}=x_{i+1} \xco \xcO x_{j}:j>i+1$

as in Yablo's construction.

 \xEj

See Diagram \ref{Diagram DD-Essential} (page \pageref{Diagram DD-Essential}).

\clearpage

\begin{diagram}

Diamond Essentials, Recursion on the Left, Simple Contradiction on the
Right

\label{Diagram DD-Essential}
\index{Diagram DD-Essential}

\unitlength0.5mm
\begin{picture}(280,220)(0,0)

\put(20,205){Essentials of Double Diamond}

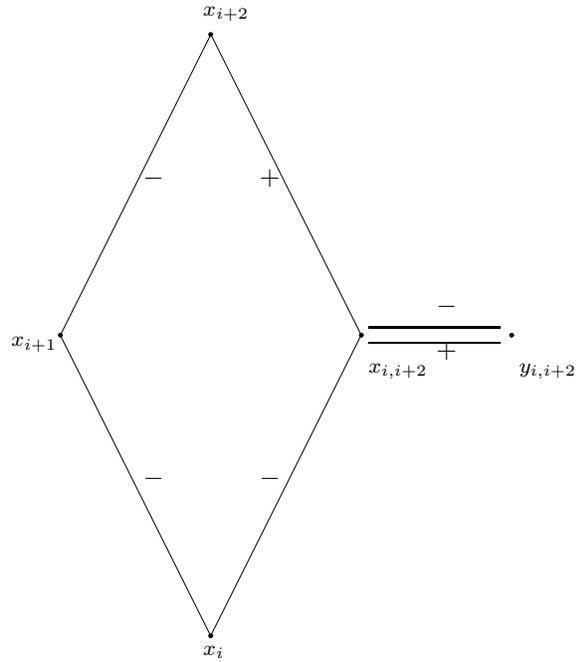
\put(120,20){\circle*{1}}
\put(118,15){\footnotesize $x_i$}
\put(20,5){Diagonal lines point upwards, the others to the right.}

\put(120,20){\line(-1,2){40}}
\put(120,20){\line(1,2){40}}

\put(102,60){$-$}
\put(133,60){$-$}

\put(80,100){\circle*{1}}
\put(160,100){\circle*{1}}
\put(200,100){\circle*{1}}

\put(67,97){\footnotesize $x_{i+1}$}
\put(162,90){\footnotesize $x_{i,i+2}$}
\put(202,90){\footnotesize $y_{i,i+2}$}

\put(162,102){\line(1,0){35}}
\put(162,98){\line(1,0){35}}

\put(80,100){\line(1,2){40}}
\put(160,100){\line(-1,2){40}}

\put(102,140){$-$}
\put(133,140){$+$}

\put(180,106){$-$}
\put(180,94){$+$}

\put(120,180){\circle*{1}}
\put(118,185){\footnotesize $x_{i+2}$}

\end{picture}

\end{diagram}

\vspace{10mm}

\clearpage

We now consider recursion on the other side of the diamond,
see Diagram \ref{Diagram Matrix-5} (page \pageref{Diagram Matrix-5}).

If $x+,$ then $x_{2}-,$ $x_{3}-,$ and $x_{4}-.$ Thus, also $x \xcP x_{2}
\xcP x_{2,2}$ has to be made
contradictory, here by $x \xcP x_{3} \xcp x_{2,2}.$ But, as $x_{3}-,$ we
cannot be sure
$x_{2,2}$ is chosen at $x_{3},$ it might just es well be that $x_{3,2}$ is
chosen,
unless this leads to a contradiction. So we introduce $x_{4}$ and
add $x \xcP x_{4} \xcp x_{3,2},$ to contradict $x \xcP x_{3} \xcP
x_{3,2}.$ But now we have the
same problem again with $x_{4}-.$ So this does not work.

\clearpage

\begin{diagram}

Recursion on the Right

\label{Diagram Matrix-5}
\index{Diagram Matrix-5}

\unitlength0.8mm
\begin{picture}(280,220)(0,0)


\put(40,50){\circle*{1}}
\put(38,45){\footnotesize $x$}

\put(20,70){\circle*{1}}
\put(60,70){\circle*{1}}
\put(100,70){\circle*{1}}
\put(140,70){\circle*{1}}

\put(18,65){\footnotesize $x_1$}
\put(58,65){\footnotesize $x_2$}
\put(98,65){\footnotesize $x_3$}
\put(138,65){\footnotesize $x_4$}

\put(40,50){\line(-1,1){20}}
\put(40,50){\line(1,1){20}}
\put(40,50){\line(3,1){60}}
\put(40,50){\line(5,1){100}}

\put(40,90){\circle*{1}}
\put(70,90){\circle*{1}}
\put(100,90){\circle*{1}}
\put(110,90){\circle*{1}}
\put(140,90){\circle*{1}}
\put(150,90){\circle*{1}}
\put(180,90){\circle*{1}}

\put(38,85){\footnotesize $y$}

\put(68,85){\footnotesize $x_{2,1}$}
\put(98,85){\footnotesize $x_{2,2}$}
\put(108,85){\footnotesize $x_{3,1}$}
\put(138,85){\footnotesize $x_{3,2}$}
\put(148,85){\footnotesize $x_{4,1}$}
\put(178,85){\footnotesize $x_{4,2}$}

\put(20,70){\line(1,1){20}}

\put(60,70){\line(-1,1){20}}
\put(60,70){\line(1,2){10}}
\put(60,70){\line(2,1){40}}

\put(100,70){\line(0,1){20}}
\put(100,70){\line(1,2){10}}
\put(100,70){\line(2,1){40}}

\put(140,70){\line(0,1){20}}
\put(140,70){\line(1,2){10}}
\put(140,70){\line(2,1){40}}

\put(85,110){\circle*{1}}
\put(125,110){\circle*{1}}
\put(165,110){\circle*{1}}

\put(83,105){\footnotesize $y_{2}$}
\put(123,105){\footnotesize $y_{3}$}
\put(163,105){\footnotesize $y_{4}$}

\put(70,90){\line(3,4){15}}
\put(100,90){\line(-3,4){15}}

\put(110,90){\line(3,4){15}}
\put(140,90){\line(-3,4){15}}

\put(150,90){\line(3,4){15}}
\put(180,90){\line(-3,4){15}}

\put(53,80){$+$}
\put(95,80){$+$}
\put(135,80){$+$}
\put(95,100){$+$}
\put(135,100){$+$}
\put(175,100){$+$}

\put(20,30){All lines stand for upward pointing arrows, unmarked ones}
\put(20,22){for negative arrows.}

\end{picture}

\end{diagram}

\vspace{10mm}

\clearpage

\subsection{
Infinite Branching
}

\label{Section Branch}

We argue now that we have infinite branching, perhaps not immediately,
at $x_{0}$ as in Construction 
\ref{Construction Details} (page 
\pageref{Construction Details}), but after finitely
many steps.

 \xEh
 \xDH
We analysed the possible contradictions in
Section \ref{Section VariousTypes} (page \pageref{Section VariousTypes}).

 \xDH
We saw in Remark 
\ref{Remark Rhombus-Basic} (page 
\pageref{Remark Rhombus-Basic})  that doubly (on both sides)
recursive
Diamonds do not work, and in
Section \ref{Section How-Much} (page \pageref{Section How-Much})  that simply
recursive diamonds are either trivial modifications of the Yablo Cell,
or do not work.
So, wlog., we may assume that all contradictions have the form of Yablo
Cells (with maybe paths instead of single arrows).

This applies to procrastinated contradictions, see
Diagram 
\ref{Diagram Escape} (page 
\pageref{Diagram Escape}), too, resulting in ``Yablo Cells with a
stem''.

 \xDH
Diagram \ref{Diagram X1} (page \pageref{Diagram X1})  and
Diagram \ref{Diagram X2} (page \pageref{Diagram X2})
show that all arrows (or paths) have to be negative.

 \xDH
Construction \ref{Construction Details} (page \pageref{Construction Details})
shows that the construction is recursive, and we have infinite
width.

 \xEj
\clearpage

$ \xCO $

\end{document}